\documentclass[10pt]{amsart}
\usepackage{amsmath,amssymb,amsfonts} 
\usepackage{graphics}                 
\usepackage[square, numbers]{natbib}
\usepackage[french,english]{babel}
\usepackage{tikz} 
\usepackage[utf8]{inputenc}
\usepackage[all]{xy}	
\usepackage[T1]{fontenc}
\usepackage{pifont}

\usepackage[margin=1.2in]{geometry}

\newcommand{\Mod}{\mathfrak{M}}

\newcommand{\U}{\mathcal{U}}

\newcommand{\C}{\mathbb{C}}

\newcommand{\N}{\mathbb{N}}
\newcommand{\Z}{\mathbb{Z}}
\newcommand{\R}{\mathbb{R}}

\newcommand{\Spec}{\rm{Spec}}

\renewcommand{\L}{\mathcal{L}}

\theoremstyle{definition}
\newtheorem{defn}{Définition}[subsection]
\newtheorem{rem}[defn]{Remarque}

\newtheorem{cons}[defn]{Construction}
\newtheorem{exam}[defn]{Exemple}

\newtheorem{nota}[defn]{Notation}

\theoremstyle{plain}
\newtheorem{Prop}[defn]{Proposition}
\newtheorem{Thm}[defn]{Théorème}
\newtheorem{cor}[defn]{Corollaire}
\newtheorem{lem}[defn]{Lemme}

\newtheorem*{asser}{Assertion}

\newtheorem*{ThmInt}{Théorème}
\newtheorem*{lemInt}{Lemme}

\def\Ind{\setbox0=\hbox{$x$}\kern\wd0\hbox to 0pt{\hss$\mid$\hss}
\lower.9\ht0\hbox to 0pt{\hss$\smile$\hss}\kern\wd0}
\def\Notind{\setbox0=\hbox{$x$}\kern\wd0\hbox to 0pt{\mathchardef
\nn=12854\hss$\nn$\kern1.4\wd0\hss}\hbox to
0pt{\hss$\mid$\hss}\lower.9\ht0 \hbox to
0pt{\hss$\smile$\hss}\kern\wd0}
\def\ind{\mathop{\mathpalette\Ind{}}}
\def\nind{\mathop{\mathpalette\Notind{}}}

\renewcommand{\Spec}{\mathrm{Spec}}

\newtheoremstyle{exampstyle}
{3pt} 
{3pt} 
{\itshape} 
{} 
{\bfseries} 
{.} 
{.5em} 
{} 

\theoremstyle{exampstyle}\newtheorem{thmx}{Théorème}

\begin{document}

\title[Corps différentiels et flots géodésiques I]{Corps différentiels et flots géodésiques I : \\  \footnotesize\mdseries Orthogonalité aux constantes pour les équations différentielles autonomes}
\author{Rémi Jaoui}
\thanks{Recherche soutenue en partie par le contrat ANR-13-BS01-0006}
\address{Rémi Jaoui, Pure Mathematics Department, University of Waterloo, Ontario, Canada}
\email{rjaoui@uwaterloo.ca}
\date\today
\selectlanguage{french}
\begin{abstract}
L'orthogonalité aux constantes est une propriété issue de l'étude modèle-théorique des équations différentielles algébriques et qui traduit des propriétés d'indépendance algébrique remarquables pour ses solutions.

Dans cet article, on étudie la propriété d'orthogonalité aux constantes dans un langage algebro-différentiel pour les équations différentielles autonomes ainsi que des méthodes effectives pour établir cette propriété. Le résultat principal est un critère d'orthogonalité aux constantes (et sa version en famille) pour les $D$-variétés réelles absolument irréductibles $(X,v)$ s'appuyant sur la dynamique du flot réel associé $(M,\phi)$. Plus précisément, on montre que s'il existe des parties compactes $K$ de $M$, Zariski-dense dans $X$ telle que la restriction du flot à $K$ est topologiquement faiblement mélangeante, alors le type générique de $(X,v)$ est orthogonal aux constantes.

Ce critère sera appliqué dans \cite{moi2} à l'étude modèle-théorique du flot géodésique sur les variétés riemanniennes compactes à courbure strictement négative, présentées algébriquement. 
\end{abstract}

\maketitle

L'étude des propriétés d'indépendance algébrique des solutions d'une équation différentielle algébrique est un thème majeur de la théorie des équations différentielles, qui se situe au coeur des travaux de nombreux illustres mathematiciens --- notamment Newton, Leibniz, Picard, Painlevé et Poincaré. En particulier, à la fin du XIX$^{eme}$ siecle, Drach étudie dans sa thèse \cite{Drach-thesis}, l'existence d'un cadre algébrique (désincarné de ses réalisations analytiques) où étudier les propriétés d'intégrabilité algébrique de ces équations différentielles, à la manière de Picard et Vessiot pour les équations différentielles linéaires (voir \cite{Archibald} pour un aperçu historique).

Néanmoins, les problèmes de formalisme rencontrés par Drach  --- comparables à ceux rencontrés en géométrie algébrique à la même époque --- ont poussé les mathématiciens du début du XX$^{\text{ème}}$ siècle à privilégier un cadre analytique pour étudier les équations différentielles algébriques alors même que ces dernières vérifient (très vraisemblablement) des propriétés structurelles plus fortes que leurs analogues analytiques.

A la fin des années 1970, les développements conjoints de la théorie des modèles des théories stables et de l'algèbre différentielle ont permis de développer un cadre géométrique noethérien, propice à l'étude des équations différentielles algébriques. Une caractéristique essentielle  est le recours à des corps différentiels ``universels'' appelés \textit{corps differentiellement clos}.

L'étude de ces ``compagnons de route de la théorie des modèles'' (voir \cite{Poi2}) est une illustration remarquable de l'efficacité des méthodes de la théorie géométrique de la stabilité pour étudier certaines géometries noetheriennes, innaccessibles auparavant. Un résultat essentiel est le théorème de Hrushovski-Sokolovic \cite{Sok} --- une incarnation du théorème de trichotomie de Hrushovski et Zilber pour les corps differentiellement clos --- qui décrit \textit{les équations différentielles minimales}, dont la résolution ne peut être ramenée à la résolution successive d'équations différentielles ``plus simples''.\\

L'objectif de cet article est de présenter quelques conséquences du théorème de Hrushovski-Sokolovic pour les équations différentielles autonomes (dans un langage géométrique) ainsi que d'établir des critères effectifs pour étudier les propriétés mises en jeu par ce théorème  --- \textit{orthogonalité aux constantes et désintégration}. En s'appuyant sur des résultats d'Anosov \cite{Ano} concernant la dynamique des champs de vecteurs hyperboliques, ces résultats seront appliqués dans \cite{moi2} et \cite{moi3} à l'étude de propriétés d'indépendance algébrique des solutions géodésiques  d'une variété riemannienne compacte (présentée algébriquement) à courbure strictement négative.
\\

\textbf{Équations différentielles algébriques autonomes.} Dans ce texte, on travaille uniquement au dessus de corps $k$ \textit{de caractéristique $0$} et on considère des \textit{équations différentielles algébriques autonomes} à paramètres dans $k$. Une telle equation peut être décrite ``sous forme explicite'' comme une paire $(X,v)$ où \textit{$X$ est une variété algébrique au dessus d'un corps $k$ de caractéristique $0$} --- que l'on pourra toujours supposée lisse --- \textit{munie d'un champ de vecteurs $v$}.

Les relations algébriques entre $n$ solutions de $(X,v)$ sont représentées par certaines \textit{sous-variétés fermées} de $(X,v)^n$ dites \textit{invariantes}.

Ici, une sous-variété fermée (et plus généralement tout sous-schéma fermé) de $(X,v)$ est \textit{invariant} si le faisceau d'idéaux associé $\mathcal I \subset (\mathcal O_{X}, \delta_v)$ est stable par la dérivation $\delta_v$ induite par le champ de vecteurs $v$ sur $X$. Géométriquement, une sous-variété fermée $Z$ de $(X,v)$ est invariante lorsque la sous-variété localement fermée $Z_{reg}$ de $X$ (constituée des points réguliers de $Z$) est tangente en tout point au champ de vecteurs $v$.

Les sous-variétés fermées invariantes de $(X,v)^n$ sont définies, tout simplement, comme les sous-variétés fermées de $X^n$ invariante pour le champ de vecteurs produit $v \times \cdots \times v$ induit par $v$ sur $X^n$. En notant $\mathcal I_n(X,v)$ la collection des sous-variétés fermées irréductibles invariantes de $(X,v)^n$, la suite des  $\mathcal I_n(X,v)$ lorsque $n$ parcourt $\mathbb{N}^\ast$ définit alors un langage $\mathcal L$ (au sens de la théorie des modèles).

L'ensemble des solutions analytiques de l'équation différentielle $(X,v)$ est naturellement munie d'une $\mathcal L$-structure. Comme il est d'usage en théorie des modèles, il convient alors de compléter cette structure en \textit{une structure existentiellement close} notée $(X,v)^{(\U,\delta)}$  en adjoignant toutes les solutions de cette équation différentielle dans un corps différentiel universel $(\U,\delta)$.

La $\mathcal L$-structure $(X,v)^{(\U,\delta_\U)}$, produite par la construction précédente, jouit de propriétés particulièrement agréables : c'est une \textit{structure $\omega$-stable de rang fini} (borné par la dimension de $X$), interprétable dans un corps differentiellement clos, qui admet \textit{l'élimination des quantificateurs dans le langage $\mathcal L$}. Les méthodes de la théorie géométrique de la stabilité appliquées à  la structure $(X,v)^{(\U,\delta)}$  permettent alors d'obtenir:

(a) Une \textit{classification des équations différentielles  minimales}. Intuitivement ce sont celles dont la résolution ne peut être ramenée à la résolution d'équations différentielles plus simples.

(b) Un dévissage des équations différentielles de rang fini à partir d'équations différentielles minimales selon un procédé appelé \textit{analyse semi-minimale des types}. 

\hspace{0.01cm}

Dans le paragraphe qui suit, on décrit plus précisément ces deux procédés (en commençant par le second) dans le langage de \textit{la théorie géométrique de la stabilité}.\\

\textbf{Théorie géométrique de la stabilité. }  Plutôt que la structure associée à l'équation différentielle $(X,v)$, les résultats de décomposition dont il est question ici, concernent en réalité \textit{les types}  --- qui sont les ultrafiltres sur l'algèbre booléenne des ensembles définissables --- de cette structure. Lorsque la variété $X$ est irréductible, il existe un type privilégié parmi ces derniers (par élimination des quantificateurs), appelé \textit{type générique de $(X,v)$} dont les réalisations sont les solutions génériques de l'équation différentielle $(X,v)$, c'est-à-dire, celles qui sont Zariski-dense dans $X$. \\

Commençons par décrire les deux situations extrêmes qui peuvent se produire pour les relations entre un type de rang fini $q$ et un type minimal $p$ dans une théorie stable $T$: soit il est possible de reconstruire (après une possible extension des paramètres) une réalisation du type $q$ à  partir de réalisations de $p$ (et de ses conjugués). On dit alors que le type \textit{$q$ est interne à $p$} (et ses conjugués). A l'inverse, le \textit{type $q$ est orthogonal à $p$} lorsque $p$ et $q$ sont sans relation  dans la théorie $T$ et le demeurent après toute extension des paramètres, ou autrement dit après tout changement de base.

En général, on peut aussi observer des situations intermédiaires, où le type $q$ est non-orthogonal au type minimal $p$ sans pour autant lui  être interne. Dans ce cas, il existe alors un \textit{facteur non-trivial}\footnote{Le type $q_0 \in S(\emptyset)$ est un facteur non-trivial de $q \in S(\emptyset)$ s'il existe des réalisations $a \models q$ et $b \models q_0$ telles que $b \in \mathrm{dcl}(a) \setminus \mathrm{acl}(\emptyset)$.} $q_0$ de $q$ qui est interne à $p$ (et ses conjugués). Quitte  à remplacer le type $q$ par un facteur non-trivial $q_0$, on peut donc toujours se ramener à une des deux situations extrêmes exposées ci-dessus. 

De plus, on peut alors itérer ce processus pour obtenir un dévissage complet du type $q$ en types semi-minimaux (c'est-à-dire interne à un type minimal et ses conjugués) appelé \textit{analyse semi-minimale de $q$}. Appliqué aux équations différentielles algébriques, ce dévissage permet ainsi de décomposer toute equation différentielle algébrique $(X,v)$, au voisinage de son point générique, en une suite de facteurs semi-minimaux. Décrivons maintenant la classification des équations différentielles minimales, mentionnée ci-dessus.

L'étude des types minimaux  --- pour toute théorie stable et donc dans un cadre bien plus général que celui des équations différentielles --- a conduit Zilber  à les classifier en trois différentes classes   selon la nature de la notion de dimension qui leur est naturellement associée. De façon informelle, la notion de dimension se comporte ou bien comme le cardinal dans un ensemble infini sans structure, ou bien comme la dimension linéaire dans un espace vectoriel, ou bien comme le degré de transcendance dans un corps algébriquement clos. On dit alors respectivement que \textit{le type minimal est désintégré, localement modulaire (non désintégré) ou non-localement modulaire}.

Pour les équations différentielles algébriques minimales, le cas désintégré correspond à celles dont les solutions (génériques) possèdent \textit{les plus fortes propriétés d'indépendance algébrique} (voir le paragraphe suivant). Le théorème de Hrushovski-Sokolovic (voir \cite{Sok}) est une classification complète des équations différentielles minimales (pas nécessairement autonomes) non-désintégrées.

{\ThmInt[Hrushovski-Sokolovic] Soit $p \in S(A)$ un type minimal non-désintégré de la théorie des corps differentiellement clos. Alors: 
\begin{itemize}
\item[(i)] Ou bien le type $p$ est non-localement modulaire. Dans ce cas, le type $p$ est non-orthogonal au type générique du corps des constantes.

\item[(ii)] Ou bien le type $p$ est localement modulaire. Il existe alors une variété abélienne simple  $A$, qui ne descend pas au corps des constantes, tel que $p$ est non-orthogonal au type générique du ``noyau de Manin'' $A^\sharp$ associé  à $A$. 
\end{itemize}}

Dans \cite{Itai}, Hrushovski et Itai étudient les conséquences de ce théorème pour les équations différentielles autonomes $(X,v)$ de \textit{dimension $1$ où l'hypothèse de minimalité est automatiquement vérifiée}. En dimension supérieure, il est très difficile de formuler une description géométrique (raisonnable en pratique) de cette propriété (voir par exemple \citep[Chapitre 3, Théorème 1.3.12]{jaoui} pour une formulation géométrique). \\

\textbf{Désintégration et orthogonalité aux constantes.} Dans ce paragraphe, nous décrivons en termes géométriques quelques conséquences du théorème de Hrushovski-Sokolovic pour des équations différentielles autonomes $(X,v)$ (de dimension supérieure et donc sans hypothèse de minimalité \textit{a priori} pour leur type générique) à travers une analyse semi-minimale de leur type générique.  \\

La propriété d'\textit{orthogonalité aux constantes} pour (le type générique d') une équation différentielle $(X,v)$ peut être formulée à  l'aide de la notion d'intégrale première rationnelle dans $(X,v)$ et ses puissances. Ici, \textit{une intégrale première rationnelle} de $(X,v)$ est simplement une fonction rationnelle $f \in k(X)$ sur $X$  constante le long des orbites de $v$, c'est-a-dire vérifiant $v(f) = 0$.  Le type générique de $(X,v)$ est alors orthogonal aux constantes si et seulement si:
$$(O): \text{ Pour tout entier } n \geq 0, \text{ } (X,v)^n \text{ n'admet pas } \text{d'intégrale première rationnelle non constante}.$$

La propriété de \textit{désintégration} (au point générique) de $(X,v)$ concerne, quant à elle, les relations algébriques entre un nombre arbitraire de solutions génériques de $(X,v)$. Intuitivement, elle exprime que l'ensemble des relations algébriques entre $n$ solutions de l'équation différentielle $(X,v)$ est aussi petit que possible dès que $n \geq 3$.

Formellement, notons $\mathcal I_n^{gen}(X,v)$ l'ensemble des sous-variétés fermées irréductibles invariantes \textit{génériques} (c'est-à-dire celles se projetant génériquement sur chacun des facteurs) de $(X,v)^n$. L'équation différentielle $(X,v)$ est \textit{génériquement désintégrée} si et seulement si:   

$$(D):\begin{cases} \text{ Pour tout } n \geq 3 \text{ et tout } Z \in \mathcal I_n^{gen}(X,v), \text{ il existe } (Z_{i,j})_{i \neq j \leq n} \in \mathcal I_2^{gen}(X,v) \text{  pour tout } i \neq j, \\ \hspace{1.8cm} \text{ telles que } Z \text{ est une composante irreductible de } \bigcap_{1 \leq i \neq j \leq n} \pi_{i,j}^{-1}(Z_{i,j}).
\end{cases}$$ 

Pour toute équation différentielle algébrique autonome, \textit{la propriété (D) implique la propriété (O)}. En apparence, la propriété (D) --- qui concerne des sous-variétés fermées invariantes de $(X,v)^n$ de dimension arbitraire --- semble bien plus restrictive que la propriété (O). Néanmoins, nous montrons que la propriété $(O)$ est suffisante pour garantir l'existence d'un facteur rationnel (non-trivial) génériquement désintégré.

{\thmx\label{intro1} Soient $k$ un corps de caractéristique $0$ et $(X,v)$ une variété absolument irréductible au dessus de $k$ munie d'un champ de vecteurs $v$.

Si le type générique de $(X,v)$ est orthogonal aux constantes alors il existe une variété absolument irréductible $Y$ au dessus de $k$ de dimension $> 0$, un champ de vecteurs $w$ sur $Y$ et un morphisme rationnel dominant 
$$ \pi : (X,v) \dashrightarrow (Y,w)$$
tels que $(Y,w)$ est génériquement désintégrée.}

Le morphisme $\pi$ apparaissant dans le théorème \ref{intro1} est un morphisme de la catégorie des $D$-variétés (voir Appendice A), c'est-à-dire que le morphisme rationnel $\pi: X \dashrightarrow Y$ vérifie la règle de compatibilité évidente aux champs de vecteurs: $d\phi(v) = w$. Le théorème \ref{intro1} est obtenu dans la première partie de ce texte, en appliquant le théorème de Hrushovski-Sokolovic à un facteur semi-minimal de l'équation différentielle $(X,v)$. \\

Pour résumer, la propriété d'\textit{orthogonalité aux constantes} pour une équation différentielle algébrique concerne donc aussi bien \textit{l'absence d'intégrales premières rationnelles pour l'équation différentielle sous étude et ses puissances} que \textit{l'existence d'un facteur rationnel génériquement désintégré}. \\

\textbf{Un critère dynamique d'orthogonalité aux constantes.} Paradoxalement, on ne connait que très peu d'exemples explicites d'équations différentielles algébriques (en dimension $> 1$) ou cette propriété de désintégration a été établie. Mentionnons deux exemples récents: en dimension $2$, les équations différentielles (non-autonomes) de Painlevé à  parametres génériques étudiées par Pillay et Nagloo dans \cite{Nag} et en dimension $3$, l'équation différentielle satisfaite par la fonction $j$ et ses $Gl_2(\mathbb{C})$-conjugués, étudiée par Freitag et Scanlon dans \cite{Frei}.

Dans ces deux cas, leurs démonstrations s'appuient sur des propriétés très spécifiques des solutions des équations différentielles considérées. Nous démontrons dans ce texte, \textit{un critère dynamique d'orthogonalité aux constantes pour les équations différentielles réelles} au caractère bien plus général, permettant notamment de construire des familles non-limitées\footnote{c'est-à-dire ne provenant pas toutes d'une même famille ne dépendant d'un nombre fini de paramètres.} d'équations différentielles algébriques dont le type générique est orthogonal aux constantes (voir aussi le paragraphe suivant concernant d'autres applications de ce critère).

Notre critère d'orthogonalité aux constantes concernent certaines équations différentielles algébriques $(X,v)$  \textit{définies sur le corps $\R$ des nombres réels} (munie de la dérivation triviale). Sous des hypothèses naturelles de lissité, on peut considérer l'espace analytique réel $X(\R)^{an}$ \textit{des conditions initiales réelles} ainsi que le flot analytique réel $\phi$ du champs de vecteurs $v$  (défini sur un voisinage connexe $U$ de $X(\R)^{an} \times \lbrace 0 \rbrace$ dans $X(\R)^{an} \times \R$) 
$$ \phi : U \subset X(\R)^{an} \times \R \longrightarrow \R$$
décrivant pour tout $(t,x_0) \in U$, la position $\phi(x_0,t)$ de la condition initiale $x_0$ au temps $t$.

En général, le flot analytique réel $\phi$ ne peut être étendu à $X(\R)^{an} \times \R$ tout entier --- on parle alors d'explosion en temps fini des solutions --- et il est nécessaire de se restreindre à  un voisinage connexe de $X(\R)^{an} \times \lbrace 0 \rbrace$ dans $X(\R)^{an} \times \R$.

Lorsque c'est le cas et que le flot analytique réel $\phi$ peut être étendu  à  $X(\R)^{an} \times \R$ tout entier, on dit alors que le flot $\phi$ est \textit{complet}. Dans ce cas, on l'identifie à une action continue du groupe additif $(\R, +)$ par automorphismes analytiques sur $X(\R)^{an}$ ou en d'autres termes à un \textit{$(\R, +)$-système dynamique analytique réel}. Cette condition est automatique lorsque $X(\R)^{an}$ est compact ou plus généralement lorsque l'on se restreint à travailler sur un à compact $K \subset X(\R)^{an}$ invariant.

Le critère d'orthogonalité aux constantes, que nous démontrons ici, est extension modèle-théorique des arguments de  non-intégrabilité \textit{à la Poincaré}, reposant sur la complexité de la dynamique topologique réelle des équations différentielles considérées.

{\thmx\label{intro2} Soient $X$ une variété absolument irréductible sur $\R$ et $v$ un champ de vecteurs rationnel sur $X$. On note $(M, \phi)$ le flot régulier réel de $(X,v_X)$.
Supposons qu'il existe une partie compacte $K $ de $M$, Zariski-dense dans $X$ et invariante par le flot $\phi$.

Si $(K,(\phi_{t |K})_{t \in \R})$ est faiblement topologiquement mélangeant alors le type générique de $(X,v)$ est orthogonal aux constantes. \vspace{3pt}}

Par définition, $M = X(\R) \setminus (\mathrm{Sing}(X)\cup \mathrm{Sing}(v))$ est ici la variété analytique réelle vérifiant les hypothèses de lissité nécessaires à l'intégration analytique du champ de vecteurs $v$. 
La restriction de ce flot à toute partie compacte invariante $K$ de $M$ définit alors un système dynamique réel. Ce système dynamique est dit \textit{faiblement topologiquement mélangeant} si tous ses produits sont topologiquement transitifs (voir la partie 3.2 pour des formulations équivalentes de cette condition).

Le théorème \ref{intro2} est prouvé dans la quatrième partie de ce texte  à partir d'une étude consciencieuse (dans la troisième partie de ce texte) des compatibilités entre les notions d'invariance dans la catégorie ``algébrique'' des $D$-variétés et la catégorie ``dynamique'' des systèmes dynamiques analytiques. \\

\textbf{Application aux équations de la mécanique classique.} Dans ce paragraphe, nous décrivons brièvement les applications du théorème \ref{intro2} aux équations géodesiques à courbure negative qui sont l'objet de \cite{moi2} et \cite{moi3}. \\
 
Nous considérons des  équations différentielles modélisant des \textit{systèmes autonomes conservatifs de la mécanique classique}, décrits de la façon suivante: Si $(M,g)$ est une variété riemannienne présentée algébriquement et $V$ une fonction algébrique sur $M$ représentant l'énergie potentielle, ce sont les équations hamiltoniennes associé au hamiltonien: 

$$ H(x,p) = \frac 1 2 g_x(p,p) + V(x)$$

(voir \cite{moi2} pour une description du formalisme hamiltonien dans ce cadre algébrique).

Pour un tel système hamiltonien, les relations algébriques entre $n$ solutions (resp. génériques) d'un même niveau d'énergie $H = E_0$ sont décrites par l'ensemble $\mathcal I_n(H;E_0)$ (resp. $\mathcal I^{gen}_n(H,E_0)$)  des sous-variétés fermées $Z \subset T_{M/\R}^n$  qui s'écrivent 
$$Z = \overline{ \lbrace (\gamma_1(t), \dot \gamma_1(t), \ldots , \gamma_n(t),\dot \gamma_n(t)) \text{ | } t \in \mathbb{D} \rbrace} .$$ 
comme la clôture de Zariski de $n$ solutions analytiques (resp. génériques) $\gamma_1, \ldots , \gamma_n$ (définies sur un disque complexe $\mathbb{D}$), d'énergie $E_0$, de l'équation hamiltonienne sous étude.

Dans \cite{moi2}, nous étudions le cas des équations géodésiques d'une variété riemannienne compacte à courbure strictement négative --- c'est-à-dire lorsque $V = 0$ et la variété riemannienne $(M,g)$ est compacte a courbure strictement négative.
Dans ce cas, les résultats d'Anosov \cite{Ano} mettent en évidence des proprietés hyperboliques globales pour la dynamique du flot analytique réél --- regroupées sur le terme de \textit{flot d'Anosov} --- sur les niveaux d'énergie non nul du hamiltonien $H$. En utilisant conjointement ces résultats d'Anosov, le théorème \ref{intro1} et le théorème \ref{intro2}, nous obtenons dans \cite{moi2}, \textit{l'existence d'un facteur désintégré pour le système de relations algébriques:} 
$$\mathcal L^{gen}(H,E_0) = \lbrace \mathcal I^{gen}_n(H,E_0) \text{ } | \text{  } n \geq 1 \rbrace \text{ lorsque } E_0 \neq 0.$$

Plus récemment et de façon complémentaire, j'ai étudié dans \cite{moi3} les facteurs rationnels d'une équation différentielle réelle $(X,v)$  dont le flot analytique est un flot d'Anosov mélangeant, lorsque la dimension de $X$ est $3$. Ces résultats impliquent que \textit{pour une surface riemannienne $(M,g)$ (présentée algébriquement) compacte à courbure strictement négative,} \textit{le système de relations algébriques $\mathcal L^{gen}(H,E_0)$ de l'équation différentielle géodésique de $(M,g)$ sur un niveau d'énergie non nul est lui-même désintégré}. \\

%
%

\textbf{Familles d'équations différentielles algébriques.} Outre les équations différentielles ``concrètes'' évoquées ci-dessus, nous étudions aussi les propriétés d'orthogonalité aux constantes en famille pour les familles lisses d'équations différentielles algébriques.

En termes géométriques, une famille lisse  d'équations différentielles algébriques $f: (\mathcal X,v) \longrightarrow (S,0)$ (à paramètres constants) est décrite par la donnée d'une famille lisse $f: \mathcal X \longrightarrow S$ de variétés algébriques paramétrée par $S$ au dessus d'un corps $k$ de caractéristique $0$ et d'un champ de vecteurs $v$ sur $\mathcal X$ tangent aux fibres de $f$.

Le champ de vecteurs $v$ se restreint ainsi le long de chaque fibre et on pense au morphisme $f$ comme à la famille $ \lbrace (X,v)_s \text{ | } s \in S \rbrace$ de ses fibres. Lorsque toutes les fibres de $f$ sont absolument irréductibles (et donc que le type générique de chacune des fibres est stationnaire), on peut alors considérer 
$$ S^{\perp 0} = \lbrace s \in S(\overline{l}) \text{ | }  
 \text{ le type generique de } (X,v)_s \text{ est orthogonal aux constantes} \rbrace.$$  
où $\overline{l}$ est une extension algébriquement close fixée saturée de $k$.

Les exemples de dimension $1$ de Hrushovski et Itai (voir \cite{Itai}) montrent qu'en général, le sous-ensemble $S^{\perp 0} \subset S(\overline{l})$ n'est ni définissable, ni $\infty$-définissable. En revanche, des considérations très générale de théorie de la stabilité  permettent de présenter $S^{\perp 0}$ comme un sous-ensemble $\mathrm{Aut}(\overline{l}/k)$-invariante de $S(\overline{l})$ --- ou autrement dit, comme l'image inverse d'un sous-ensemble du schéma $S$.

Pour les familles d'équations différentielles algébriques complexes, cette propriété admet la conséquence suivante:

{\lemInt Soit $S$ une variété irréductible algébrique complexe et $f: (\mathcal X,v) \longrightarrow (S,0)$ une famille d'équations différentielles algébriques paramétrée par $S$ dont les fibres sont toutes absolument irréductibles.

Si la fibre générique de $f$ est orthogonale aux constantes alors, il existe un ensemble dénombrable de sous-variétés algébriques fermées propres $\lbrace Z_n \text{ | } n \in \mathbb{N}\rbrace$  de $S$ tel que:
$$ S(\mathbb{C}) \setminus \bigcup_{i \in \mathbb{N}} Z_i (\mathbb{C}) \subset  S^{\perp 0}(\C).$$}   

Plus généralement, cet énoncé est valide si l'on remplace les nombres complexes par tout corps algébriquement clos. Dans le cas complexe qui est formulé ici, le théorème de Baire permet de garantir la non-vacuité de $S(\mathbb{C}) \setminus \bigcup_{i \in \mathbb{N}} Z_i (\mathbb{C})$. Plus précisément, si la fibre générique de $f$ est orthogonale aux constantes alors pour presque tout $s \in S(\mathbb{C})$ au sens de Baire,  l'équation différentielle $(X,v)_s $ est orthogonale aux constantes.

Le lemme précédent sera appliqué conjointement avec le théorème de spécialisation suivant concernant la propriété de non-orthogonalité aux constantes pour les familles lisses d'équations différentielles algébriques.

{\thmx\label{intro3} Soient $S$ une variété algébrique lisse et irréductible au dessus d'un corps $k$ et  $f : (\mathcal X,v) \longrightarrow (S,0)$ une famille lisse de $D$-variété à paramètres dans $S$.

On suppose que toutes les fibres de $f$ sont absolument irréductibles, on fixe $s \in S(k)$ et on dénote par $\eta$ le point générique de $S$. Si le type générique de  $(X,v)_\eta$ est non-orthogonal aux constantes alors le type générique de $(X,v)_s$ est non-orthogonal aux constantes. \vspace{0.2cm}}

Pour les familles d'équations différentielles complexes $f : (\mathcal X,v) \longrightarrow (S,0)$ vérifiant les hypothèses d'irréductibilité et de lissité du théorème \ref{intro3}, le théorème \ref{intro3} et le lemme ci-dessus impliquent donc que:
$$\text{ ou bien } S^{\perp 0} = \emptyset \text{, ou bien }  S^{\perp 0} = S(\mathbb{C}) \setminus \bigcup_{i \in \mathbb{N}} Z_i (\mathbb{C}) $$
pour un ensemble dénombrable de sous-variétés algébriques fermées propres $\lbrace Z_n \text{  }|  n \in \mathbb{N}\rbrace$  de $S$ (voir le corollaire \ref{specialisationgeometrique}).
%

La preuve du théorème \ref{intro3} repose sur l'équivalence entre la propriété d'orthogonalité aux constantes et sa formulation géométrique ( la propriété (O)) Sous les hypothèses de lissité, les énoncés de spécialisation correspondant pour les intégrales premières rationnelles sont prouvés dans la seconde partie cet article. 
\subsection*{Équations différentielles à paramètres très génériques} Pour conclure, nous formulons maintenant deux applications du théorème \ref{intro3} aux équations différentielles algébriques à paramètres très génériques.

Le premier résultat concerne les champs de vecteurs de degré $d$ sur $\mathbb{A}^n_\mathbb{C}$ avec $n \geq 1$.

{\thmx\label{intro4} Soit $d \geq 3$ et $n \geq 1$. Considérons un champs de vecteurs 
$$v(x_1,\ldots , x_n) = f_1(x_1,\ldots, x_n) \frac {d} {dx_1} + \cdots f_n(x_1,\ldots, x_n) \frac {d} {dx_n}.$$ 
sur l'espace affine complexe de dimension $n$, où $f_1,\ldots, f_n \in K[X_1,\ldots, X_n]_{\leq d}$ sont des polynômes de degré $\leq d$.

Si les coefficients de $f_1, \ldots, f_n$ sont $\mathbb{Q}$-algébriquement indépendants\footnote{Cela implique, en particulier, que les $f_i$ sont des polynômes de degré $d$ dont tous les coefficients sont distincts et non nuls.} alors le champ de vecteurs $v$ est orthogonal aux constantes.}

Lorsque $n = 1$, le theoreme \ref{intro4} est conséquence immédiate de résultats de Rosenlicht \cite{Rosenlicht}. En dimension supérieure, notre démonstration  (voir Section 2.4) s'appuie le cas de dimension $1$ et sur le théorème \ref{intro3}. \\

Pour conclure, formulons une seconde application du théorème \ref{intro3} qui est une variante du critère dynamique d'orthogonalité aux constantes pour la fibre générique d'une famille d'équations différentielles algébriques. Cet énoncé obtenu en utilisant conjointement le théorème \ref{intro2} et le théorème \ref{intro3}.

{\thmx\label{intro5} Soient $k$ un sous-corps des nombres réels, $S$ une variété algébrique lisse et irréductible au dessus de $k$ et  $f : (\mathcal X,v) \longrightarrow (S,0)$ une famille lisse de $D$-variétés absolument irréductibles à paramètres dans $S$. 

Supposons qu'il existe un point $p \in S(\R)$ et un compact $K \subset \mathcal X_p(\R)^{an}$ Zariski-dense dans $X$ et invariant par le flot $\phi$ du champ de vecteurs $v_{|X_p}$.

Si $(K,(\phi_{t |K})_{t \in \R})$ est faiblement topologiquement mélangeant alors il existe un ensemble dénombrable $\lbrace Z_i : i \in \mathbb{N} \rbrace$ de sous-variétés fermées algébriques strictes (sur le corps des nombres réels) $Z_i$ de $S$ tel que:

$$\forall s \in S(\mathbb C) \setminus \bigcup_{i \in \mathbb{N}} Z_i(\mathbb C)\text{, } (\mathcal X,v)_s \text{ est orthogonal aux constantes}.$$}

De façon similaire aux applications du théorème \ref{intro2} à l'étude des équations géodésiques sur une variété riemannienne compacte à  courbure strictement négative, le théorème \ref{intro4} peut être appliqué à l'étude \textit{des équations hamiltoniennes pour un potentiel $V$ à paramètres très génériques sur une variété riemannienne compacte à courbure strictement négative}. 

Plus précisément, considérons $(M,g)$ est une variété riemannienne compacte à courbure strictement négative (présentée algébriquement) et $\mathcal V \in \mathrm{H}^0(S \times M, \mathcal O_{S \times M})$ une famille algébrique de potentiels sur $M$, paramétrée par une variété algébrique réelle $S$. Supposons de plus que pour un point $s_0 \in S(\R)$, on ait $V(s_0,0) = 0$.

On obtient alors une famille (au sens précédent) $ f_\mathcal V:(T_{M/\R}, v_\mathcal V) \longrightarrow (S,0)$ d'équations différentielles hamiltoniennes paramétrée par $S$ dont la fibre au dessus de $s_0$ est l'équation différentielle géodésique de $(M,g)$. Si la famille $f_\mathcal V$ vérifie les hypothèses du théorème \ref{intro4} alors il existe un ensemble dénombrable $\lbrace Z_i : i \in \mathbb{N} \rbrace$ de sous-variétés fermées algébriques strictes  $Z_i$ de $S$ tel que \textit{pour tout $s \in S(\mathbb C) \setminus \bigcup_{i \in \mathbb{N}} Z_i(\mathbb C)$, le système hamiltonien avec potentiel  $V(s,-) $  est orthogonal aux constantes}. \\

\textbf{Organisation de l'article.} Cet article est organisé de la manière suivante:

La première partie est consacrée à l'étude des conséquences du théorème de Hrushovski-Sokolovic pour les équations différentielles autonomes.  On etudie notamment la propriété de désintégration (propriété (D)) ainsi que le théorème \ref{intro1}.

La seconde partie concerne l'orthogonalité aux constantes, en tant que telle. On y démontre l'équivalence avec sa formulation géométrique (propriété (O)) et le théorème \ref{intro3} de spécialisation (ou sa formulation contraposée de générisation, comme l'on préfère).

La troisième partie se concentrera l'aspect dynamique de l'étude des équations différentielles et jouera un rôle ``d'intermédiaire'' essentiel entre la catégorie des équations différentielles algébriques et la catégorie des systèmes dynamiques topologiques réels. 

La quatrième et dernière partie est consacrée à la preuve du critère dynamique d'orthogonalité aux  constantes (Théorème \ref{intro2})  et de sa version en famille (Théorème \ref{intro4}).

Enfin, cet article contient en appendice, des resultats formels concernant la catégorie des $D$-variétés, qui joue un rôle central dans les traductions entre les propriétés modèle-théoriques et leurs formulations géométriques. \\

\textbf{Remerciements.} Les résultats de cet article constituent une partie de ma thèse de doctorat, réalisée sous la direction de Jean-Benoît Bost (Orsay) et de Martin Hils (Paris VII -- Münster). Outre mes directeurs de thèse, je tiens à remercier Elisabeth Bouscaren et Zoe Chatzidakis pour leurs précieuses remarques sur le contenu présenté dans ce texte, lors des exposés que j'ai donnés à Paris VII et à Orsay. Je tiens aussi à remercier Rahim Moosa (Waterloo) pour de nombreuses discussions au sujet de cet article.

\tableofcontents

\section{Orthogonalité aux constantes}

Dans cette partie, nous discutons la propriété d'orthogonalité aux constantes pour un type stationnaire de la théorie $\textbf{DCF}_0$ défini sur un corps différentiel constant et sa relation avec la propriété de désintégration.

Il est bien connu que, dans toute théorie stable, un type stationnaire désintégré est orthogonal à tout type minimal non-localement modulaire (voir par exemple \citep[Lemma 2.3 pp158]{GST}. En particulier, un type stationnaire désintégré est toujours orthogonal aux constantes. En revanche, des constructions élémentaires utilisant des noyaux de Manin illustrent simplement que la réciproque n'est pas toujours vérifiée (voir $\S 14$ de l'introduction et le théorème 4.24 du chapitre 3 de \cite{jaoui}).  

Notre résultat principal concerne les équations différentielles autonomes: \textit{Une équation différentielle autonome dont le type générique est orthogonal aux constantes admet toujours un facteur rationnel} (non-trivial, c'est-à-dire de dimension strictement positive) \textit{désintégré}.

Ce résultat, propre à la théorie $\textbf{DCF}_0$  --- le corollaire \ref{desintegration} et sa formulation géométrique, le corollaire \ref{desintegrationgeometrique} --- est obtenu à l'aide de méthodes pures de théorie géométrique de la stabilité et du théorème de trichotomie de Hrushovski-Sokolovic. 

Nous supposerons une familiarité avec les résultats élémentaires de théories des modèles et de théorie de la stabilité. La lecture des chapitres 1-3 et 5-8 de \cite{Tent} est suffisante à la compréhension de cette partie. La première section est consacrée la relation d'orthogonalité entre deux types d'une même théorie stable. Nous exposons brièvement dans la seconde, quelques résultats sur la théorie $\textbf{DCF}_0$ dont nous nous servirons, dans la suite. Enfin, la troisième partie de cette section est consacrée à la preuve du théorème \ref{desintegrationgeometrique}.  

\subsection{Orthogonalité dans une théorie stable} Dans ce paragraphe, fixons $T$ une théorie du premier ordre, stable, complète, dans un langage $\mathcal L$ ainsi que $\Mod \models T$ un modèle monstre (relativement à un cardinal fixé $\kappa$) de la théorie $T$.

Un \textit{(petit) ensemble de paramètres} désigne un sous-ensemble $A \subset M$ de cardinal strictement inférieur à $\kappa$. Si $A$ est un ensemble de paramètres, on désigne par $S(A)$ l'ensemble des types (avec un nombre fini arbitraire de variables) au sens de la théorie $T$ à paramètres dans $A$.

Si $a,b$ des uplets de $M$, on rappelle que les deux uplets $a$ et $b$ sont dits \textit{indépendants au  dessus de $A$} si $\mathrm{tp}(a/A,b)$ est une extension non-déviante de $\mathrm{tp}(a/A)$. On notera alors $a \ind_A b$ pour désigner que $a$ et $b$ sont indépendants au  dessus de $A$. Cette relation vérifie les règles propriétés usuelles d'une relation d'indépendance \citep[Théorème 8.5.5]{Tent}. En particulier, elle est symétrique, i.e. $a \ind_A b$ si et seulement si $b \ind_A a$.

\subsubsection{Orthogonalité entre deux types stationnaires} 

{\nota Rappelons qu'un type $p \in S(A)$ est dit \textit{stationnaire} lorsqu'il admet une unique extension non-déviante à toute extension des paramètres $A \subset B$. Si $T$ est une théorie stable qui élimine les imaginaires alors tout type à paramètres dans un ensemble algébriquement clos est stationnaire.

Si $p$ est un type stationnaire et $A \subset B$ une extension des paramètres, on note $p|B$ son unique non-déviante à $B$.} 

{\defn Soient $p,q \in S(A)$ deux types stationnaires sur $A$. On dit que $p$ et $q$ sont \textit{faiblement orthogonaux} et on note $p \perp^a q$ si  deux  réalisations quelconques $a,b \in M$ de $p$ et $q$ respectivement sont indépendantes au dessus de $A$ i.e.
$$\text{Si  } a \models p \text{ et } \models q \text{ alors } a \ind_A b.$$}
Par symétrie de la relation d'indépendance dans une théorie $\omega$-stable, la notion d'orthogonalité faible est une relation symétrique.

{\rem De façon équivalente, deux types stationnaires $p(x) \in S(A)$ et $q(y) \in S(A)$ sont faiblement orthogonaux si et seulement si le type partiel $\pi(x,y) = p(x) \cup q(y)$ admet une unique complétion en un type complet.}

{\defn\label{ortho} Soient $A,B \subset M$ des ensembles, $p \in S(A)$ et $q \in S(B)$ deux types stationnaires. On dit que $p$ et $q$ sont \textit{orthogonaux} et on note $p \perp q$ si pour tout ensemble de paramètres $C \supset A \cup B$, les extensions non-déviantes respectives $p|C$ et $q|C$ de $p$ et $q$ à $C$ sont  faiblement orthogonales.}

\subsubsection{Principe de réflexivité de Shelah}

{\defn Soient $p,q \in S(A)$ deux types stationnaires. On appelle \textit{produit tensoriel de $p$ et $q$}, le type complet à paramètres dans $A$
$$p \otimes q = \mathrm{tp}(a,b/A)$$
où $a \models p$ réalise $p$ et $b \models q|A,a$ réalise l'unique extension non déviante de $q$ à $A,a$.
De même, pour tout $n \in \mathbb{N}$, on note $p^{\otimes n} = p \otimes \cdots \otimes p$ le produit tensoriel de $p$ avec lui même $n$ fois.}

Nous utiliserons l'avatar suivant du principe de réflexivité de Shelah pour la notion d'orthogonalité.

{\Prop[{\citep[Chapitre 1, Lemme 4.3.1]{GST}}] \label{ortho2} Soient $\Mod$ un modèle de $T$, $A \subset M$ un ensemble de paramètres et $p,q \in S(A)$ deux types stationnaires. On a équivalence entre : 
\begin{itemize}
\item[(i)] Les types $p$ et $q$ sont orthogonaux.
\item[(ii)] Pour tout $n \in \mathbb{N}$ et tout $m \in \mathbb{N}$, les types $p^{\otimes n}$ et $q^{\otimes m}$ sont faiblement orthogonaux.
\end{itemize}}

\subsubsection{Analyse  semi-minimale et non-orthogonalité} Rappelons qu'un type stationnaire $p \in S(A)$ dans une théorie stable $T$ est \textit{minimal} s'il est non-algébrique et toute extension déviante de $p$ est algébrique. 

{\Prop[{\citep[Chapitre 2, Remarque 2.10]{GST}}] Soit $\Mod$ un modèle $\kappa$-saturé de $T$. La relation "être non-orthogonal"  est une relation d'équivalence sur l'ensemble des types minimaux à paramètres dans $M$.}

La relation de non-orthogonalité cesse d'être une relation d'équivalence lorsqu'on ne se restreint plus au types minimaux dès que la théorie $T$ est multidimensionnelle --- i.e. admet au moins deux types minimaux orthogonaux. Néanmoins, elle permet de définir un dévissage des types $p$ de rang de Lascar fini en \textit{types semi-minimaux}, que nous décrivons maintenant.

{\lem[{\citep[Chapitre 2, Lemme 2.5.1]{GST}}]\label{orthominimal} Soient $\Mod$ un modèle $\kappa$-saturé de $T$, $A \subset M$ un ensemble de paramètres petit et $p \in S_n(A)$ un type de rang $\rm{RU}$ fini. Il existe une extension $A \subset B$ de paramètres et un type stationnaire minimal $q \in S(B)$ non-orthogonal à $p$.}

{\Prop\label{analyse} Soient $\Mod$ un modèle de $T$, $A \subset M$ un ensemble de paramètres et $p \in S_n(A)$ un type stationnaire de rang $\rm{RU}$ fini égal à $r \in \mathbb{N}$.
Il existe une suite d'extensions de paramètres $A \subset A_1 \subset A_2 \subset \cdots \subset A_r$, une suite d'extension $p \subset p_1 \subset p_2 \subset \cdots \subset p_{r}$ où $p_i \in S(A_i)$ est un type stationnaire et des types minimaux stationnaires $q_i \in S(A_i)$ vérifiant
\begin{itemize}
\item[(i)] Pour tout $i \leq r$, le rang de Lascar $RU(p_i)$ de $p_i$ est $r + 1 - i$. En particulier, le type $p_1$ est une extension non déviante de $p$ et $p_r$ est un type minimal.
\item[(ii)]Pour $i = r$, on a $p_r = q_r$ et pour tout $ 1 \leq i < r$, il existe des réalisations $a_i \models p_i$ et $b_i \models q_i$ telles que 
$ b_i \in \mathrm{acl}(A_i,a_i)$  et $p_{i+1}$ est une extension stationnaire non déviante à $A_{i+1}$ de $\mathrm{tp}(a_i/A_i,b_i)$.
\end{itemize}}

\begin{proof}
On raisonne par récurrence sur le rang de Lacar $r= \mathrm{RU}(p) \in \mathbb{N}$ de $p$.

Pour $r = 1$, il suffit de poser $A_1 = A$ et $p_1 = q_1 = p$. 

Supposons le résultat montré pour les types de rang de Lascar $r \in \mathbb{N}$. Considérons $p \in S(A)$ un type stationnaire de rang de Lascar $r + 1$.

D'après le lemme \ref{orthominimal}, il existe un ensemble $B \subset M$ de paramètres et $q \in S(B)$ un type stationnaire minimal non orthogonal à $p$.
Considérons $C \supset A \cup B$ un ensemble de paramètres et des réalisations $a \models p|C$ et $b \models q|C$ vérifiant $a \nind_C b$.

Le type $\mathrm{tp}(b/C)$ étant minimal, on en déduit que $b \in \mathrm{acl}(C,a)$. On pose alors $A_1 = C$, $q_1= q|C$, $p_1 = p|C$. Par additivité du rang de Lascar, le type $p = \mathrm{tp}(a/b,C)$ est de rang de Lascar $r \in \mathbb{N}$. On applique alors l'hypothèse de récurrence à une extension stationnaire de $p = \mathrm{tp}(a/b,C)$.
\end{proof}

{\rem Les données $\mathcal A = ((A_i)_{i \leq r}, (p_i)_{i \leq r}, (q_i)_{i \leq r})$ où  $A \subset A_1 \subset A_2 \subset \cdots \subset A_r$ est une extension de paramètres, et de types stationnaires $p_i,q_i \in S(A_i)$ dont l'existence est assurée par la proposition \ref{analyse} est appelée \textit{une analyse du type $p \in S(A)$ par les types minimaux $q_1 , \cdots q_r$.}

On note alors $[p]_\mathcal A = \lbrace [q_1] , \cdots , [q_r] \rbrace$ l'ensemble des classes de non-orthogonalité des types minimaux intervenant dans l'analyse de $p$.}

{\Prop\label{facteurs minimaux}  Soient $\Mod$ un modèle de $T$, $A \subset M$ un ensemble de paramètres et $p \in S_n(A)$ un type stationnaire de rang $\rm{RU}$ fini égal à $r \in \mathbb{N}$. Si $\mathcal A$ et $\mathcal A'$ sont deux analyses de $p$ alors $[p]_\mathcal A = [p]_{\mathcal A'}$.

Plus précisément, pour toute analyse $\mathcal A$ de $p$, on a 
$$ [p]_\mathcal A = \lbrace [q]\text{ } | \text{ } q \text{ est un type minimal non orthogonal à une extension de $p$} \rbrace.$$}

Si $p \in S(A)$ est un type de rang de Lascar fini, on notera $[p] = [p]_\mathcal A$ pour toute analyse $\mathcal A$ de $p$ dont l'existence est assurée par la proposition \ref{analyse} et qui est bien défini d'après la proposition \ref{facteurs minimaux}.

\begin{proof}
La deuxième partie de la proposition implique la première. La condition (ii) de la proposition \ref{analyse} montre que :
$$ [p]_\mathcal A \subset \lbrace [q]\text{ } | \text{ } q \text{ est un type minimal non orthogonal à une extension de $p$} \rbrace.$$
Il suffit donc de montrer l'inclusion réciproque : Considérons $\mathcal A$ une analyse de $p$ et un type $q \in S(B)$ minimal orthogonal à $[p]_\mathcal A$ (i.e. orthogonal à tous les éléments de $[p]_\mathcal A$).
Soit $A \subset C$ une extension de paramètres et $\widehat{p} \in S(C)$ une extension de $p$. Montrons que $q$ est orthogonal à $\widehat{p}$.

D'après \citep[Chapitre 8, Lemme 1.2(i)]{GST}, il existe une analyse $\widehat{ \mathcal A}$ de $\widehat{p}$ avec 
$$[\widehat p]_{\widehat{ \mathcal A}} \subset [p]_{\mathcal A}.$$   

Par hypothèse, le type $q$ est orthogonal à tous les éléments de  $[\widehat p]_{\widehat{ \mathcal A}}$ et donc d'après \citep[Chapitre 8, Lemme 1.2(iv)]{GST}, le type $q$ est orthogonal à tous les types analysables dans $[\widehat p]_{\widehat{ \mathcal A}}$ et en particulier à $\widehat{p}$.
\end{proof}

{\exam\label{produitpur} Soient $\Mod$ un modèle de $T$, $A \subset M$ un ensemble de paramètres et
$p_1 , \cdots , p_n \in S(A)$ des types minimaux stationnaires. On vérifie facilement que 
$$[p_1 \otimes \cdots \otimes p_n] = \lbrace [p_1] , \cdots , [p_n] \rbrace.$$}

\subsection{Corps différentiellement clos} On rappelle qu'un corps différentiel (de caractéristique $0$) est un couple $(K,\delta)$ où $K$ est un corps de \textit{caractéristique $0$} et $\delta : K \longrightarrow K$ une dérivation. Dans tout le texte, on travaillera en caractéristique $0$, on parlera donc de corps différentiel pour désigner un corps différentiel de caractéristique $0$. 

\subsubsection{Définition} 

Du point de vue syntaxique, un corps différentiel est une $\L_\delta$-structure où $\L_\delta = \lbrace 0,1,+,.,-,\delta \rbrace$ est appelé \textit{langage des corps différentiels} vérifiant les axiomes des corps différentiels (qui sont du premier ordre dans le langage $\L_\delta$).

{\rem Pour le langage $\mathcal L_\delta$ des corps différentiels, les formules sans quantificateurs à paramètres dans un corps différentiel $(K,\delta)$ sont les combinaisons booléennes d'équations différentielles algébriques, c'est-à-dire d'équations différentielles de la forme : 
$$ (E) : 
P(x_1 , \cdots , x_r ,  \delta(x_1) , \cdots , \delta(x_r), \cdots,  \delta^k(x_1), \cdots,  \delta^k(x_r))  = 0 $$
où $P \in K[X_1^{(0)},\cdots, X_r^{(0)}, \cdots , X_1^{(k)}, \cdots, X_r^{(k)}]$ est un polynôme.}

Parmi les modèles de la théorie des corps différentiels, certains possèdent une importance particulière, ce sont les modèles existentiellement clos.

{\defn Soient $T$ une théorie dans un langage $\mathcal L$ et $\Mod$ un modèle de $T$. On dit $\Mod$ est \textit{un modèle existentiellement clos de $T$} si pour toute extension de modèles $\Mod \subset \Mod'$ de $T$, toute formule sans quantificateurs $\phi(\overline{x},\overline{y})$ et tout $\overline{b} \in M^n$, on a : 
$$ \Mod' \models \exists \overline{x} \phi(\overline{x},\overline{b}) \Longleftrightarrow  \Mod \models \exists \overline{x} \phi(\overline{x},\overline{b}).$$}

{\defn On appelle \textit{corps différentiellement clos}, tout modèle existentiellement clos de la théorie des corps différentiels.}

\subsubsection{Propriétés globales de la théorie $\textbf{DCF}_0$}

{\Thm[Blum, Poizat] La classe des corps différentiellement clos est axiomatisable dans le langage $\mathcal L_\delta$ des corps différentiels par une théorie du premier ordre notée $\textbf{DCF}_0$ et appelée \textit{théorie des corps différentiellement clos}. Une axiomatisation de cette théorie est donnée par : 
\begin{itemize}
\item[$\textbf{DCF}_{0}^I$] $(K,\delta)$ est un corps différentiel (de caractéristique 0).
\item[$\textbf{DCF}_{0}^{II}$] Pour tout $n \in \mathbb{N}$, et tous polynômes $P \in K[X_1^{(0)},\cdots, X_1^{(n)}] \setminus K[X_1^{(0)},\cdots, X_1^{(n-1)}]$ et $Q \in K[X_1^{(0)},\cdots, X_1^{(n - 1)}] $ non nul, il existe $y \in K$ tel que 
$$ P(y,\delta(y),\cdots , \delta^n(y)) = 0 \wedge Q(y,\delta(y),\cdots , \delta^{n-1}(y)) \neq 0.$$
\end{itemize}
De plus, la théorie $\textbf{DCF}_0$ est complète et élimine les quantificateurs et les imaginaires dans le langage $\mathcal L_\delta$ des corps différentiels.}

{\cor[Description des types]\label{descriptiondestypesaffine}  Soient $(\mathcal U,\delta_U)$ un corps différentiellement clos et $(K,\delta) \subset (\mathcal U,\delta_U)$ un sous-corps différentiel. L'application définie par
$$ I : \begin{cases} S_n(K) \longrightarrow & \mathrm{Spec}_\delta(K\lbrace X_1 , \cdots X_n \rbrace) \\
 p \mapsto & \lbrace P \in K \lbrace X_1 , \cdots X_n \rbrace \text{ $|$ } ``P(x) = 0'' \in p \rbrace
\end{cases}$$
est une bijection, où $\mathrm{Spec}_\delta(K\lbrace X_1 , \cdots , X_n\rbrace)$ désigne l'ensemble des idéaux premiers différentiels de la $(K,\delta)$-algèbre différentielle libre $K\lbrace X_1 , \cdots , X_n\rbrace$ avec variables $X_1,\ldots X_n$.}

{\cor La théorie $\mathbf{DCF}_0$ est $\omega$-stable. }  

\begin{proof}
D'après le théorème de Ritt-Raudenbush (see \citep[Partie 2, Théorème 1.16]{MTF}), on peut choisir pour tout idéal premier différentiel $I \subset K\lbrace X_1,\cdots, X_n \rbrace$, un système générateur $S(I)$ de l'idéal $I$ en tant qu'idéal différentiel radical.
La description précédente de l'ensemble $S_n(A)$ des types à paramètres dans $A$ montre que l'application 
$$\begin{cases}
S_n(A)  \longrightarrow & \mathbb{Q} \langle A \rangle \lbrace X_1 , \cdots , X_n \rbrace^{(\mathbb{N})} \\
p   \mapsto & S(I(p)) 
\end{cases}$$
est une injection. On en déduit que $S_n(A)$ est dénombrable dès que $A$ est dénombrable et donc que la théorie $\mathbf{DCF}_0$ est $\omega$-stable. 
\end{proof}

\subsubsection{Structure induite sur le corps des constantes}

{\defn Soit $(K,\delta)$ un corps différentiel. Le \textit{corps des constantes} de $(K,\delta)$ est le sous-corps de $K$ noté $K^\delta$ et défini par:
$$ K^\delta = \lbrace x \in K \text{ } | \text{ } \delta(x) = 0 \rbrace.$$}
  
{\cor[{\citep[Partie 2, Lemme 5.10]{MTF}}]\label{corpsdesconstantes} Soit $(\mathcal U,\delta_\U)$ un corps différentiellement clos. Le corps des constantes $\mathcal C$ de $(\mathcal U,\delta_U)$ est un corps algébriquement clos et la structure induite par $\U$ sur le sous-corps définissable $\mathcal C$ est celle d'un pur corps algébriquement clos\footnote{Autrement dit, pour tout sous-ensemble de paramètres petit $B\subset \U$, les sous-ensembles $B$-définissables de $\mathcal C$ sont les sous-ensembles $\mathrm{dcl}(B)\cap \mathcal C$-définissables de $\mathcal C$ dans le langage des anneaux.}.}

Un corps différentiellement clos induit donc une structure de pur corps algébriquement clos sur son corps des constantes. En particulier le corps des constantes d'un corps différentiellement clos est un sous-ensemble définissable fortement minimal.

{\defn Soient $(\U,\delta_\U)$ un corps différentiellement clos et $p \in S(A)$ un type à paramètres dans un sous-ensemble $A \subset \U$. On dit que le type $p$ est \textit{orthogonal aux constantes} si $p$ est orthogonal au type générique de l'ensemble fortement minimal
$$\mathcal U^\delta = \lbrace x \in \U \text{ } | \text{ } \delta(x) = 0 \rbrace.$$}
\subsection{Orthogonalité aux constantes et désintégration}

\subsubsection{Propriété de désintégration} On fixe $k$ un corps de caractéristique $0$.

{\rem On travaille dans ce qui suit dans la catégorie des $D$-schémas $(X,\delta_X)$ au dessus de $(k,0)$ au sens de l'appendice A --- c'est-à-dire des schémas $X$ au dessus de $k$ munis d'une dérivation $\delta_X : \mathcal O_X \longrightarrow \mathcal O_X$ sur le faisceau structurel de $X$. Lorsque $X$ est un schéma intègre lisse, une $D$-variété $(X,\delta_X)$ est \textit{tout simplement une variété irréductible $X$ munie d'un champ de vecteurs $v$}.}

{\nota Soit $(X,v)$ une $D$-variété absolument irréductible au dessus de $(k,0)$. Pour tout $n > 0$, on note $\mathcal I^{gen}_n(X,v)$, \textit{l'ensemble des sous-schémas fermés intègres invariants de $(X,v)^n$ qui se projettent génériquement sur tous les facteurs}.

Notons que d'après la proposition \ref{schemainvariantreduit}, un sous-schéma fermé réduit  $Y \subset (X,v)^n$ est invariant si et seulement si ses composantes irréductibles le sont.}

{\defn\label{desintegrationdefinition} Soit $(X,v)$ une variété absolument irréductible $X$ munie d'un champ de vecteurs $v$ au dessus de $k$. On dit que $(X,v)$ est \textit{génériquement désintégrée} si pour tout $n \geq 3$, tout élément $Z \in \mathcal I^{gen}_n(X,v)$ s'écrit comme une composante irréductible (se projetant génériquement sur tous les facteurs) de:
\begin{eqnarray}\label{equationdesintegration}
\bigcap_{1 \leq i < j \leq n} \pi_{i,j}^{-1}(Z_{i,j}).
\end{eqnarray}
où $\pi_{i,j} : X^n \longrightarrow X^2$ est la projection sur les $i^{eme}$ et $j^{eme}$ coordonnées et $Z_{i,j} \in \mathcal I^{gen}_2$ pour tout $i \neq j$.}

{\rem Les propriétés suivantes sont des conséquences immédiate de la définition précédente: 
\begin{itemize}

\item[(1)] Dans la définition \ref{desintegrationdefinition}, on peut remplacer l'identité (\ref{equationdesintegration}) par:
$$\bigcap_{1 \leq i \neq j \leq n} \pi_{i,j}^{-1}(Z_{i,j}).$$
Cette forme est particulièrement agréable lorsque l'on travaille avec des notations multi-indices (par exemple Lemme \ref{desintegrationproduit}).

\item[(2)] Si $(X,v)$ est génériquement désintégrée, alors pour tout $n \geq 3$, tout $Z \in \mathcal I^{gen}_n(X,v)$ s'écrit comme une composante irréductible de  
$$\bigcap_{1 \leq i < j \leq n} \pi_{i,j}^{-1}(\overline{\pi_{i,j}(Z)}).$$

\item[(3)] Dans la définition \ref{desintegrationdefinition}, il n'est pas nécessaire de supposer que les $Z_{i,j}$ sont irréductibles. En effet, si $Z_{1,2}$ se décompose en une union de sous-variétés fermés stricts  $Z_{1,2} = F \cup G$ alors ces dernières sont invariantes (voir Appendice A) et:
$$ \bigcap_{1 \leq i < j \leq n} \pi_{i,j}^{-1}(Z_{i,j}) =  \Big( \pi_{1,2}^{-1}(F) \cap  \bigcap_{(i,j) \neq (1,2)} \pi_{i,j}^{-1}(Z_{i,j})\Big)  \cup \Big( \pi_{1,2}^{-1}(G) \cap \bigcap_{(i,j) \neq (1,2)} \pi_{i,j}^{-1}(Z_{i,j}) \Big).$$
\end{itemize}
En particulier, les composantes irréductibles du membre de gauches sont les composantes irréductibles des deux facteurs du membre de gauche.}

{\lem\label{desintegrationminimal} Soit $(X,v)$ une variété absolument irréductible $X$ munie d'un champ de vecteurs $v$ au dessus de $k$. Si le type générique de $(X,v)$ est minimal et désintégré alors $(X,v)$ est génériquement désintégrée.}

Cet énoncé est une variante des résultats bien connus de théorie des modèles. Par manque de référence, nous formulons une preuve rapide.

\begin{proof}
Considérons $Z \in \mathcal I_n^{gen}(X,v)$ une sous-variété fermée invariante de $(X,v)^n$ se projetant génériquement sur tous les facteurs.

Considérons $n_0 \in \mathbb{N}$ maximal tel que $Z$ se projette génériquement sur $X^{n_0}$ via les projections coordonnées. Sans perte de généralité, on peut supposer que $Z$ se projette génériquement sur les $n_0$ premières coordonnées.

Fixons maintenant un corps differentiellement clos $(\U,\delta_\U)$ et considérons $(a_1,\ldots, a_n)$ une réalisation du type générique de $Z$ (comme $Z$ se projette génériquement sur chacun des facteurs, chaque $a_i$ réalise le type générique de $(X,v)$).

Comme le type générique de $X$ est minimal et désintégré, on en deduit que pour tout $k > n_0$, 
$$a_k \in \mathrm{acl}(a_1,\ldots, a_{n_0}) = \mathrm{acl}(a_1) \cup \cdots \cup \mathrm{acl}(a_{n_0}).$$

Pour tout $k > n_0$, il existe donc $i(k) \leq n_0$ et une correspondance génériquement finie $Y_k \subset X^2$  telle que $(a_k, a_{i(k)}) \in Y_k$. Posons alors
$$ Y = \bigcap_{k > n_0} \pi_{k,i(k)}^{-1}(Z_k).$$

Par construction, la fibre générique de $Y$ par la projection sur les $n_0$ premières coordonnées est finie (bornée par le produit des fibres des $Y_k$ pour $k > n_0$). De plus, on a $(a_1, \ldots, a_n) \in Y$ et donc $Z \subset Y$. On en déduit que:
$$\mathrm{dim}(Z) = \mathrm{dim}(Y) = n_0.\mathrm{dim}(X)$$
et donc que $Z$ est une composante irréductible de $Y$.
\end{proof} 

\subsubsection{Quelques propriétés structurelles de la désintégration}

{\lem\label{desintegrationproduit} Soit $(X,v)$ une variété absolument irréductible $X$ munie d'un champ de vecteurs $v$ au dessus de $k$. Si $(X,v)$ est génériquement désintégrée et $n \geq 1$ alors $(X,v)^n$ l'est aussi.}

\begin{proof}
Considérons $Z \in \mathcal I_m^{gen}((X,v)^n)$ une sous-variété fermée invariante de $((X,v)^n)^m$ se projetant génériquement sur tous les facteurs. On utilise une notation double-indice: pour tout $1 \leq i \leq n$ et $1 \leq k \leq m$, $(i,k)$ désigne la $i$-ème coordonne du $k$-ème facteur $(X,v)^n$.

Puisque $(X,v)$ est désintégré, il existe $Z_{i,j,k,l} \in \mathcal I_2^{gen}(X,v)$ pour tout $(i,k) \neq (j,l)$ telles que $Z$ est une composante irréductible de
$$ Y = \bigcap_{(i,k) \neq (j,l)} \pi_{i,k,j,l}^{-1} (Z_{i,j,k,l}).$$  

Comme $Z$ se projette génériquement sur tous les facteurs on peut supposer $Z_{i,j,k,k} = X^2$ pour $i \neq j$. On peut donc réécrire:
$$ Y = \bigcap_{k \neq l} \Big( \bigcap_{i,j} \pi_{i,k,j,l}^{-1} (Z_{i,j,k,l}) \Big) = \bigcap_{k \neq l} Y_{k.l}.$$ 

Fixons $k \neq l$. En notant $\pi_{k,l}$ la projection sur le $k$-ème et $l$-ème facteur et $\pi_{i,j}$ la projection de $(X,v)^n \times (X,v)^n$ sur la $i$-ème coordonnée du premier facteur et la $j$-ème du second, on a $\pi_{i,j,k,l} = \pi_{i,j} \circ \pi_{k,l}$. On en déduit que:
$$ Y_{k,l} = \bigcap_{i,j} \pi_{k,l}^{-1}(\pi_{i,j}^{-1}(Z_{i,j,k,l})) = \pi_{k,l}^{-1}( \bigcap_{i,j} \pi_{i,j}^{-1}(Z_{i,j,k,l})).$$

Posons alors $Z_{k,l} = \bigcap_{i,j} \pi_{i,j}^{-1}(Z_{i,j,k,l})$ qui est une sous-variété fermée invariante de $(X,v)^n \times (X,v)^n$ se projetant génériquement sur les deux facteurs. On conclut que $Z$ est une composante irréductible de $Y$ qui s'écrit:
$$ Y = \bigcap_{k \neq l} Y_{k,l} =  \bigcap_{k \neq l} \pi_{k,l}^{-1}(Z_{k,l}).$$
\end{proof}

{\Prop\label{generiquementfinidesintegration}  Soient $(X,v)$ et $(Y,w)$ deux variétés absolument irréductibles munies respectivement de champs de vecteurs $v$ et $w$ au dessus de $k$ et $f: (X,v) \dashrightarrow (Y,w)$ un morphisme rationnel génériquement fini de $D$-variétés au dessus de $k$.

L'équation différentielle $(X,v)$ est génériquement désintégrée si et seulement si $(Y,w)$ est génériquement désintégrée.}

La preuve de la proposition précédente est construire a partir du lemme de géométrie algébrique pure (non différentielle) suivant:

{\lem  Soit $f : X \longrightarrow Y$ est un morphisme étale au dessus d'un corps $k$ de caractéristique $0$, $Z_Y$ une sous-variété fermée irréductible de $Y$ et $Z_X$ une composante irréductible de $f^{-1}(Z_Y)$. Pour toute sous-variété fermée $Z$  de $Y$ les propriétés suivantes sont équivalentes:
\begin{itemize}
\item[(i)] $Z_X$ est une composante irréductible de $f^{-1}(Z)$.
\item[(ii)] $Z_Y$ est une composante irréductible de $Z$.
\end{itemize}}

\begin{proof} 
Supposons $(i)$.  Il suffit de montrer que $Z_Y$ est une sous-variété fermée irréductible maximale de $Z$.

Considérons une sous-variété fermée irréductible $T$ de $Y$ telle que $Z_Y \subset T \subset Z$. On en déduit que $Z_X$ est inclus dans une composante irréductible $T_1$ de $f^{-1}(T)$. Par maximalité de $Z_X$, on a $Z_X = T_1$. On en déduit que:
$$ Z_Y = \overline{f(Z_X)} = \overline{f(T_1)} = T.$$ 

Réciproquement supposons (ii) et considérons une sous-variété fermée irréductible $T$ de $X$  telle que $Z_X \subset T \subset f^{-1}(Z)$

En appliquant $f$, on obtient $Z_Y \subset \overline{p(T)} \subset \overline Z$ et donc $Z_Y = \overline{p(T)}$ (par maximalité de $Z_Y$). Il suit que $T$ est une composante irreductible de $ f^{-1}(\overline{p(T)}) = f^{-1}(Z_Y)$ qui contient $Z_X$ et donc $T = Z_X$. 
\end{proof}
\begin{proof}[Demonstration de la proposition \ref{generiquementfinidesintegration}]
La propriété de désintégration est préservée en replaçant $X$ et $Y$ par des ouverts non vides, on peut donc supposer que $f: (X,v) \longrightarrow (Y,w)$ est un morphisme étale surjectif.

Pour tout $n \geq 1$, on note $f_n = f \times \cdots \times f : X^n \longrightarrow Y^n$ qui est aussi un morphisme étale. On a un diagramme commutatif:
 $$ \xymatrix{
    X^n \ar[r]^{f_n} \ar[d]_{\pi^X_{i,j}} & Y^n \ar[d]^{\pi^Y_{i,j}} \\
    X^2 \ar[r]^{f_2} & Y^2
  }$$

Considérons $Z_Y \in  \mathcal I_n^{gen}(Y,v)$ et considérons $Z_X$ une composante irréductible de $f^{-1}(Z_Y)$. Nous montrons que les propriétés de désintégration pour $Z_X$ et $Z_Y$ sont équivalentes:

Posons $X_{i,j} = \overline{\pi^X_{i,j}(Z_X)} \in \mathcal I^{gen}_2(X,v)$ et  $Y_{i,j} = \overline{f_2(X_{i,j})} \in \mathcal I_2^{gen}(Y,w)$. Comme le morphisme $f_2$ est étale,  $X_{i,j}$ est une composante irréductible $f_2^{-1}(\overline{f_2(X_{i,j})})$.


Comme le morphisme $f_2$ est étale,  $X_{i,j}$ est une composante irréductible $f_2^{-1}(\overline{f_2(X_{i,j})})$. En utilisant le point (3) de la remarque après la définition \ref{desintegrationdefinition}, on en deduit que l'ensemble des composantes irréductibles de $$\bigcap_{1 \leq i \neq j \leq n} (\pi_{i,j}^{X})^{-1}(X_{i,j})$$ est inclus dans l'ensemble des composantes irréductibles de:
$$\bigcap_{1 \leq i \neq j \leq n} (\pi_{i,j}^{X})^{-1}(f_2^{-1}(Y_{i,j})) = f_n^{-1} \Big( \bigcap_{1 \leq i \neq j \leq n}(\pi_{i,j}^{Y})^{-1}(Y_{i,j}) \Big). $$

Ainsi $Z_X$ vérifie la propriété de désintégration si et seulement si c'est une composante irréductible du membre de droite si et seulement $Z_Y$ est une composante irréductible de 
$$\bigcap_{1 \leq i \neq j \leq n}(\pi_{i,j}^{Y})^{-1}(Y_{i,j}). \qedhere$$
\end{proof}

\subsubsection{De la minimalité à la semi-minimalité} Dans cette section, on étend le lemme \ref{desintegrationminimal} au cas où le type générique de l'équation différentielle $(X,v)$ est seulement semi-minimal et désintégré. Ce raffinement sera necessaire pour la preuve du corollaire \ref{desintegrationgeometrique}.

{\defn Soient $T$ une théorie stable, $A$ un ensemble de paramètres et $p \in S(A)$ un type à paramètres dans $A$.

On dit que $p$ est \textit{semi-minimal} s'il existe une extension des paramètres $A \subset B$ et $q \in S(B)$ un type minimal tels que le type $p$ est presque-interne à l'ensemble des $\mathrm{acl}(A)$-conjugués de $q$. On dit que le type $p$ est \textit{semi-minimal et désintégré} si, de plus, le type minimal $q \in S(B)$ est désintégré.}   

On a une description particulièrement simple des types semi-minimaux désintégrés si l'on s'autorise à remplacer le type $q$ par un autre type minimal qui lui est non-orthogonal:

{\Prop\label{parametredesintegre} Soient $T$ une théorie stable, $A$ un ensemble de paramètres et $p \in S(A)$ un type semi-minimal et désintégré à paramètres dans $A$.

Il existe un type minimal désintégré $r \in S(A)$ (non-orthogonal à $q$) et des réalisations $A$-indépendantes $a_1,\ldots, a_j$ de $r$ et une réalisation $b$ de $p$ telles que:
$$ \mathrm{acl}_A(b) = \mathrm{acl}_A(a_1, \ldots a_j).$$}  

\begin{proof} 
On travaille dans un modèle saturé $\U$. Le type $q$ est un type minimal et désintégré et donc en particulier localement modulaire. Cela implique que tout type analysable dans des conjugués de $q$ est monobasé. En particulier, le type $p$ est un type monobasé.

Soit $b$ une réalisation de $p$. Considérons $e$ tel que $\mathrm{RU}(b/eA) = \mathrm{RU}(b/A) - 1$ et tel que $e = \mathrm{Cb}(b/e)$. Comme le type $p$ est monobasé, on a donc $e \in \mathrm{acl}(A,b)$.
On en déduit facilement que:
$$\mathrm{RU}(e/A) = \mathrm{RU}(be/A) - \mathrm{RU}(b/Ae) = 1.$$

On en déduit que $r = \mathrm{tp}(e/A)$ est un type minimal non-orthogonal à $p$ et donc à $q$ (car le type $p$ est $q$-semi minimal). On a donc construit un type minimal et désintégré $r \in S(A)$ tel que $p$ est presque interne à $r$. Quitte a remplacer $p$ par un type interalgebrique, on peut donc supposer que $p$ est interne à $r$. \\

Montrons maintenant que $b \in \mathrm{acl}_A(e_1,\ldots, e_n)$ pour des réalisations $e_1,\ldots, e_n$ de $r$, c'est-à-dire que $b$ a une orbite finie sous l'action du groupe $\mathrm{Aut}(\U / \mathcal R)$ des automorphismes de $\U$ fixant point par point l'ensemble $\mathcal R$ des réalisations de $r$.

Comme le type $p$ est interne à $r$, le groupe de liaison $\mathrm{Aut}_A(\mathcal P / \mathcal R)$ des permutations de $\mathcal P$ provenant d'un automorphisme de $\U$ fixant $\mathcal R$ est isomorphe à un groupe définissable dans $\mathcal R^{eq}$ (voir \citep[Theorem 3]{Hru6}). La propriété de désintégration de $r$ garantit la finitude du groupe $G$ et donc de l'orbite de $b$. \\

On a donc obtenu $b \in \mathrm{acl}_A(e_1,\ldots, e_n)$ pour des réalisations $e_1,\ldots, e_n$ de $r$ et donc que $\mathrm{acl}_A(b) \subset \mathrm{acl}_A(e_1,\ldots ,e_n)$. Par minimalité de $r$, cet ensemble $\mathrm{acl}$-clos admet une base de transcendance de réalisations de $r$: il existe des réalisations $A$-indépendantes $a_1,\ldots a_j$ de $r$ telles que
$$\mathrm{acl}_A(b) = \mathrm{acl}_A(a_1,\ldots a_j).$$     
\end{proof}

{\cor\label{desintegrationsemiminimale} Soit $Y$ variété absolument irréductible au dessus d'un corps $k$ munie d'un champ de vecteurs $w$. Si le type générique de $(Y,w)$ est semi-minimal et désintégré alors $(Y,w)$ est génériquement désintégrée.}

\begin{proof}
Considérons un corps différentiellement clos $(\U,\delta_\U)$ et $b$ réalisant le type générique de $(X,v)$. D'après la proposition \ref{parametredesintegre}, il existe un type minimal et désintégré $q \in S(k)$ et des réalisations $k$-indépendantes $(a_1,\ldots, a_j)$ de $q$ telles que:  
$$(\ast) \qquad \mathrm{acl}_k(b) = \mathrm{acl}_k(a_1, \ldots a_j).$$

Comme le type $q$ est un type de rang fini à paramètres dans $k$, il existe une variété absolument irréductible $Y$ et champ de vecteurs $w$ sur $Y$ tel que $q$ est le type générique de $(Y,w)$.

Comme les $a_i$ sont des réalisations $k$-indépendantes de $q$, on en déduit que $a_1,\ldots, a_j$ réalise le type générique de $(Y,w)^j$.
La propriété $(\ast)$ exprime donc que les types génériques de $(Y,w)^j$ et $(X,v)$ sont interalgebriques, c'est-à-dire qu'il existe une correspondance génériquement finie entre $(X,v)$ et $(Y,w)^j$. On conclut à l'aide du lemme \ref{desintegrationproduit} et de la proposition \ref{generiquementfinidesintegration}.
\end{proof}

\subsubsection{Désintégration et orthogonalité aux constantes} Le théorème de trichotomie dans les corps différentiellement clos accorde une place toute particulière au corps des constantes d'un corps differentiellement clos: il s'agit d'un représentant ``canonique'' de l'unique classe de non-orthogonalité non-localement modulaire dans un corps differentiellement clos.

{\Thm[Théorème de trichotomie pour $\textbf{DCF}_0$, \cite{Sok}] Pour la théorie $\textbf{DCF}_0$, on a la description suivante des classes d'équivalence de non-orthogonalité pour les types minimaux: 
\begin{itemize}
\item[(i)] Il existe une unique classe d'équivalence non-localement modulaire représentée par le type générique du corps des constantes.
\item[(ii)] Les classes de non-orthogonalité, localement modulaires non-désintégrées sont en correspondance biunivoque avec les variétés abéliennes simples définies sur $\U$ modulo isogénie, qui ne descendent pas aux constantes.
\item[(iii)] Il existe une infinité non-dénombrable de classes de non-orthogonalité désintégrées. 
\end{itemize}}

{\cor\label{desintegration} Soient $(k,0)$ un corps différentiel constant et $p \in S(k)$ un type non-algébrique et de rang fini.

On suppose que le type $p$ est orthogonal aux constantes. Alors il existe un type semi-minimal désintégré $q \in S(k)$ et des réalisations $a \models p$ and $b \models q$ tels que:
$$ b \in \mathrm{dcl}(k,a) \setminus \mathrm{acl}(k).$$}
\begin{proof}
En effet, puisque le type $p$ est non-algébrique, il existe une extension des paramètres $(k,0) \subset (L,\delta)$ et un type minimal $r \in S(L)$ non-orthogonal à $p$.

Considérons $a \models p$, une réalisation de $p$. D'après le lemme de décomposition (voir, par exemple, \citep[Chapter 7, Lemma 4.5]{GST}, il existe $ b \in \mathrm{dcl}(k,a) \setminus \mathrm{acl}(k)$ tel que  $\mathrm{tp}(b/k)$ est (non-algébrique) et interne a l'ensemble des $\mathrm{acl}(k)$-conjugues de $r$.

On en déduit que le type $q = \mathrm{tp}(b/k)$ est $r$-semi-minimal. Il suffit alors de montrer que le type $r$ est un type minimal désintégré:

On sait que le type $p$ est orthogonal aux constantes. Comme $r$ est un type minimal non-orthogonal à $p$, on en déduit que $r$ est orthogonal aux constantes. D'après le théorème de Hrushovski-Sokolovic, le type $r$ est \textit{localement modulaire}.

Supposons par l'absurde que $r$ est localement modulaire et non désintégré. La classification fine des types minimaux localement modulaires non désintégrés obtenue par Hrushovski et Sokolovic implique que $r$ est orthogonal à tout type défini sur un corps différentiel constant (voir \citep[Theorem 2.8]{Sok}).

Or, on sait que le type $p$ est défini sur un corps différentiel constant et que $r$ est un type minimal non-orthogonal à $p$, ce qui est une contradiction. On en déduit que $r$ est un \textit{type minimal désintégré}.
\end{proof}

Comme expliqué dans l'introduction, le corollaire \ref{desintegration} permet d'obtenir l'existence de facteurs rationnels désintégrés pour les équations différentielles autonomes (définies sur un corps différentiel constant):

{\cor\label{desintegrationgeometrique} Soient $k$ un corps de caractéristique $0$ et $(X,v)$ une variété absolument irréductible au dessus de $k$ munie d'un champ de vecteurs $v$.

Si le type générique de $(X,v)$ est orthogonal aux constantes alors il existe une variété absolument irréductible $Y$ au dessus de $k$ de dimension $> 0$, un champ de vecteurs $w \in \mathrm{H}^0(Y, \Theta_{Y/k})$ et un morphisme rationnel dominant 
$$ \pi : (X,v) \dashrightarrow (Y,w)$$
tels que $(Y,w)$ est génériquement désintégrée.}

Quand les conclusions du corollaire \ref{desintegrationgeometrique} sont vérifiées, on dit que $(X,v)$ admet \textit{un facteur rationnel désintégré}.
\begin{proof}
On note $p$ le type générique de $(X,v)$ et on travaille dans un corps differentiellement clos $(\mathcal U, \delta)$.

D'après le corollaire \ref{desintegration}, il existe un type $q \in S(k)$ semi-minimal et désintégré  et des réalisations $a \models p$ and $b \models q$ tels que:
$$ b \in \mathrm{dcl}(k,a) \setminus \mathrm{acl}(k).$$

Maintenant, il existe une $D$-variété $(Y,w)$  au dessus de $(k,0)$ tel que $q$ est (interdéfinissable avec) le type générique de $(Y,w)$. Le corollaire \ref{desintegrationsemiminimale} assure que la $D$-variété $(Y,w)$  est génériquement désintégrée.

En utilisant que $b \in \mathrm{dcl}(a,k)$, on obtient, à l'aide du lemme \ref{morphismerationnel}, un morphisme rationnel dominant $ \pi : (X,v) \dashrightarrow (Y,w)$. Enfin, puisque $a \notin \mathrm{acl}(k)$, la variété $Y$ est de dimension $> 0$. 
\end{proof}

\section{Intégrale première rationnelle et orthogonalité aux constantes}

Dans cette partie, on étudie la relation entre la propriété d'orthogonalité aux constantes pour le type générique d'une équation différentielle autonome et les propriétés classiques de ``non-intégrabilité algébrique'' qui concernent l'absence de quantités algébriques conservées le long du mouvement (\textit{intégrales premières rationnelles}).

Étant donnée une équation différentielle algébrique (absolument irréductible) $(X,v)$,  une des caractéristiques nouvelles des méthodes issues de la théorie des modèles est de considérer l'équation différentielle  $(X,v)$ et toutes ses puissances $(X,v)^n$ sur le même plan.  Le critère d'orthogonalité aux constantes (Théorème \ref{orthogonaliteconstantes}) que nous démontrons dans la première section est une illustration parfaite de ce phénomène: \textit{Le type générique de $(X,v)$ est orthogonal est constantes si et seulement si pour tout $n > 0$, toute intégrale première rationnelle de $(X,v)^n$ est constante}. \\

La suite de cette partie est alors consacrée à l'étude des conséquences de ce premier résultat pour les familles lisses d'équations différentielles autonomes absolument irréductibles. 
Plutôt qu'une équation différentielle autonome, on considère alors une famille lisse $f: (\mathcal X,v) \longrightarrow (S,0)$ d'équations différentielles algébriques à paramètres constants au dessus d'un corps différentiel constant $(k,0)$ --- ou, en d'autres termes, une famille lisse de variétés algébriques $f : \mathcal X \longrightarrow S$ munie d'un champ de vecteurs $v \in \mathrm{H}^0(\mathcal X, \Theta_{X/S})$, tangent au fibres de $f$.

Si $s$ est un point de $S$, le champ de vecteurs $v$ sur $\mathcal X$ induit alors un champ de vecteurs sur la fibre $\mathcal X_s$. Autrement dit, la fibre $(\mathcal X,v)_s$ est une $D$-variété au dessus du corps différentiel (constant) $(k(s),0)$ que l'on supposera absolument irréductible. On étudie l'ensemble:
$$ S^{\not\perp 0} = \lbrace s \in S(\overline{l}) \text{ } |  \text{ le type générique de } (X,v)_s \text{ est non-orthogonal aux constantes} \rbrace$$     
où $\overline{l}$ est une extension algébriquement close de $k$. A ce sujet, il est bien connu que $S^{\not\perp 0}$  est, en général, plus complexe qu'un simple sous-ensemble définissable de $S(\overline{l})$. En revanche, des considérations très générales permettent de décrire $S^{\not\perp 0}$ comme l'image inverse d'un  sous ensemble de $S$ par l'application:
$$ \pi:  S(\overline{k})  \longrightarrow S$$
qui envoie un point de $S(\overline{k})$ sur son image dans le schéma $S$.

Dans la section 2.3, nous démontrons que, sous des hypothèses nécessaires de lissité, l'ensemble  $S^{\not\perp 0}$ est stable par spécialisation: si $s_\ast, s \in S(\overline{l})$ sont deux points et $\mathrm{tp}(s_\ast /k)$ se spécialise (au sens usuel de la théorie $\textbf{ACF}_0$) en $\mathrm{tp}(s/k)$ alors:
$$s_\ast \in S^{\not\perp 0} \Longrightarrow s \in S^{\not\perp 0}.$$

En outre, nos résultats s'appuient sur des énoncés classiques de spécialisation pour les intégrales premières rationnelles que nous rappelons (faute de référence appropriée) dans la deuxième sous-section de cette partie.

\subsection{Un premier critère d'orthogonalité aux constantes} On donne une caractérisation des  $D$-variétés absolument irréductibles dont le type générique est orthogonal aux constantes.

{\defn Soient $(K,\delta)$ un corps différentiel et $(X,\delta_X)$ une $D$-variété irréductible sur $(K,\delta)$. On appelle \textit{intégrale première rationnelle de $(X,\delta_X)$}, toute fonction rationnelle $f \in K(X)$ telle que $\delta_X(f) = 0$.}

{\rem Soient $(K,\delta)$ un corps différentiel et $(X,\delta_X)$ une $D$-variété irréductible sur $(K,\delta)$.
L'ensemble des intégrales premières rationnelles de $X$ s'identifie au corps des constantes du corps différentiel $(K(X),\delta_X)$.

Nous dirons que $(X,\delta_X)$ est \textit{sans intégrale première non constante} si $K(X)^\delta = K^\delta$.}

{\lem\label{sansintegralepremiererationnelle} Soient $(K,\delta)$ un corps différentiel et $(X,\delta_X)$ une $D$-variété irréductible au-dessus de $(K,\delta)$. On a équivalence entre : 
\begin{itemize}
\item[(i)] La $D$-variété irréductible $(X,\delta_X)$ est sans intégrale première rationnelle.
\item[(ii)] Pour toute variété $Y$ de dimension $> 0$ au dessus du corps des constantes $k = K^\delta$ de $(K,\delta)$, il n'existe pas de morphisme rationnel dominant de $D$-schémas au dessus de $(K,\delta)$
$$f : (X,\delta_X) \dashrightarrow (Y,0) \times_{(k,0)} (K,\delta).$$
\end{itemize}}

\begin{proof}
On raisonne par contraposition. Soit $f \in (K(X),\delta_X)$ une intégrale première rationnelle non constante. La fonction $f$ détermine un morphisme rationnel dominant de $D$-schémas au dessus de $(K,\delta)$
$$f : (X,\delta_X) \longrightarrow (\mathbb{A}^1,0) \times_{(k,0)} (K,\delta) $$
On en déduit que $(ii) \Longrightarrow (i)$.

Réciproquement, supposons qu'il existe un morphisme $(X,\delta_X) \dashrightarrow (Y,0) \times_{(k,0)} (K,\delta)$ rationnel dominant. Considérons $f \in k(Y)$ une fonction rationnelle non constante. Le morphisme rationnel dominant obtenu par composition
$$ (X,\delta_X) \dashrightarrow (Y,0) \times_{(k,0)} (K,\delta) \dashrightarrow 
(\mathbb{A}^1,0) \times_{(k,0)} (K,\delta)$$
est un morphisme de $D$-schémas. On en déduit qu'il correspond à une intégrale première rationnelle non constante de $(X,\delta_X)$.
\end{proof}

{\Thm\label{orthogonaliteconstantes} Soient $(K,\delta)$ un corps différentiel et $(X,\delta_X)$ une $D$-variété absolument irréductible au dessus de $(K,\delta)$. On a équivalence entre : 
\begin{itemize}
\item[(i)] Le type générique de $(X,\delta_X)$ est orthogonal aux constantes.
\item[(ii)] Pour toute extension de corps différentiels $(K,\delta) \subset (L,\delta_L)$, le changement de base $(X,\delta_X)_{(L,\delta_L)}$ est sans intégrale rationnelle non constante.
\item[(iii)]Pour tout $n \in \mathbb{N}$, la $D$-variété produit $(X,\delta_X)^n$ est sans intégrale première rationnelle non constante.
\end{itemize}}

\begin{proof}
$(i) \Longrightarrow (ii)$ : Supposons qu'il existe une extension de corps différentiels $(K,\delta) \subset (L,\delta_L)$ telle que le changement de base $(X,\delta_X)_{(L,\delta_L)}$ admette une intégrale première non constante $f \in (L(X),\delta_X)$.
On obtient un morphisme rationnel dominant de $D$-schémas au dessus de $(L,\delta_L) $ : 
$$f : (X,\delta_X)_{(L,\delta_L)} \dashrightarrow (\mathbb{A}^1,0) \times_{(L^\delta,0)} (L,\delta_L)$$
Les types génériques de $(X,\delta_X)_{(L,\delta_L)}$ et de $(\mathbb{A}^1,0) \times_{(L^\delta,0)} (L,\delta)$  ne sont pas faiblement orthogonaux d'après le lemme \ref{morphismerationnel}.

De plus, d'après la propriété (ii) du lemme \ref{typegenerique}, le type générique de $(X,\delta_X)_{(L,\delta_L)}$ est l'unique extension non-déviante à $L$ du type générique de $(X,\delta_X)$ et le type générique  de $(\mathbb{A}^1,0) \times_{(L^\delta,0)} (L,\delta)$ est l'unique extension non déviante à $L$ du type générique des constantes.
On en déduit que $p_{(X,\delta_X)}$ est non-orthogonal au type générique des constantes.
 
$(ii) \Longrightarrow (iii)$ Considérons $f$ une intégrale première rationnelle de $(X,\delta_X)^n$ avec $n \geq 1$  minimal :
$$f : (X,\delta_X)^n \dashrightarrow (\mathbb{A}^1,0) \times_{(k,0)} (K,\delta)$$
 
On note $(L,\delta_L)$ le corps des fractions de  la $D$-variété irréductible $(X,\delta_X)^{n-1}$. Par minimalité de $n \in \mathbb{N}$, la $D$-variété $(X,\delta_X)^{n-1}$ n'admet pas d'intégrales premières rationnelles non constantes et donc $L^{\delta_L} = k$. Par changement de base, $f$ induit une intégrale première

$$\tilde{f} : (X,\delta_X)_{(L,\delta_L)} \dashrightarrow (\mathbb{A}^1,0) \times_{(k,0)} (L,\delta_L)$$
 
qui est non constante car $f \notin k$ et $L^{\delta_L} = k$.
 
$(iii) \Longrightarrow (i)$: On note $q$ le type générique des constantes.

{\lem Soit $(K,\delta)$ un corps différentiel et $p \in S(K)$ un type non-orthogonal aux constantes. Il existe un entier $n \in  \mathbb{N}$, une réalisation $a$ de $p^{\otimes n}$ et une réalisation $b$ de $q|K$ telle que $b \in \mathrm{dcl}(a,K)$.}

\begin{proof}
D'après le lemme \ref{ortho2}, il existe des entiers $n,m \in \mathbb{N}$ tels que les types $p^{\otimes m}$ et $q^{\otimes n}$ ne sont pas faiblement orthogonaux. Considérons des réalisations $a \models p^{\otimes m}$ et $b \models q^{\otimes n}$ telles que $a \nind_K b$.
Il existe donc une formule $\phi(x,y)$ à paramètres dans $K$ telle que
$$b \in \phi(\U, a) \subset \mathcal C^n \text{ et }  \mathrm{RM}(\phi(\U,a)) < n.$$ 

Le sous-ensemble $\phi(\U, a) \subset \mathcal C^n$ est un sous-ensemble définissable à paramètres dans $(\U,\delta_U)$. D'après le corollaire \ref{corpsdesconstantes}, le corps des constantes est un pur corps algébriquement clos stablement plongé (et donc en particulier élimine les imaginaires dans le langage des anneaux). Il existe donc un paramètre canonique $c \in \mathcal C^l$  pour la formule $\phi(\U,a)$, c'est-à-dire que pour tout automorphisme de corps $\sigma \in \mathrm{Aut}(\U,\delta_\U)$, on a 
$$ \sigma(c) = c \text{ si et seulement si } \sigma(\phi(\U,a)) = \phi(\U,a).$$    

En particulier, on a $c \in \mathrm{dcl}(K,a)$. Montrons que $\mathrm{RU}(c/K) > 0$.
En effet, le rang de Lascar et le rang de Morley coïncident dans un corps algébriquement clos (example \ref{RU-ACF0}). On a alors :
$$ n = \mathrm{RU}(b/K) > \mathrm{RM}(\phi(x,a)) \geq  \mathrm{RM}(b/K,c) = \mathrm{RU}(b/K,c).$$   
On en déduit que $c \nind_K b$ et donc que le rang de Lascar de $c$ sur $K$ est $\geq 1$. Il suffit alors de remplacer $c$ par une coordonnée non algébrique de $c$.
\end{proof}
Supposons que $p = p_{(X,\delta_X)}$ est non-orthogonal aux constantes. D'après le lemme précédent, il existe  un entier $n \in  \mathbb{N}$, une réalisation $a$ de $p^{\otimes n}$ et une réalisation $b$ de $q|K$ telle que $b \in \mathrm{dcl}(a,K)$.

Comme $p^{\otimes n}$ est le type générique de $(X,\delta_X)^n$, il suffit d'appliquer le lemme \ref{morphismerationnel} pour obtenir une intégrale première non constante de $(X,\delta_X)^n$, ce qui montre par contraposition que $(iii) \Longrightarrow (i)$. 
\end{proof}

{\rem Modulo la proposition \ref{morphismerationnel}, l'équivalence entre (i) et (iii) se réduit à l'énoncé suivant (valable dans toute théorie stable) : Soient $T$ une théorie stable éliminant les imaginaires et $\Mod$ un modèle de $T$ $\kappa$-saturé.

Considérons $A$ un ensemble de paramètres petit, $D$ un ensemble fortement minimal $A$-définissable éliminant les imaginaires et $p \in S(A)$ un type stationnaire. Alors le type $p$ est non-orthogonal au type générique de $D$ si et seulement s'il existe $\overline{a} \in M^n$ réalisant  $p^{(n)}$ tel que
$\mathrm{dcl}(K) \cap D \subsetneq \mathrm{dcl}(K\overline{a}) \cap D.$}

\subsection{Spécialisation pour les intégrales premières rationnelles} Dans cette sous-section, on fixe $k$ un corps de caractéristique $0$.

{\defn Soit $S$ un schéma au dessus de $k$. On appelle \textit{famille de $D$-schémas paramétrée par $S$}, tout morphisme de $D$-schémas $f : (\mathcal X,v) \longrightarrow (S,0)$.

On dira que $f$ est \textit{une famille lisse de $D$-variétés paramétrée par $S$} lorsque $S$ et $\mathcal X$ sont des schémas intègres lisses, de type fini au dessus de $k$ et le morphisme $f$ est lisse.}

{\rem De façon équivalente, une famille lisse de $D$-variétés au dessus de $k$ est la donnée d'une famille lisse $f: \mathcal X \longrightarrow S$ de variétés au dessus de $k$ et d'un champ de vecteurs $v \in \mathrm{H}^0(\mathcal X , \Theta_{X/S})$ sur $\mathcal X$, tangent aux fibres de $f$.}

\nota Soit $f : (\mathcal X,v) \longrightarrow (S,0)$  une famille lisse de $D$-variétés paramétrée par $S$ et $s \in S$. On note $(\mathcal X,v)_s$, \textit{la fibre de $(\mathcal X,v)$ en $s$}.

Le champ de vecteurs $v$ étant une section du fibré tangent relatif de $f$, il se restreint un champ de vecteurs sur chacune des fibres de $f$. La fibre $(\mathcal X,v)_s$ est une variété lisse $X_s$ au dessus de $k$ munie d'un champ de vecteurs $v_s$.

Les fibres de $f$ seront toujours munies de cette structure de $D$-variétés. Le corps des intégrales premières rationnelles de $(\mathcal X,v)_s$ sera alors noté $ k(\mathcal X_s)^\delta$ (et donc sans référence explicite au champ de vecteurs sur $\mathcal X_s$).

{\Thm\label{specialisationtheoremerationnelle} Soit  $f : (\mathcal X,v) \longrightarrow (S,0)$  une famille lisse de $D$-variétés paramétrées par une variété lisse et irréductible $S$ au dessus de $k$.

Supposons, de plus, que toutes les fibres de $f$ soient absolument irréductibles. Pour tout point $s \in S(k)$, on a: 
$$ \mathrm{degtr}(k(\mathcal X)^\delta/k(S)) \geq 1 \Longrightarrow \mathrm{degtr}(k(\mathcal X_s)^\delta/k) \geq 1.$$}

{\rem Il est vraisemblable qu'une inégalité plus forte que l'implication du théorème \ref{specialisationtheoremerationnelle} soit en fait valide et que sous les hypothèses du théorème \ref{specialisationtheoremerationnelle}, on ait en réalité: 
\begin{equation*}
\mathrm{degtr}(k(\mathcal X_s)^\delta/k)  \geq \mathrm{degtr}(k(\mathcal X)^\delta/k(S))\qquad \qquad (*)
\end{equation*} 

Il semble en effet que les résultats de \cite{Julliard} --- dans le cas hamiltonien --- s'appuyant sur le lemme de Ziglin puisse être adaptés ici en unr preuve de l'inégalité $(\ast)$.} 

La suite de cette sous-section est consacrée à une preuve du théorème \ref{specialisationtheoremerationnelle}.

\subsubsection{Terme initial} Dans le cas où $S = \mathbb{A}^1$, la notion essentielle est celle de terme initial qui permet de restreindre une fonction rationnelle $f$ à un diviseur irréductible (possiblement inclus dans le lieu d'indétermination de $f$).

{\cons\label{terme-initial}  Soit $f : \mathcal X \longrightarrow \mathbb{A}^1$ une famille lisse de variétés paramétrée par $\mathbb{A}^1$ telle que $\mathcal X_0$ soit absolument irréductible.

Par hypothèse, la fibre $\mathcal X_0 = f^{-1}(0)$ est un diviseur irréductible de $\mathcal X$ et correspond à un anneau de valuation discrète $\mathcal O_0$ sur le corps $k(\mathcal X)$. On note $\mathrm{ord}_0$ la valuation discrète associée. Par hypothèse, la fonction $f$ est lisse et donc:
$$ \mathrm{ord}_0(f) = 1.$$

Considérons maintenant une fonction rationnelle $g \in k(\mathcal X)$ sur $\mathcal X$ et notons $n_g = \mathrm{ord}_0(g)$. La fonction rationnelle

$$g_0 =  g.f ^{-n_g} \in k(\mathcal X)\text{ vérifie donc }\mathrm{ord}_0(g_0) = 0.$$

On en déduit que la fonction $g_0$ vérifie $g_0 \in \mathcal O_Z \setminus m_Z$ et se restreint donc en une fonction rationnelle  non identiquement nulle sur $\mathcal X_0$.}

{\defn\label{terme-initialdef} Soient $f : \mathcal X \longrightarrow \mathbb{A}^1$ une famille lisse de variétés paramétrée par $\mathbb{A}^1$.
Pour toute fonction rationnelle $g \in k(\mathcal X)$ non identiquement nulle sur $\mathcal X$, on appelle \textit{terme initial de $g$ le long de $\mathcal X_0$}, la fonction $g_0 \in k(\mathcal X)$ associée à $g$ par la construction \ref{terme-initial}.} 

{\rem Les deux remarques suivantes découlent de la construction \ref{terme-initial}:
\begin{itemize}
\item Le terme initial $g_0$ de $g$ se restreint en une fonction non-identiquement nulle sur $\mathcal X_0$.
\item Si $f : (\mathcal X,v) \longrightarrow (\mathbb{A}^1,0)$ une famille lisse de $D$-variétés  à paramètres dans $\mathbb{A}^1$ et g est une intégrale première de $(\mathcal X,v)$ alors son terme initial $g_0 \in k(\mathcal X)$ est aussi une intégrale première de $(\mathcal X,v)$.
\end{itemize}}

\subsubsection{Familles de $D$-variétés paramétrée par $\mathbb{A}^1_k$}
Pour les fonctions rationnelles, l'indépendance algébrique et fonctionnelle coïncident:
 
{\rem[voir, par exemple, {\citep[Appendice III.6]{Aud}}]\label{independance algebrique}  Soit $X$ une variété algébrique au dessus de $k$ $f_1, \ldots , f_p \in k(X)$ des fonctions rationnelles. Les propriétés suivantes sont équivalentes:
\begin{itemize}
\item[(i)] Les fonctions $f_1, \ldots, f_p$ sont $k$-algébriquement indépendantes (indépendance algébrique).
\item[(ii)] La $p$-forme $df_1 \wedge \cdots df_p \neq 0$ est non nulle dans $\mathrm{H}^0(X, \Lambda^p \Omega^1_{X/k} \otimes k(X))$ (indépendance fonctionnelle).
\end{itemize}}

{\lem\label{specialisationtheoreme} Soient $k$ un corps de caractéristique $0$ et  $f : (\mathcal X,v) \longrightarrow (\mathbb{A}^1_k,0)$  une famille lisse de $D$-variétés à paramètres dans $\mathbb{A}^1_k$.

On suppose que les variétés algébriques $\mathcal X$ et $\mathcal X_0$ (la fibre en $0$) sont absolument irréductibles. Alors:
$$ \mathrm{degtr}(k(\mathcal X)^\delta/k(\mathbb{A}^1)) \geq 1 \Longrightarrow \mathrm{degtr}(k(\mathcal X_0)^\delta/k) \geq 1.$$}
\begin{proof}
Considérons $g \in k(\mathcal X)^\delta$ une intégrale première rationnelle algébriquement indépendante de $f$. Notons que:
\begin{itemize}
\item Quitte à réduire $\mathcal X$ par un ouvert affine rencontrant $\mathcal X_0$, on peut de plus supposer que $\mathcal X = \mathrm{Spec}(A)$ est affine.

\item Quitte à remplacer $g$ par son terme initial le long de $\mathcal X_0$ (en passant éventuellement à un nouvel ouvert affine), on peut supposer que $g$ est bien définie sur $\mathcal X$ et non-identiquement nulle sur $\mathcal X_0$.
\end{itemize}
Par hypothèse, $f$ est lisse et donc l'idéal $I$ définissant le sous-schéma fermé $X\mathcal X_0$ s'écrit  $I = (f)$.

On raisonne par l'absurde et on suppose que $\mathrm{degtr}(k(\mathcal X_0)^\delta/k) = 0$. Puisque $\mathcal X_0$ est absolument irréductible, cela implique que $k(\mathcal X_0)^\delta = k$ et donc que toute intégrale première rationnelle  de $(\mathcal X_0,v)$ est constante.

{\asser Pour tout $n \geq 1$, il existe $c_0, \ldots , c_{n-1} \in k$ tels que:
$$ g - \sum_{i = 0}^{n-1} c_k.f^k \in I^{n}.$$}

\begin{proof}[Preuve de l'assertion]
On raisonne par récurrence sur $n \in \mathbb{N}$:
\begin{itemize}
\item Pour $n = 1$, on sait que $g_{|\mathcal X_0}$ est une intégrale première de $(\mathcal X_0,v)$ et que $k(\mathcal X_0)^\delta = k$. Il existe donc une constante $c_0 \in k$ telle que $g - c_0 \in I$.

\item Supposons qu'il existe  $c_0, \ldots , c_{n-1} \in k$ tel que $g  = \sum_{k = 0}^{n-1} c_k.f^k + h.f^n$ avec $h \in A$.

Puisque $h \in k(g,f)$, c'est une intégrale première de $\mathcal X$ et en appliquant le cas $n = 1$, il existe $c_n \in k$ tel que $h = c_n + f.h_2$ avec $h_2 \in A$. On en déduit que
$g - \sum_{k = 0}^{n} c_k.f^k = f^{n+1}.h_2 \in I^{n+1}$.
\end{itemize}
\end{proof}

Considérons maintenant l'élément $df \wedge dg$ du $A$-module $\Lambda^2 \Omega^1_{A/k}$. Pour tout $n \geq 1$, on a

$$ df \wedge dg = df \wedge d(g - \sum_{k = 0}^{n-1} c_k.f^k) \in \Omega^1_{A/k} \wedge d(I^{n}) \subset I^{n-1}. \Lambda^2 \Omega^1_{A/k}.$$

Cela implique donc que:
$$df \wedge dg \in \bigcap_{n \in \mathbb{N}} I^n. \Omega^1_{A/k} = \lbrace 0 \rbrace.$$

En utilisant la remarque \ref{independance algebrique}, cela implique que $f$ et $g$ ne sont pas algébriquement indépendantes, ce qui contredit le choix de $g$.
\end{proof}

\subsubsection{Preuve du théorème \ref{specialisationtheoremerationnelle}} On commence par étendre le lemme \ref{specialisationtheoreme} au cas des familles paramétrées par $\mathbb{A}^n_k$ en raisonnant par induction. Pour cela on considère le drapeau:

$$ \lbrace 0 \rbrace  \subset \lbrace 0 \rbrace \times \mathbb{A}^1_k \subset \cdots \subset \lbrace 0 \rbrace \times \mathbb{A}^{n-1}_k \subset \mathbb{A}^{n}_k.$$

{\lem\label{avaleursdansAn} Soient $k$ un corps de caractéristique $0$ et  $f : (\mathcal X,v) \longrightarrow (\mathbb{A}^n_k,0)$  une famille lisse de $D$-variétés paramétrée par $\mathbb{A}^n$ avec $n > 0$.

Pour $r \leq n$, on note $\eta_r$ le point générique de $\lbrace 0 \rbrace \times \mathbb{A}^r_k$ et on suppose que les fibres de $f$ en $0 = \eta_0 \ldots, \eta_n$ sont absolument irréductibles (et donc, en particulier, non vides). Alors:
$$ \mathrm{degtr}(k(\mathcal X)^\delta/k(\mathbb{A}^n) \geq 1 \Longrightarrow \mathrm{degtr}(k(\mathcal X_0)^\delta/k) \geq 1.$$}
\begin{proof}
On raisonne par récurrence sur $n \in \mathbb{N}^\ast$
\begin{itemize}
\item Pour $n = 1$, c'est le contenu du lemme \ref{specialisationtheoreme}.
\item Supposons maintenant le résultat établi pour un $n > 0$ et considérons $f : (\mathcal X,v) \longrightarrow (\mathbb{A}^{n+1},0)$  une famille lisse de $D$-variétés à paramètres dans $\mathbb{A}^{n +1}$ dont toutes les fibres sont absolument irréductibles.
\end{itemize}

Considérons l'hyperplan $H = \lbrace 0 \rbrace \times \mathbb{A}^n$ et la restriction de la famille
$$f_H : f^{-1}(H) \longrightarrow H$$
à cet hyperplan. La famille $f_H$ paramétrée par $H$ vérifie les hypothèses du lemme \ref{avaleursdansAn}. En appliquant l'hypothèse de récurrence à cette famille, il suffit de vérifier que:
\begin{eqnarray} \label{equation} \mathrm{degtr}(k(\mathcal X)^\delta/k(\mathbb{A}^{n+1}) \geq 1  \Longrightarrow \mathrm{degtr}(k(f^{-1}(H))^\delta/k(\mathbb{A}^{n})) \geq 1.
\end{eqnarray} 

\textit{Preuve de l'implication (\ref{equation}):}
Notons $\pi: \mathbb{A}^{n+1} \longrightarrow \mathbb{A}^n$ la projection sur les $n$ dernières coordonnées  et considérons 
$$ (\mathcal X,v) \overset f \longrightarrow (\mathbb{A}^{n+1},0) \overset \pi \longrightarrow (\mathbb{A}^n,0).$$

Posons $K = k(x_1, \ldots , x_n)$. Après changement de base par le point générique de $\mathbb{A}^n$, on obtient une famille de $D$-variétés à paramètres dans $\mathbb{A}^1$ au dessus de $K$
$$ f_K: \mathcal (X,v)_{(K,0)} \longrightarrow (\mathbb{A}^1_K, 0).$$
L'implication (\ref{equation}) est alors conséquence de la proposition \ref{specialisationtheoreme} appliquée à cette famille  de  $D$-variétés à paramètres dans $\mathbb{A}^1$ au dessus de $K$. 
\end{proof}

\begin{proof}[Preuve du théorème \ref{specialisationtheoremerationnelle}] On se ramène au cas des familles  de $D$-variétés paramétrée par $\mathbb{A}^n_k$ à l'aide de coordonnées étales.

Soient $S$ une variété algébrique lisse et irréductible au dessus d'un corps $k$ et  $f : (\mathcal X,v) \longrightarrow (S,0)$  une famille lisse de $D$-variétés à paramètres dans $S$ dont toutes les fibres sont absolument irréductibles.

Notons $n = \mathrm{dim}(S)$ et fixons $s \in S(k)$. Puisque $S$ est supposée lisse, il existe un voisinage $U$ de $s$ dans $S$ et un morphisme étale $g : U \longrightarrow \mathbb{A}^{n}$ tel $g(s) = 0$. On considère la composition
$$ \tilde{f} : (f^{-1}(V),v) \overset f \longrightarrow (V,0) \overset g \longrightarrow  (\mathbb{A}^{n},0).$$
qui est une famille lisse de $D$-variétés paramétrée par $\mathbb{A}^n$. 

On modifie maintenant la famille $\tilde{f}$ de sorte que les fibres au dessus de $0 = \eta_0 \ldots, \eta_n$ soient absolument irreductible.

Considérons l'hyperplan $H = \lbrace 0 \rbrace \times \mathbb{A}^{n-1}$. On peut alors écrire une décomposition en composantes irréductibles:
$$g^{-1}(H) = Z_1 \cup \ldots \cup Z_n.$$

Comme le morphisme $g$ est étale, toutes les composantes irréductibles de $g^{-1}(H)$ ont même dimension et pour tout point $p \in H$ tel que $g^{-1}(p) \in H$, il existe une unique composante irréductible de $H$ contenant $p$. En particulier:
\begin{itemize}
\item Il existe une unique composante irréductible qui contient $s$. Sans perte de généralité, on peut supposer que $s \in Z_1$.

\item Pour tout $r \leq n$, le point $\eta_r$ admet un unique antécédent noté $\gamma_r$ dans $Z_1$.  
\end{itemize}

On pose alors $V = U \setminus \lbrace Z_2 \ldots Z_n \rbrace$ et on considère la composition:
$$ \tilde{f}_V : (f^{-1}(V),v) \overset f \longrightarrow (V,0) \overset g \longrightarrow  (\mathbb{A}^{n},0).$$ 

Avec les notations précédentes, le morphisme $\tilde{f}_V$ satisfait:
\begin{itemize}
\item[(i)] Le morphisme $\tilde{f}_V$ définit une famille lisse de $D$-variétés paramétrée par $\mathbb{A}^n_k$.
\item[(ii)] La fibre $(\tilde{f}^{-1}_{|V}(\eta_r), v_{|\tilde{f}^{-1}_{|V}(\eta_r)})$ est isomorphe (en tant que $D$-variété) à $(\mathcal X,v)_{\gamma_r}$. En particulier, cette fibre est absolument irréductible.
\end{itemize}
Il suffit alors d'appliquer le lemme \ref{specialisationtheoreme} à la famille $\tilde{f}_V$ à paramètres dans $\mathbb A^n_k$. 
\end{proof}

\subsection{Spécialisation et non-orthogonalité aux constantes} 

\subsubsection{Orthogonalité aux constantes en famille} Fixons $T$ une theorie stable avec l'élimination des imaginaires.

{\defn Soit $r(\overline{y})$ un type partiel à paramètres dans $A$. On appelle \textit{famille de types stationnaires à paramètres dans $r(\overline{y})$}, tout type partiel $\pi(\overline{x},\overline{y})$ à paramètres dans $A$ vérifiant: pour tout $\overline{a} \models r(\overline{y})$, le type $\pi(\overline{x},\overline{a}) \in S(A\overline{a})$ est un type (complet) stationnaire.}

{\rem Soit $f : (X,v) \longrightarrow (S,0)$ une famille de $D$-variétés absolument irréductibles au dessus d'un corps différentiel constant $(k,0)$.
Le type partiel (à paramètres dans $k$):
$$\pi(x,y) := \lbrace y \in S \rbrace  \cup \lbrace x \text{ realise le type generique de } (X,v)_y \text{ au dessus de }k(y) \rbrace$$
est une famille de types stationnaires à paramètres dans $\lbrace y \in S \rbrace$.}    

{\cor\label{typedefinissable} Soient $r(\overline{y}) \in S(A)$ un type partiel et $\pi_1(\overline{x},\overline{y}),\pi_2(\overline{x},\overline{y}) \in S(A)$ deux familles de types stationnaires à paramètres dans $r(\overline{y})$. Pour toute réalisation $\overline a \models r(y)$, la propriété:
$$ \mathcal P(\overline a) :  \pi_1(\overline{x},\overline{a})\text{ et } \pi_2(\overline{x},\overline{a}) \text{ sont orthogonaux}$$
ne dépend que du type $\mathrm{tp}(\overline a/A)$ de $\overline a$ sur $A$.}

\begin{proof}
En utilisant la proposition \ref{ortho2}, il suffit de démontrer que la propriété:
$$ \mathcal P'(a) :\pi_1(\overline{x},\overline{a})\text{ et } \pi_2(\overline{x},\overline{a}) \text{ sont faiblement  orthogonaux}$$
ne dépend que du type $\mathrm{tp}(\overline{a}/A)$.

Soient $\overline{a}$ et $\overline{a}'$ deux réalisations du type $\mathrm{tp}(\overline{a}/A)$.

Supposons que la propriété $\mathcal P(\overline a)$ ne soit pas vérifiée. Il existe une formule $\phi(x,y,\overline{a})$ qui dévie au dessus de $A$ et telle que:
$$ \pi(x,y,\overline{a}) = p(x,\overline a) \cup q(y,\overline a) \cup \phi(x,y,\overline{a}) \text{ est coherent. }$$

Par invariance par automorphisme de la déviation, on en déduit que $\pi(x,y,\overline a')$ est cohérent  et que $\phi(x,y,\overline{a}')$ dévie au dessus de $A$. Il suit que $\mathcal P(\overline a')$ n'est pas vérifiée.
\end{proof}

{\rem  Soit $f : (\mathcal X,v) \longrightarrow (S,0)$ une famille de $D$-variétés absolument irréductibles au dessus d'un corps différentiel constant $(k,0)$. Considérons l'ensemble:
$$ S^{\not\perp 0} = \lbrace s \in S(\overline l) \text{ } |  \text{ le type générique de } (\mathcal X,v)_s \text{ est non-orthogonal aux constantes} \rbrace$$     
où $\overline{l}$ est une extension algébriquement close saturée, fixée de $k$. 

Le corollaire \ref{typedefinissable} implique ici que $S^{\not \perp 0}$ s'écrit comme l'image inverse d'un sous-ensemble du schéma $S$ par l'application 
$$ \pi:  S(\overline{k})  \longrightarrow S$$
qui envoie un point de $S(\overline{k})$ sur son image dans le schéma $S$.}

En revanche, on sait que en général l'ensemble $S^{\not\perp 0}$ n'est pas un sous-ensemble constructible (ni même un sous ensemble $\infty$-definissable de $S(\overline{l})$). Cela découle de l'étude de Hrushovski et Itai des champs de vecteurs sur $\mathbb{A}^1$ orthgonaux aux constantes (voir \cite{Itai}, Exemple 2.20 ainsi que l'exemple \ref{Hrushovski-Itai} dans la troisième section de cette partie).

De même, on note $S^{\perp 0}$, le complémentaire de $S^{\not\perp 0}$ dans $S(\overline{l})$ dont les éléments sont les points $s \in S(\overline{l})$ tels que $(\mathcal X,v)_s$ est orthogonal aux constantes.

{\lem\label{compacitytrick} Soit $S$ une variété irréductible algébrique complexe et $f: (\mathcal X,v) \longrightarrow (S,0)$ une famille d'équations différentielles algébriques paramétrée par $S$ dont les fibres sont toutes absolument irréductibles.

Si la fibre générique de $f$ est orthogonale aux constantes alors, il existe un ensemble dénombrable de sous-variétés algébriques fermées propres $\lbrace Z_n \text{ | } n \in \mathbb{N}\rbrace$  de $S$ tel que:
$$ S(\mathbb{C}) \setminus \bigcup_{i \in \mathbb{N}} Z_i (\mathbb{C}) \subset  S^{\perp 0}(\C).$$}   

\begin{proof}
Il existe un corps $k$ finement engendré sur $\mathbb{Q}$ (donc dénombrable) sur lequel la variété $S$, le morphisme $f$ et le champ de vecteurs $v$ sont définis. Le morphisme $f$ est donc obtenu par changement de base à partir d'un morphisme $f_k : (\mathcal X_k,v_k) \longrightarrow S_k$ défini au dessus de $k$. Par hypothèse, il existe une réalisation $s$ du point générique de $S_k$ telle que $(\mathcal X,v)_s$ est orthogonal aux constantes. 

Le corollaire \ref{typedefinissable} assure que c'est le cas pour toutes les réalisations dans $S(\mathbb{C})$ du point générique de $S_k$. Il suffit donc de choisir pour $\lbrace Z_n \text{ | } n \in \mathbb{N} \rbrace$, l'ensemble des sous-variétés fermées propres de $S$ définies sur $k$. Comme $k$ est dénombrable, cet ensemble est bien dénombrable.
\end{proof}

%

\subsubsection{Un énoncé de spécialisation} A l'aide du premier critère d'orthogonalité aux constantes (Théorème \ref{orthogonaliteconstantes}), nous montrons maintenant que $S^{\perp 0}$ est stable par spécialisation pour les familles lisses de $D$-variétés dont les fibres sont absolument irréductibles.  

{\Thm[Spécialisation et non-orthogonalité aux constantes]\label{specialisationtheoreme2} Soient $S$ une variété algébrique lisse et irréductible au dessus d'un corps $k$ et  $f : (\mathcal X,v) \longrightarrow (S,0)$ une famille lisse de $D$-variété à paramètres dans $S$.

On suppose que toutes les fibres de $f$ sont absolument irréductibles, on fixe $s \in S(k)$ et on dénote par $\eta$ le point générique de $S$. Si le type générique de  $(\mathcal X,v)_\eta$ est non-orthogonal aux constantes alors le type générique de $(\mathcal X,v)_s$ est non-orthogonal aux constantes.}

\begin{proof}
Considérons $f : (\mathcal X,v) \longrightarrow (S,0)$ une famille lisse de $D$-variété à paramètres dans $S$ dont les fibres sont absolument irréductibles. Posons:
$$(\mathcal X_n, v_n) = \mathcal (X,v) \times_{(S,0)} \cdots \times_{(S,0)}  \mathcal (X,v)  \text{ and } 
 f_n : (\mathcal X_n, v_n) \longrightarrow (S,0).$$

Les hypothèses sur la famille $f$ impliquent que $f_n : \mathcal X_n \longrightarrow S$ est une famille lisse de $D$-variétés à paramètres dans $S$ dont les fibres sont absolument irréductibles. De plus, si $s \in S$, la fibre de $f_n$ au dessus de $s$ est la puissance $n$-ieme $(\mathcal X,v)_{s}^n$ de la fibre de $f$ en $s$.

Supposons maintenant que le type générique de $(\mathcal X,v)_\eta$ est non-orthogonal aux constantes. D'après le théorème \ref{orthogonaliteconstantes}, il existe un entier $n > 0$ telle que $(\mathcal X,v)^n_\eta$ admet une intégrale première rationnelle non constante. On en déduit que: 
$$ \mathrm{degtr}(k(\mathcal X_n)^\delta/k) \geq \mathrm{dim}(S) + 1 .$$

Pour tout $s \in S$, la proposition \ref{specialisationtheoreme} implique alors que 
$$ \mathrm{degtr}(k(\mathcal X_{n, s})^\delta/k) \geq 1.$$ 

Cela implique que $(\mathcal X,v)_{s}^n$ admet une intégrale première non-constante et donc que le type générique de $(\mathcal X,v)_{s}$ est non-orthogonal aux constantes
\end{proof}

{\cor\label{specialisationgeometrique} Soient $S$ une variété algébrique lisse et irréductible au dessus du corps $\mathbb{C}$ des nombres complexes et  $f : (\mathcal X,v) \longrightarrow (S,0)$ une famille lisse de $D$-variété à paramètres dans $S$.

On suppose que toutes les fibres de $f$ sont absolument irréductibles. Exactement  un des deux cas suivants est réalisé:
\begin{itemize}
\item[(i)] $\forall s \in S(\mathbb{C})$, le type générique de $(X,v)_s$ est  non-orthogonal aux constantes.
\item[(ii)] Le sous-ensemble $S^{\not \perp 0}(\mathbb{C}) \subset S(\mathbb{C})$ des paramètres complexes $s \in S(\mathbb{C})$ tels que $(X,v)_s$ est non-orthogonal aux constantes s'écrit: 
$$ S^{\not \perp 0}(\mathbb{C}) = \bigcup_{n \in \mathbb{N}} Z_n(\mathbb{C})$$
où $\lbrace Z_n \text{ } | n \in \mathbb{N} \rbrace$ est un ensemble dénombrable de sous-varietés fermées (irréductibles) propres de $S$. 
\end{itemize}}

\begin{proof}
On peut toujours supposer que la famille $f : (\mathcal X,v) \longrightarrow (S,0)$ est définie au dessus d'un sous-corps dénombrable $k$ du corps des nombres complexes. Pour  toute sous-variété fermée irréductible $Z$ de $S$, note $\eta_Z$, le point générique de $Z$ (au dessus de $k$). 

Considérons l'ensemble dénombrable $\mathcal E$ des sous-variétés fermées irréductibles $Z$ de $S$ (au dessus de $k$) telles que l'équation différentielle $(X,v)_{\eta_Z}$ est non-orthogonal aux constantes.

Supposons que le cas (i) n'est pas réalisé. D'après le théorème \ref{specialisationtheoreme2}, l'équation différentielle $(X,v)_{\eta_S}$ est orthogonale aux constantes et donc $S \notin \mathcal E$. Montrons que: 
$$ S^{\not \perp 0}(\C) = \bigcup_{Z \in \mathcal E} Z (\mathbb{C}).$$

Considérons $Z \in \mathcal E$  et $\pi: \widehat{Z} \rightarrow Z$ une désingularisation de $Z$ (qui existe car $k$ est un corps de caractéristique $0$). Notons de plus $\widehat{f_Z} : \mathcal X_Z \longrightarrow \widehat{Z}$ la famille d'équations différentielles induite par changement de base de la restriction $f_{|f^{-1}(Z)}$  par $\pi$.

Puisque la lissité et l'irréductibilité des fibres d'un morphisme $f: X \longrightarrow S$ sont deux propriétés invariantes par restriction et par changement de base, on en déduit que  $\widehat{f_Z} : \mathcal X_Z \longrightarrow \widehat{Z}$ vérifie les hypothèses du théorème \ref{specialisationtheoreme2}.

De plus, comme $\pi^{-1}(\eta_Z) = \eta_{\widehat{Z}}$, la fibre générique de $f_{\widehat{Z}}$ est isomorphe à $(\mathcal X,v)_{\eta_Z}$ et donc non-orthogonale aux constantes. D'après le théorème \ref{specialisationtheoreme2}, pour tout point $s \in \widehat{Z}(\mathbb{C})$, l'équation différentielle $(\mathcal X_Z,v)_s$ est orthogonale aux constantes. On en déduit que pour tout $z \in Z(\mathbb{C})$, l'équation différentielle $(\mathcal X,v)_z$ est orthogonale aux constantes et donc que:
$$ \bigcup_{Z \in \mathcal E} Z (\mathbb{C}) \subset S^{\not \perp 0}(\C).$$

Réciproquement, si $s \in  S^{\not \perp 0}(\C)$ alors $Z = \mathrm{loc}_k(s)$ est une sous-variété fermée de $S$ au dessus de $k$ et l'équation différentielle $(\mathcal X,v)_{\eta_Z}$ est orthogonale aux constantes (car $(X,v)_s$ est orthogonale aux constantes et $s$ réalise le point générique de $Z$). On conclut que $Z \in \mathcal E$ et donc que $s \in \bigcup_{Z \in \mathcal E} Z (\mathbb{C})$.
\end{proof}
\subsection{Champ de vecteurs polynomiaux très génériques} On applique maintenant les résultats précédents a l'étude des champs de vecteurs polynomiaux  complexes sur l'espace affine complexe $\mathbb{A}^n_\mathbb{C}$ dimension $n$. Le résultat principal de cette section  concerne le cas très générique où les coefficients du champs de vecteurs sont $\mathbb{Q}$-algébriquement indépendants.

{\cor\label{caseofAn} Soit $d \geq 3$ et $n \geq 1$. Considérons un champs de vecteurs 
$$v(x_1,\ldots , x_n) = f_1(x_1,\ldots, x_n) \frac {d} {dx_1} + \cdots f_n(x_1,\ldots, x_n) \frac {d} {dx_n}.$$ 
sur l'espace affine complexe de dimension $n$, où $f_1,\ldots, f_n \in K[X_1,\ldots, X_n]_{\leq d}$ sont des polynômes de degré $\leq d$.

Si les coefficients de $f_1, \ldots, f_n$ sont $\mathbb{Q}$-algébriquement indépendants\footnote{Cela implique, en particulier, que les $f_i$ sont des polynômes de degré $d$ dont tous les coefficients sont distincts et non nuls.} alors le champ de vecteurs $v$ est orthogonal aux constantes.}

Formulons aussi la version géométrique suivante qui se déduit immédiatement à l'aide du lemme \ref{compacitytrick}. 

{\cor\label{geometriquecaseofAn} Soit $d \geq 3$ et $n \geq 1$. Notons $\mathcal A_{n,d}$ l'ensemble des champs de vecteurs complexes de degré $d$ sur $\mathbb{A}^n_\mathbb{C}$, qui est un $\mathbb{C}$-espace vectoriel de dimension $n^{d+1}$.

Il existe un ensemble dénombrable $\lbrace Z_n \text{ } | n \in \mathbb{N} \rbrace$ de sous-variétés fermées propres de $\mathcal A_{n,d}$ tel que:
$$ \forall v \in \mathcal A_{n,d} \setminus \bigcup_{i \in \mathbb{N}} Z_i(\mathbb{C}) \text{ , le type generique de } (\mathbb{A}^n_\mathbb{C},v) \text{ est orthogonal aux constantes}.$$}

\subsubsection{Cas de la droite affine} En dimension $1$, des travaux de Rosenlicht (voir \cite{Rosenlicht}) décrivent très precisement les champs de vecteurs polynomiaux orthogonaux aux constantes. Le corollaire \ref{caseofAn} est donc en dimension $1$, une conséquence directe \cite{Rosenlicht}. La formulation que nous donnons ici est celle de \cite{Itai} (voir Exemple 2.20).

{\Thm[Rosenlicht] Soit $v(x) = P(x) \frac d {dx}$ un champ de vecteurs polynomial sur la droite affine. Le type générique de $(\mathbb{A}^1,v)$ est orthogonal aux constantes si et seulement si l'un des deux cas suivants est réalisé:
\begin{itemize}
\item[(i)] soit le polynôme $P(x)$ admet au moins une racine multiple et une racine simple.
\item[(ii)] soit toutes les racines de $P(x)$ sont simples et il existe deux racines $x_1$ et $x_2$ de $P(x)$ telles que les résidus de $\frac 1 {P(x)}$ sont $\mathbb{Q}$-linéairement indépendants. 
\end{itemize}}

{\exam Considérons la famille $\mathcal A_{1, d}^\ast$ des champs de vecteurs polynomiaux unitaires de degré $d$ complexes, de dimension $1$.  On a donc $\mathcal A_{1, d}^\ast \simeq \mathbb{C}^d$.

A l'aide du théorème de Rosenlicht, on peut décrire ici explicitement les sous-variétés fermées $\lbrace Z_n \text{ | } n \in \mathbb{N} \rbrace$ apparaissant dans le corollaire \ref{geometriquecaseofAn}:
\begin{itemize}
\item Le premier cas est contenu dans la sous-variété fermée propre $Z \subset \mathcal A^\ast_{1,d}$  décrite par l'équation $$\mathrm{Res}(P,P') = 0.$$

\item Notons $\Delta \subset \mathbb{C}^d$ l'ensemble des diagonales de $\mathbb{C}^d$. Le complémentaire $\mathcal A_{1,d}^\ast \setminus Z$ peut alors être identifié au quotient $(\mathbb{C}^d \setminus \Delta) / \Sigma_d$ via:
$$\pi: \begin{cases}
\mathbb{C}^d \setminus \Delta \rightarrow  \mathcal A_{1,d}^\ast \setminus Z \\
(x_1, \ldots , x_d) \mapsto (x - x_1) \cdots (x - x_d) \frac d {dx}
\end{cases}$$
\end{itemize}
Notons que en particulier $\pi$ est fermée. Pour tout choix $i \neq j \leq n$ et $q \in \mathbb{Q}$, l'équation
$$  \mathrm{Res}_{x_i} (1/ P(x)) =  q. \mathrm{Res}_{x_j} (1/ P(x))$$
définit une sous-variété algébrique fermée $Y_{i,j,q}$ (au dessus de $\mathbb{Q}$) de $\mathbb{C}^n \setminus \Delta$. Un calcul immédiat montre que cette sous-variété fermée est une \textit{sous-variété propre} dès que le degré $d \geq 3$.

Supposons donc $d \geq 3$. En posant $Z_{i,j,q} = \pi(Y_{i,j,q})$ (qui est une sous-variété fermée propre de $\mathcal A^\ast_{1,d} \setminus Z$), le théorème de Rosenlicht implique donc que:

$$ \forall v \in \mathcal A^\ast_{1,d} \setminus \Big( \bigcup_{i \neq j, q \in \mathbb{Q}} Z_{i,j,q}(\mathbb{C}) \cup Z(\mathbb{C}) \Big) \text{, le type generique de } (\mathbb{A}^1_\mathbb{C},v) \text{ est orthogonal aux constantes}.$$}

\subsubsection{Démonstration du corollaire \ref{caseofAn}} La démonstration du corollaire \ref{caseofAn} s'appuie sur le théorème \ref{specialisationtheoreme2} plutôt que sur le théorème de Rosenlicht (valable uniquement sur $\mathbb{A}^1$).

{\asser Il existe un champ de vecteurs  $v_d$ ``universel de degré $d$''  sur $\mathcal X = \mathbb{A}^n_\mathbb{C} \times   \mathcal A_{n,d}$ défini sur $\mathbb{Q}$ et tangent à la seconde projection $\pi: (\mathcal X,v_d) \longrightarrow (\mathcal A_{n,d},0)$ tel que pour tout champ de vecteurs $w$ de degré $d$ sur $\mathbb{A}^n$, on ait:
$$ (\mathcal X , v_d)_w \simeq (\mathbb{A}^n,w).$$}
\begin{proof}
Il suffit de poser $\mathcal X = \mathbb{A}^n_\mathbb{C} \times   \mathcal A_{n,d}$ et de considérer la section $v_d$ du faisceau coherent $\pi_1^\ast \Theta_{\mathbb{A}^n/k}$ définie par
$ v_d(x,v) = v(x).$
%
\end{proof}

{\asser Il existe un champ de vecteurs $w$ de degré $3$ sur $\mathbb{A}^n$ tel que le type générique de  $(\mathbb{A}^n,w)$ est orthogonal aux constantes.}

\begin{proof}
Considérons le champ de vecteurs $v(x) = x^2(x-1)\frac d {dx}$ de degré $3$ sur $\mathbb{A}^1$.

Le théorème de Rosenlicht implique que le type générique de $(\mathbb{A}^1,v)$ est orthogonal aux constantes.  Cela implique que pour tout $n \in \mathbb{N}$, le type générique de $(\mathbb{A}^1,v)^n = (\mathbb{A}^n, v \times v \ldots \times v)$ est orthogonal aux constantes.

Le champ de vecteurs
$$ w = (v \times v \ldots \times v) (x_1,\ldots , x_n) = x_1^2(x_1-1)\frac \partial {\partial x_1} + \cdots + x_n^2(x_n-1)\frac \partial {\partial x_n}$$
est donc un champ de vecteurs de degré 3 sur $\mathbb{A}^n$ vérifiant les conclusions de l'assertion.
\end{proof}

\begin{proof}[Démonstration du corollaire \ref{caseofAn}:]
La famille $\pi: (\mathcal X,v_d) \longrightarrow (\mathcal A_{n,d},0)$ vérifie les hypothèses du théorème \ref{specialisationtheoreme2}. En effet, le morphisme $\pi$ est clairement lisse et les fibres de $\pi$ sont absolument irréductibles (car isomorphes à $\mathbb{A}^n$).

De plus, l'assertion précédente assure que le type générique de la fibre de $\pi$ au dessus de $w$
$$(\mathcal X, v_d)_w \simeq (\mathbb{A}^n, w)$$
est orthogonal aux constantes.

D'après le théorème \ref{specialisationtheoreme2}, la fibre generique de $\pi$ est orthogonal aux constantes. On en déduit qu'un champ de vecteurs dont les coefficients réalisent le type générique de $\mathcal A_{n,d}$ --- c'est à dire dont les coefficients sont $\mathbb{Q}$-algébriquement indépendants --- est orthogonal aux constantes. 
\end{proof}

\section{Propriétés dynamiques du flot réel analytique d'une $D$-variété réelle}

On s'intéresse désormais aux $D$-variétés $(X,v_X)$ définies sur le corps $\R$ des nombres réels muni de la dérivation triviale. On dispose d'un foncteur d'analytification réel vers la catégorie des espaces analytiques réels munis de champs de vecteurs. Sous des hypothèses de lissité, on peut alors, par les théorèmes classiques d'existence et d'unicité des solutions d'équations différentielles analytiques, intégrer ce champ de vecteur et obtenir un flot réel analytique.
Lorsque ce flot est complet (ce qui est automatique sous des hypothèses de compacité), il définit une action continue du groupe additif $(\R,+)$ sur l'espace topologique métrisable $X(\R)^{an}$.

Dans le première section, on étudie le foncteur d'analytification réel (ainsi que son analogue complexe) et son effet sur les sous-variétés fermées invariantes. Dans la deuxième section, on donne une présentation auto-contenue des résultats de dynamique topologique que nous utiliserons, et en particulier la notion de \textit{flot faiblement topologiquement mélangeant}.

\subsection{Des $D$-variétés réelles aux flots réels} On présente les foncteurs d'analytification réel (resp. complexe) pour les $D$-variétés réelles (resp. complexe). Sous des hypothèses de lissité, on construit alors le flot réel associé (resp. le flot complexe associé). On supposera que le lecteur est familier avec les résultats élémentaires de la théorie des équations différentielles (voir \cite{Arn} pour le cas réel et la partie 1 du chapitre 1 de \cite{Ily} pour le cas complexe).

\textit{Contrairement au cas algébrique, on suivra la terminologie analytique usuelle c'est-à-dire que les variétés analytiques seront toujours supposées lisses et on parlera d'espaces analytiques (réels ou complexes) lorsque l'on travaillera sans hypothèse de lissité.}

\subsubsection{Flot réel d'une $D$-variété réelle} On présente d'abord le foncteur d'analytification  pour les $D$-variétés réelles.

{\defn On appellera \textit{$D$-espace analytique réel}, tout couple $(M,v_M)$ où $M$ est un espace analytique réel et $v_M$ un champ de vecteurs analytique sur $M$. Si l'espace analytique réel sous-jacent $M$ est lisse, on dit que $(M,v_M)$ est \textit{une $D$-variété analytique réelle}}

Si $(M,v_M)$ et $(N,v_N)$ sont deux espaces analytiques réels, on définira les morphismes de $(M,v_M)$ vers $(N,v_N)$ comme les applications analytiques $f : M \longrightarrow N$ vérifiant $df(v_M) = f^\ast v_N$.

.

{\cons\label{analytification} Soit $(X,v_X)$ une $D$-variété au dessus de $(\mathbb{R},0)$. On munit l'ensemble des points réels  $M = X(\mathbb{R})$ de $X$ de sa structure analytique $(M,\mathcal O_M)$, i.e. de la topologie  métrisable d'espace analytique et du faisceau des fonctions analytiques $\mathcal O_M$.

{\rem L'espace analytique réel $(M,\mathcal O_M)$ est lisse si et seulement si $X(\R)$ est contenu dans l'ouvert $X_{reg}$ des points régulier de $X$.}
\smallskip

Le champ de vecteurs $v_X$ induit alors un champ de vecteurs analytique sur $M$ noté $v_M$. On a ainsi construit un $D$-espace analytique réel $(M,v_M)$.

La correspondance $(X,v_X) \mapsto (M,v_M)$  détermine un foncteur appelé \textit{foncteur d'analytification réel} de la catégorie des $D$-variétés au dessus de $(\mathbb{R},0)$ vers la catégorie des $D$-espaces analytiques réels.}

{\cons\label{flot} Soit $(M,v_M)$ une $D$-variété réelle analytique. Le théorème d'existence et d'unicité de solutions pour les équations différentielles analytiques (voir \citep[$\S 35$]{Arn}) assure l'existence pour tout $x \in M$, d'un unique ouvert connexe $U_x \subset \R$ maximal et d'une application analytique $\phi_x : U_x \longrightarrow M$ telle que 
 $$\begin{cases}
  \phi_x(0) = x   \\
\forall t \in U_x \text{ , } \frac d {dt} \phi_x(t) = v_M(\phi_x(t)). \end{cases}$$

De plus, l'ensemble $U = \bigsqcup_{x \in M}  U_x \subset \R \times M$ est un voisinage connexe de $\lbrace 0 \rbrace \times M$ et la collection des $(\phi_x)_{x \in M}$ définit une application analytique $\phi : U \longrightarrow M$ vérifiant :

 $$\begin{cases}
\forall x\in M \text{ , }  \phi(0,x) = x    \\
\forall t \in U_x \text{ , } \frac d {dt} \phi (t,x) = v_M(\phi(t,x)). \end{cases}$$}

{\defn Soit $(M,v_M)$ une $D$-variété analytique réelle. On appelle \textit{flot réel (maximal) de $(M,v_M)$}, le couple $(U,\phi)$ maximal associé par la construction \ref{flot}.}

On dira de plus que le flot $(U,\phi)$ est \textit{complet} si $U = \R \times M$.

{\lem[{\citep[$\S 35$]{Arn}}]\label{flotcomplet} Soit $(M,v_M)$ une $D$-variété analytique réelle compacte. Le flot réel de $(M,v_M)$ est complet.}
 
{\rem Soit $(M,v_M)$ une $D$-variété analytique réelle et $(U,\phi)$ le flot réel associé. On appelle \textit{courbe intégrale ou orbite au point $x \in X$}, le sous-ensemble 
$$ \mathcal O_x = \lbrace  \phi(t,x) \text{ | } (t,x)\in U \rbrace .$$

On observera que la partition de $M$ selon les courbes intégrales du champ de vecteurs $v_M$ définit un feuilletage analytique en courbes sur $M$ dont les singularités sont les zéros du champ de vecteur $v_M$.}

{\defn Soit $(M,v_M)$ une $D$-variété analytique réelle et $A \subset M$ un sous-ensemble. On note $(U,\phi)$ le flot associé à $(M,v_M)$. On dit que $A \subset M$ est \textit{$\phi$-invariant} si 
$$ \phi((\R \times A) \cap U) \subset A.$$
Autrement dit, le sous-ensemble $A \subset M$ est invariant si et seulement si pour tout $a \in A$ et tout $t \in \R$ tel que $(t,a) \in U$, on a $\phi(t,a) \in A$.}

On montre maintenant que la propriété d'invariance est de nature locale.

{\nota Soient $M$ une variété analytique réelle et $a \in M$.
On note $(M,a)$ le germe de l'espace analytique $M$ en $a$. Si $A \subset M$ est un sous-ensemble de $M$, on note $[A]_a$ le germe de $A$ en $a$, c'est-à-dire la classe d'équivalence de $A$ modulo la relation d'équivalence définie par : 
$$ A \sim_a B \text{ si et seulement s'il existe un voisinage $V \subset M$ de } a \text{ tel que } A \cap V = B\cap V.$$
La relation d'inclusion entre les sous-ensembles de $M$ induit une relation d'inclusion sur les germes donnée pour $A,B \subset M$ par :
$$[A]_a \subset [B]_a \text{ si et seulement s'il existe un voisinage $V$ de $a$ dans $M$ tels que } A \cap V \subset B.$$} 
{\lem\label{invariancelocale} Soient $(M,v_M)$ une $D$-variété analytique réelle et $A \subset M$ un sous-ensemble fermé. On note $(U, \phi)$ le flot réel associé. On a équivalence entre : 
\begin{itemize}
\item[(i)] Le sous-ensemble $A \subset M$ est $\phi$-invariant.
\item[(ii)] Pour tout $a \in A$, le germe du flot $\phi$ en $(0,a)$ noté $\phi_a : (\R \times M , (0,a)) \longrightarrow (M,a)$ vérifie
$$\phi_a([\R \times A]_{(0,a)}) \subset [A]_a.$$
\end{itemize}}

\begin{proof}
L'implication $(i) \Longrightarrow (ii)$ est tautologique. Montrons que $(ii) \Longrightarrow (i)$. Considérons $a \in A$ et notons $U_a = \lbrace t \in \R \text{ | } (t,a) \in U \rbrace$ qui est un ouvert connexe de $\R$.

L'ensemble $G = \lbrace t \in U_a \text{ | } \phi(t,a) \in A \rbrace \subset U_a$ est fermé non vide car $A \subset M$ est un sous-ensemble fermé. Pour tout $t \in G$, la propriété (ii) appliquée en $a' = \phi_t(a)$ montre que $G$ est un voisinage de $t$. On en déduit que $G$ est ouvert.

Par connexité de $U_a$, on en déduit que $G = U_a$ et donc que $A$ est $\phi$-invariant.
\end{proof}

{\defn Soient $X$ une variété algébrique sur $(\mathbb{R},0)$ et $v_X$ un champ de vecteurs rationnel sur $X$. On appelle \textit{flot réel régulier de $(X,v_X)$}, le flot réel associé à l'analytifié de la $D$-variété lisse $(U,v_{|U})$ où $U = X \setminus ( \mathrm{Sing}(X) \cup \mathrm{Sing}(v_X))$\footnote{Ici, $\mathrm{Sing}(v_X)$ désigne le complémentaire du plus grand ouvert de définition de $v_X$.}.}

{\rem  Soit $(X,v_X)$ une $D$-variété réelle. La connaissance du flot réel régulier de $(X,v_X)$ n'est pas suffisante pour étudier les sous-variétés invariantes de $(X,v_X)$ et de ses produits : Les sous-variétés invariantes $W \subset (X,v_X)$ sans points réels n'auront aucune trace dans $X(\R)$. Il suffit par exemple que $X(\R) = \emptyset$ et on ne peut rien dire du tout.}

\subsubsection{Flot complexe d'une $D$-variété complexe} En travaillant avec un flot complexe, plutôt qu'avec le flot réel, on résout complètement la difficulté précédente. Néanmoins, contrairement au cas réel, le flot complexe n'admet pas d'ouvert maximal de définition et il faut travailler localement avec des germes de flots.

{\defn On appelle \textit{$D$-espace analytique complexe}, tout couple $(M,v_M)$ où $M$ est un espace analytique complexe et $v_M$ un champ de vecteurs analytique sur $M$.}

{\cons\label{analytificationcomplexe} Soit $(X,v_X)$ une $D$-variété au dessus de $(\mathbb{C},0)$.

On munit l'ensemble $M = X(\mathbb{C})$ de sa structure analytique complexe.
Le champ de vecteurs $v_X$ induit un champ de vecteurs analytique $v_M$ sur $M$ . On obtient ainsi un $D$-espace analytique complexe $(M,v_M)$.

Comme dans le cas réel, la construction précédente détermine un foncteur de la catégorie des $D$-variétés au dessus de $(\mathbb{C},0)$ vers la catégorie des $D$-espaces analytiques complexes appelé \textit{foncteur d'analytification complexe}. Contrairement au cas réel, ce foncteur est fidèle et réalise donc la catégorie des $D$-variétés au dessus de $(\mathbb{C},0)$ comme une sous-catégorie de la catégorie des $D$-espaces analytiques complexes.

{\defn Soit $(M,v_M)$ un $D$-espace analytique complexe. Le champ de vecteurs $v_M$ induit une dérivation $\delta_M : \mathcal O_M \longrightarrow \mathcal O_M$ sur le faisceau $\mathcal O_M$ des fonctions analytiques sur $M$.
Un sous-espace analytique fermé $Z \subset (M,v_M)$ est appelée \textit{$\delta_M$-invariant} si le faisceau d'idéaux $I_Z \subset \mathcal O_M$  définissant $Z$ est stable par la dérivation $\delta_M$.}

{\rem\label{analytique-algebrique} Dans le cas où le $D$-espace analytique complexe $(M,v_M)$ est l'analytifié d'une $D$-variété au dessus de $(\mathbb{C},0)$ et où la sous-variété analytique $Z \subset (M,v_M)$ est l'analytifié d'une sous-variété algébrique, les deux notions d'invariance (analytique et algébrique) coïncident.}

{\cons\label{flotcomplexe} Soit $(M,v_M)$ une $D$-variété analytique complexe. Le théorème d'existence et d'unicité locale de solutions pour les équations différentielles analytiques complexes \citep[Théorème 1.1]{Ily} assure que pour tout point $y \in M$, il existe un voisinage connexe $U_y \subset \mathbb{C} \times M$ de $(0,y)$ et d'une application analytique $\phi : U_y \longrightarrow M$ tels que 
 $$\begin{cases}
\forall (0,x) \in U_y \text{ , } \phi(0,x) = x   \\
\forall (t,x) \in U_y \text{ , } \frac d {dt} \phi(t,x) = v_M(\phi(t,x)). \end{cases}$$
De plus, le germe de $(U_y,\phi)$ en $(y,0)$ est unique.}

{\defn Soient $(M,v_M)$ une $D$-variété analytique complexe et $y \in M$. On appelle \textit{germe de flot de $v_M$  en $y \in M$}, le germe $\phi_y : (\mathbb{C} \times M, (y,0)) \longrightarrow (M,y)$ associé par la construction \ref{flotcomplexe}.}

{\defn Soient $(M,v_M)$ une $D$-variété analytique complexe et $A \subset M$ un sous-ensemble. On dit que $A \subset M$ est \textit{$\phi$-invariant} si pour tout $a \in A$, le germe du flot $\phi$ en $(0,a)$ noté $\phi_a : (\C \times M , (0,a)) \longrightarrow (M,a)$ vérifie
$$\phi_a([\C \times A]_{(0,a)}) \subset [A]_a.$$}

{\Prop\label{invarianceequivalence} Soient $(M,v_M)$ une $D$-variété analytique complexe et $Z \subset M$ un sous-espace analytique complexe. On a équivalence entre :
\begin{itemize}
\item[(i)] La sous-variété analytique $Z \subset (M,v_M)$ est un sous-espace analytique $\delta_M$-invariant.
\item[(ii)] Le sous-ensemble $Z \subset M$ est $\phi$-invariant.
\end{itemize}}

\smallskip
Soit $a \in Z$. On note $\mathcal O_{M,a}$ l'anneau local des germes de fonctions analytiques sur $M$ en $a$. Le champ de vecteurs $v_M$ munit l'anneau 
$\mathcal O_{M,a}$ d'une dérivation notée $\delta_a : \mathcal O_{M,a} \longrightarrow \mathcal O_{M,a}$.
En notant $\mathcal O_{M,a}\lbrace t \rbrace$ l'anneau local des germes de fonctions analytiques sur $\mathbb{C} \times M$ en $(0,a)$, le germe du flot $\phi$ en $a$ induit un morphisme d'anneaux : 
$$ \phi_a^\sharp : \begin{cases} \mathcal O_{M,a} \longrightarrow \mathcal O_{M,a}\lbrace t \rbrace \\
f \mapsto f \circ \phi_a
\end{cases}$$

Nous noterons encore $\delta_a : \mathcal O_{M,a} \lbrace t \rbrace \longrightarrow  \mathcal O_{M,a} \lbrace t \rbrace $ la dérivation déduite de $\delta_a$ définie par : 
$$ \delta_a(\sum_{n \in \mathbb{N}} f_n.t^n) = \sum_{n \in \mathbb{N}} \delta_a(f_n).t^n.$$

Remarquons que le morphisme d'anneaux $._{|t= 0} : \mathcal O_{M,a} \lbrace t \rbrace \longrightarrow  \mathcal O_{M,a}$ est alors un morphisme d'anneaux différentiels.

{\lem\label{formuledeCauchy} Avec les notations précédentes, on a pour tout $f \in \mathcal O_{M,a}$,

\begin{equation}
\begin{cases}
\phi_a^\sharp(f)_{|t = 0} = f \\
\frac d {dt} \phi_a^\sharp(f) = \delta_a(\phi_a^\sharp(f))
\end{cases} \tag{$\ast$}
\end{equation}

où $\frac d {dt}$ désigne la dérivation usuelle sur $\mathcal O_{M,a}\lbrace t \rbrace$. De plus, on a : 

\begin{equation}
 \phi_a^\sharp(f) = \sum_{n = 0}^\infty \frac{\delta_a^n(f)}{n!} t^n. \tag{$\ast\ast$}
\end{equation}
 
 }
 
Les identités $(\ast)$ découlent aussitôt de la définition du flot $\phi$ (Construction \ref{flotcomplexe}). La formule $(\ast \ast)$ est une conséquence immédiate de $(\ast)$. Des formules du type $(\ast \ast)$ apparaissent déjà chez Cauchy.

\begin{proof}[Preuve de la proposition \ref{invarianceequivalence}]
Montrons que $(ii) \Longrightarrow (i)$. Il suffit de montrer que pour tout $a \in Z$, le faisceau d'idéaux $I_{Z,a} \subset \mathcal O_{M,a}$ est stable par la dérivation $\delta_a$.

Fixons $a \in Z$. Par définition, on a $\phi_a([\C \times Z]_{(0,a)}) \subset [Z]_a$. On en déduit que pour tout $f \in I_{Z,a}$, on a  :
$$\phi_a^\sharp(f) \in I_{\C \times Z , (0,a)} = I_{Z,a}. \mathcal O_{M,a}\lbrace t \rbrace$$
Le lemme \ref{formuledeCauchy} permet de conclure que $\delta_a(f)  = \frac d {dt} \phi_a^\sharp(f)_{|t = 0}  \in I_{Z,a}$, ce qui montre que $Z \subset (M,v_M)$ est bien une sous-variété analytique invariante.

Réciproquement, supposons que $Z \subset (M,v_M)$ est une sous-variété invariante. Fixons $a \in Z$. On a donc pour tout $f \in I_{Z,a} \subset \mathcal O_{M,a}$ et tout entier $n \geq 1$, $\frac 1 {n!} \delta^n_a(f) \in I_{Z,a}$. Le lemme \ref{formuledeCauchy} permet d'affirmer que
$$f \in I_{Z,a} \Longrightarrow \phi_a^\sharp(f) \in 
 I_{Z,a}. \mathcal O_{M,a}\lbrace t \rbrace =  I_{\C \times Z , (0,a)} .$$
 
On en déduit que $\phi_a([\C \times Z]_{(0,a)}) \subset [Z]_a$, ce qu'il fallait démontrer.
\end{proof}

\subsubsection{Compatibilité des différentes notions d'invariances}

{\Prop\label{invariancefinal1} Soient $(X,v_X)$ une $D$-variété lisse au dessus de $(\mathbb{R},0)$ et $Z \subset X$ une sous-variété fermée invariante (définie sur $\R$). On note $(M,\phi)$ le flot réel de $(X,v_X)$. L'ensemble $Z(\R) \subset M$ est un fermé $\phi$-invariant.}  

Autrement dit, si $Z \subset X$ est une sous-variété fermée invariante alors pour tout $x \in Z(\R)$, l'orbite $\mathcal O_x$ de $x$ est contenue dans $Z$.
\begin{proof}
Soit $Z \subset X$ une sous-variété fermée invariante. D'après la remarque \ref{analytique-algebrique}, la sous-variété fermée analytique $Z(\mathbb{C})^{an} \subset M_\mathbb{C} = X(\C)^{an}$ est $\delta_{M_\C}$-invariante. D'après la proposition \ref{invarianceequivalence}, on en déduit que, en notant $\phi_a$ le germe du flot complexe en $a \in M_\mathbb{C}$, on a : 
$$\phi_a([\C \times Z(\mathbb{C})]_{(0,a)}) \subset [Z(\mathbb{C})]_a$$

Si $a \in M \subset M_\mathbb{C}$, le germe du flot réel en $a$ est la restriction du germe du flot complexe aux nombres réels et vérifie : 
$$\phi_a( [\R \times M]_{(0,a)}) \subset [M]_a.$$

En prenant les intersections, on obtient que pour tout $a \in X(\R)$ : 
$$\phi_a([\R \times Z(\mathbb{R})]_{(0,a)}) \subset [Z(\mathbb{R})]_a.$$
Le lemme \ref{invariancelocale} permet de conclure que $Z(\mathbb{R})$ est $\phi$-invariant. 
\end{proof}

Voici une réciproque partielle au résultat précédent. Elle ne sera pas utilisée dans la preuve du critère d'orthogonalité aux constantes.

{\Prop\label{invariancefinal2} Soit $(X,v_X)$ une $D$-variété lisse au-dessus de $(\mathbb{R},0)$. On note $(M,\phi)$ le flot réel de $(X,v_X)$. Si $A \subset X(\mathbb{R})$ est un sous-ensemble $\phi$-invariant alors la clôture de Zariski $\overline{A}$ de $A$ est une sous-variété fermée invariante de $(X,v_X)$.}

\begin{proof}
Remarquons d'abord que d'après le lemme \ref{generique}, on peut supposer que $X$ est affine.
On note $Z = \overline{A}$ la clôture de Zariski  de $A$ dans $X$.

Il suffit de vérifier que l'idéal  $I \subset \mathcal O_M(M)$ des fonctions analytiques qui s'annulent sur $A$ est stable par la dérivation $\delta_v : \mathcal O_M \longrightarrow \mathcal O_M$ induite par le champ de vecteurs $v_M$ sur $M$.
En effet, en utilisant l'inclusion $\R[X] \subset \mathcal O_M(M)$, l'idéal $I_Z \subset \mathbb{R}[X]$ définissant $Z$ est donné par
$$I_Z = I \cap \R[X] \subset \mathcal O_M(M).$$
Ces deux sous-ensembles sont stables par la dérivation $\delta_v$, et donc $Z$ est une sous-variété invariante de $(X,v_X)$.

Montrons donc que le faisceau d'idéaux $\mathcal I \subset \mathcal O_M$ des fonctions analytiques s'annulant sur $A$ est invariant par la dérivation $\delta_v$.

Soient $f \in \mathcal I$ et $x \in M$. En notant $\phi_x : (\R \times M , (0,x)) \longrightarrow (M,x)$ le germe du flot $\phi$ en $x$, le germe de $\delta_v(f)$ en $x$ est donné par : 
$$\delta_v(f)_x = \frac d {dt} [f \circ \phi_x]_{|t = 0}.$$

Comme $A$ est un sous-ensemble $\phi$-invariant, on a $f \circ \phi_x([\R \times A]_{(0,x)}) = 0$ et donc 
$$\delta_v(f)_x([A]_x) = \frac d {dt} [f \circ \phi_x]_{|t = 0}([A]_x) = 0.$$ On en déduit que le germe $\delta_v(f)_x$ de $\delta_v(f)$ en $x$ appartient à $\mathcal I_x$ pour tout $x \in M$ et donc que $\delta_v(f) \in I$. 
\end{proof}

{\cor Soit $(X,v_X)$ une $D$-variété au dessus de $(\mathbb{R},0)$. La clôture de Zariski de $X(\R)$ dans $X$ est une sous-variété fermée invariante.}

Il suffit d'appliquer successivement les deux proposition \ref{invariancefinal1} et \ref{invariancefinal2}.  Cet énoncé a la forme d'un énoncé d'algèbre différentielle mais il est propre à l'extension de corps $\mathbb{R} \subset \mathbb{C}$. Il n'est pas difficile de construire des exemples de $D$-variétés $(X,v_X)$ définie sur le corps $\mathbb{Q}$ des nombres rationnels telles que la clôture de Zariski de $X(\mathbb{Q})$ ne soit pas invariante.

Par exemple, si $C$ est une courbe affine sur le corps $\mathbb{Q}$  admettant un nombre fini de points rationnels et si $v_C$ est un champ de vecteurs régulier sur $C$ qui ne s'annule pas en tous les points de $C(\mathbb{Q})$  alors la sous-variété fermée $ \overline{C(\mathbb{Q})} = C(\mathbb{Q}) \subset (C,v_C)$ n'est pas invariante.}

\subsection{Flots topologiquement transitifs et topologiquement mélangeants}

Dans cette partie, on se restreint à l'étude de flots réels complets $(M , (\phi_t)_{t \in \mathbb{R}})$. De plus, au lieu de considérer des flots analytiques réels, on considère plus généralement des flots continus sur un espace métrisable, ce qui est le cadre naturel pour les résultats de cette partie.

Un flot réel continu s'identifie alors à une action du groupe additif $(\R , +)$ des nombres réels sur un espace métrisable $M$ et donc à un système dynamique  pour le groupe additif des nombres réels.

On renvoie aux trois premières sections du chapitre 1 de \cite{Glas} pour des résultats et notions analogues dans le cadre plus général des systèmes dynamiques sur  un groupe topologique localement compact et de cardinal $\leq 2^{\aleph_0}$. Dans les  prochains paragraphes, on donne une présentation  autonome des résultats qui nous seront utiles dans la suite.

\subsubsection{Flot topologiquement transitif} On fixe $M$ un espace topologique métrisable.

{\defn Soit $(M,(\phi_t)_{t \in \R})$ un flot continu complet. On dit que $(M,(\phi_t)_{t \in \R})$ est \textit{un flot topologiquement transitif} si pour tout couple d'ouverts non vides $U,V \subset M$, il existe $t \in \R$ tel que $\phi_t(U) \cap V \neq \emptyset$.}

{\nota Soit $(M,(\phi_t)_{t \in \R})$ un flot continu complet. Pour tout couple d'ouverts non vides $U,V \subset X$, on note 
$$ N(U,V) = \lbrace t \in \R \text{ | } \phi_t(U) \cap V \neq \emptyset \rbrace \subset \R$$

Un flot complet  $(M,(\phi_t)_{t \in \R})$ est topologiquement transitif si et seulement si pour tout couple d'ouverts non vides $U,V \subset M$, le sous-ensemble $N(U,V) \subset \R$ est non vide.}

{\rem  Pour tout couple d'ouverts non vides $U,V \subset \R$, le sous-ensemble $N(U,V) \subset \R$ est un ouvert de $\R$. On en déduit que $N(U,V)$ s'écrit comme une union (au plus) dénombrable et éventuellement vide d'intervalles ouverts disjoints de $\R$.

De plus, remarquons que $N(U,V) = - N(V,U) \subset \R$.}

{\exam\label{exempletopologiquementransitif1} Si $(M,(\phi_t)_{t \in \R})$ est le flot continu complet sur le disque unité $\mathbb{S}^1 = \mathbb{R} / 2\pi \mathbb Z$ induit par le champ de vecteurs unitaire sur $\mathbb{S}^1$ alors pour tout couple d'intervalles ouverts non vides $U,V \subset \mathbb{S}^1$, il existe $a < b \in \R$ tels que   
$$N(U,V)=  \bigcup_{k \in \mathbb{Z}}]a + 2.\pi k \mathbb{Z} ; b + 2.\pi k \mathbb{Z}[.$$

Si $(M,(\phi_t)_{t \in \R})$ est le flot continu complet sur la droite réelle induit par le champ de vecteurs constant $\frac d {dt}$ alors pour tout couple d'intervalles ouverts bornés non vides $U,V \subset \mathbb{R}$, il existe $a < b \in \R$ tels que
$$N(U,V)= ]a,b[.$$

En particulier, les deux flots précédents sont topologiquement transitifs.}

Avec des hypothèses topologiques raisonnables, on a la caractérisation suivante des flots complets continus topologiquement transitifs.

{\Prop\label{orbitedense} Soit $(M,(\phi_t)_{t \in \R})$ un flot complet continu. Supposons que $M$ soit localement compact et séparable. Les propriétés suivantes sont alors équivalentes : 
\begin{itemize}
\item[(i)] Le flot $(M,(\phi_t)_{t \in \R})$ est topologiquement transitif.
\item[(ii)] Il existe une orbite dense dans $(M,(\phi_t)_{t \in \R})$. 
\end{itemize}}

\begin{proof}
$(ii) \Rightarrow (i)$ Soient  $U,V \subset M$ un couple d'ouverts non vide de $M$. Considérons $x \in M$ dont l'orbite est dense dans $M$.

Comme l'orbite de $x$ est dense dans $M$, il existe un réel $t \in \R$ tel que $\phi_t(x) \in U$ et un entier $s \in \R$ tel que $\phi_s(x) \in V$. On en déduit que
$$ \phi_{-s}(V) \cap \phi_{-t}(U) \neq \emptyset \text{ et donc }\phi_{s - t}(U) \cap V \neq \emptyset.$$ 

$(i) \Rightarrow (ii)$ Considérons une base dénombrable d'ouverts $(U_k)_{k \in \mathbb N}$ de l'espace métrique séparable $M$.
Remarquons qu'un sous-ensemble de $M$ est dense si et seulement s'il rencontre tous les $(U_k)_{k \in \N}$.

On construit par récurrence sur $k \in \mathbb{N}^\ast$ une suite de compacts d'intérieur non vide $$V_k \subset V_{k-1} \subset \cdots \subset V_1$$ tels que l'orbite de tout point de $V_i$ rencontre les ouverts $U_1 , \cdots , U_i$.

\begin{itemize}
\item Pour $n = 1$, comme $M$ est localement compact, il existe un compact $V_1 \subset U_1$ d'intérieur non vide.

\item Supposons $V_k \subset V_{k-1} \subset \cdots \subset V_1$ construits. Par hypothèse, $\overset{\circ} {V_k}$ est un ouvert non vide. Comme le flot est topologiquement transitif, il existe $t \in \R$ tel que $\overset{\circ} {V_k} \cap \phi_t(U_{k+1}) \neq \emptyset$.
L'ensemble $\overset{\circ} {V_k} \cap \phi_t(U_{k+1})$ est un ouvert non vide  de $M$ et par construction, l'orbite de chacun de ses points rencontre $U_1, \cdots , U_{k+1}$.
Comme $M$ est localement compact, il suffit de choisir un compact $V_{k + 1}$ d'intérieur non vide tel que $V_{k+1} \subset\overset{\circ} {V_k} \cap \phi_t(U_k)$.
\end{itemize}
Une intersection de compacts non vides emboités étant toujours non vide,  considérons $x \in \bigcap_{k \in \mathbb{N}} V_k$. Par construction, l'orbite de $x$ rencontre tous les $(U_k)_{k \in \mathbb{N}}$ et est donc dense dans $M$. 
\end{proof}

{\exam\label{tore} Considérons le tore $\mathbb{T}^n = \mathbb{R}^n/\mathbb{Z}^n$ et le flot $(\phi_t)_{t \in \R}$  défini sur $\mathbb T^n$ par le système d'équations différentielles linéaires
$$ \begin{cases} \frac {dy_1} {dt}  & = \omega_1 \\
& \vdots  \\
\frac {dy_n} {dt}  & = \omega_n \end{cases}$$
où $\omega = (\omega_1 , \cdots , \omega_n) \in \mathbb T^n$. Le flot complet $(\mathbb{T}^n,(\phi_t)_{t \in \R})$ défini par le système d'équations différentielles précédent est appelé \textit{le flot constant sur le tore $\mathbb{T}^n$ de fréquences $\omega_1,\cdots, \omega_n$}.

Le flot $(\mathbb{T}^n,(\phi_t)_{t \in \R})$ est topologiquement transitif si et seulement si les fréquences $\omega_1, \cdots, \omega_n$ sont $\Z$-indépendantes \citep[Proposition 1.5.1]{KH}.}

{\rem\label{produittopologiquementtransitif} Un produit de deux flots complets topologiquement transitifs n'est pas toujours topologiquement transitif : 

Considérons les flots constants sur $\mathbb{T}^n$ de fréquences respectives $(\omega_1 , \cdots \omega_n)$ et $(\tau_1 , \cdots \tau_n)$. Le flot produit est le flot constant sur $\mathbb{T}^{2n}$ de fréquences $(\omega_1 , \cdots \omega_n, \tau_1 , \cdots \tau_n)$.
D'après l'exemple \ref{tore}, le flot produit est topologiquement transitif si et seulement si $(\omega_1 , \cdots \omega_n, \tau_1 , \cdots \tau_n)$ sont $\mathbb{Z}$-linéairement indépendants. En particulier le produit d'un flot constant sur le tore avec lui-même n'est jamais topologiquement transitif.}

\subsubsection{Flots faiblement topologiquement mélangeants}
La condition "être topologiquement transitif" n'est donc pas stable par produit. On s'intéresse donc désormais à la classe des flots complets continus $(M,(\phi_t)_{t \in \R})$ tel que pour tout $n \in \mathbb{N}$, le flot produit $(M,(\phi_t)_{t \in \R})^n$ est topologiquement transitif.

{\Prop\label{produittopologique} Soit $(M,(\phi_t)_{t \in \R})$ un flot continu. Les propriétés suivantes sont équivalentes
\begin{itemize}
\item[(i)] Le flot $(M,(\phi_t)_{t \in \R}) \times (M,(\phi_t)_{t \in \R}) $ est topologiquement transitif.
\item[(ii)] Pour tout $n \in \mathbb{N}$, le flot $(M,(\phi_t)_{t \in \R}) \times \cdots \times (M,(\phi_t)_{t \in \R})$ est topologiquement transitif.
\item[(iii)] La famille des $N(U,V)$ où $U$ et $V$ parcourent l'ensemble des ouverts non vides de $M$ est une famille filtrante de $\mathcal P(\R)$. 
\end{itemize}}

Rappelons que si $E$ est un ensemble, une partie filtrante de $\mathcal P(E)$ est une partie $\mathcal F \subset \mathcal P(E) \setminus \lbrace \emptyset \rbrace$  vérifiant : pour tout $A,B \in \mathcal F$ il existe $C \in \mathcal F$ tel que $C \subset A \cap B$.

\begin{proof}
$(ii) \Rightarrow (i)$ est immédiat. Montrons que $(iii) \Rightarrow (ii)$.
Soit $n \in \mathbb{N}$. Pour vérifier que $(M,(\phi_t)_{t \in \R})^n$ est topologiquement transitif, il suffit de vérifier que $N(U,V) \neq \emptyset$ pour $U$ et $V$ parcourant une base d'ouverts de $M^n$.

Les ouverts de la forme $U_1 \times \cdots \times U_n$, où $U_1 , \cdots , U_n \subset M$ sont des ouverts de $M$, engendrent la topologie produit sur $M^n$. On en déduit que $(M,(\phi_t)_{t \in \R})^n$ est topologiquement transitif si et seulement si 
$$ N(U_1 \times \cdots \times U_n) \cap N(V_1 \times \cdots \times V_n) = N(U_1,V_1) \cap \cdots \cap N(U_n,V_n) \neq \emptyset $$
lorsque $U_1, \cdots , U_n, V_1 \cdots , V_n$  parcourent les ouverts non vides de $M$.

Soient $U_1, \cdots, U_n, V_1 \cdots , V_n$  des ouverts non vides de $M$. D'après la propriété (iii), il existe des ouverts non vides  $A,B \subset M$ tels que $N(A,B) \subset N(U_{n-1},V_{n-1}) \cap N(U_n,V_n)$. On en déduit que : 
$$N(U_1,V_1) \cap \cdots \cap N(U_{n-2},V_{n-2}) \cap N(A,B) \subset N(U_1 \times \cdots \times U_n) \cap N(V_1 \times \cdots \times V_n).$$

Par induction, il existe des ouverts non vides $\tilde{A},\tilde{B} \subset M$ tels que : 

$$N(\tilde{A},\tilde{B}) \subset N(U_1 \times \cdots \times U_n) \cap N(V_1 \times \cdots \times V_n).$$

En particulier, on obtient que  $N(U_1 \times \cdots \times U_n) \cap N(V_1 \times \cdots \times V_n) \neq \emptyset$ car $N(\tilde{A},\tilde{B}) \neq \emptyset$. 

Montrons que $(i) \Rightarrow (iii)$. Soient $U_1,V_1,U_2,V_2$ des ouverts non vides de $M$. Le système dynamique $(M,(\phi_t)_{t \in \R}) \times (M,(\phi_t)_{t \in \R})$ étant topologiquement transitif, on a
$$ N(U_1,V_1) \cap N(U_2,V_2) \neq \emptyset.$$
En particulier, $N(U_1,V_1)$ est toujours non vide.

De plus, considérons $t \in N(U_1,V_1) \cap N(U_2,V_2)$ et posons $A = \phi_t(U_1) \cap U_2$ et $B = \phi_t(V_1) \cap V_2$. Par hypothèse sur $t$, $A$ et $B$ sont des ouverts non vides de $M$ et on vérifie facilement que   
$$N(A,B) \subset N(U_1,V_1) \cap N(U_2,V_2). \qedhere $$
\end{proof}

{\defn\label{definitionfaiblemelange} Soit $(M,(\phi_t)_{t \in \R})$ un flot complet continu. On dit que $(M,(\phi_t)_{t \in \R})$ est \textit{un flot faiblement topologiquement mélangeant} si l'une des trois conditions équivalentes de la proposition \ref{produittopologique} est vérifiée.}

{\exam Le calcul direct des $N(U,V)$ pour les exemples \ref{exempletopologiquementransitif1}  montrent que ces deux flots ne sont pas faiblement topologiquement mélangeants.

De même, la remarque \ref{produittopologiquementtransitif} montre que les flots constants sur le tore  (exemple \ref{tore}) ne sont pas faiblement topologiquement mélangeants non plus.

L'exemple fondamental de flot faiblement topologiquement mélangeant est le flot géodésique unitaire sur une variété réelle compacte munie d'une métrique à courbure strictement négative (voir \cite{Ano}, \cite{Coud}, \cite{Dalbo}).
Cet exemple sera le point de départ des applications pour le critère d'orthogonalité aux constantes présenté dans ce texte et sera discuté dans une suite de cet article. }

\section{Un critère dynamique d'orthogonalité aux constantes} Dans cette section, on établit le critère annoncé d'orthogonalité aux constantes pour les $ D$-variétés réelles ainsi que sa version en familles annoncées dans l'introduction.

Le critère d'orthogonalité aux constantes que nous présentons  (Théorème \ref{montheoreme}) consiste en \textit{une obstruction topologique à la présence d'intégrale première rationnelle} pour ce flot complet et tous ses produits. L'argument essentiel est donné par une notion élémentaire de la théorie des systèmes dynamiques topologiques : les flots faiblement topologiquement mélangeants. 

En s'appuyant sur les résultats de spécialisation de la deuxième partie de ce texte, nous démontrons ensuite une variante de ce critère d'orthogonalité aux constantes (Théorème \ref{montheoremeenfamille}) pour \textit{les membres ``très génériques'' de familles  lisses de $D$-variétés absolument irréductibles a partir de l'étude dynamique d'une des fibres réelles}.

\subsection{Indétermination et intégrale première rationnelle} On étudie plus généralement les conséquences de l'existence d'une intégrale première rationnelle pour les $D$-variétés définies sur des corps de constantes. On fixe $k$ un corps de caractéristique $0$.

{\nota Soient $X$ une variété algébrique irréductible sur le corps $k$ et $f \in k(X)$ une fonction rationnelle. Alors $f$ définit un morphisme rationnel de variétés $f : X \dashrightarrow \mathbb{P}^1$ qui admet un plus grand ouvert de définition $U \subset X$.

Le graphe $G_f \subset U \times \mathbb{P}^1$ de $f$ s'identifie à une sous-variété fermée irréductible de $U \times \mathbb{P}^1$. Sa clôture de Zariski dans $X \times \mathbb{P}^1$ sera notée $\overline{G_f}$.}

{\lem\label{graphe} Soient $(X,v_X)$ une $D$-variété irréductible au dessus de $(k,0)$ et $f \in k(X)$ une fonction rationnelle. On a équivalence entre : 
\begin{itemize}
\item[(i)]La fonction rationnelle $f$ est une intégrale première rationnelle de $(X,v_X)$.
\item[(ii)]La sous-variété fermée $\overline{G_f} \subset (X,v_X) \times (\mathbb{P}^1,0)$ est une sous-variété fermée invariante.
\end{itemize}}

\begin{proof}
Considérons $U \subset X$ un ouvert affine tel que $f_{|U} : U \longrightarrow \mathbb{A}^1$ est une fonction régulière. D'après les lemmes \ref{generique} et \ref{Zariskiouvert}, il suffit de vérifier l'équivalence pour $f_{|U}$. On peut donc supposer que $U = X$, c'est à dire que $X$ est affine et que $f \in k[X]$.

Le graphe $G_f \subset X \times \mathbb{A}^1$ de $f$ est alors l'hypersurface de $X \times \mathbb{A}^1$ d'équation  $Y - f = 0.$

On en déduit que $G_f$ est une sous-variété invariante de $(X,v_X) \times (\mathbb{A}^1,0)$ si et seulement s' il existe $h(Y) \in k[X][Y]$ tel que 
$$ \delta_X(f) = h(Y) (Y - f)\in k[X][Y].$$

En comparant les dégrés en $Y$, on voit que nécessairement $h(X,Y) = 0$ et donc que la seconde condition est équivalente à $\delta_X(f) = 0$, autrement dit que $f$ est une intégrale première rationnelle de $(X,v_X)$.
\end{proof}

{\defn Soient $X$ une variété algébrique irréductible et $f \in k(X)$ une fonction rationnelle. On appelle \textit{lieu d'indétermination de $f$} et on note $\mathrm{Ind}(f)$, le complémentaire du plus grand ouvert de définition de $f : X \dashrightarrow \mathbb{P}^1$.}
Le lieu d'indétermination de $f$ est donc une sous-variété fermée de $X$.

{\rem\label{Maintheorem} Soient $X$ une variété algébrique irréductible et lisse et $f : X \dashrightarrow \mathbb{P}^1$ une fonction rationnelle sur $X$. En notant $\pi_f : \overline{G_f} \longrightarrow X$ la restriction de la première projection à $\overline{G_f}$, d'après le Main Theorem de Zariski \citep[Théorème 3.20]{Mum1}, on a la description suivante de $\mathrm{Ind}(f)$ :
$$ \mathrm{Ind}(f) = \lbrace x \in X \text{ | } \pi_f^{-1}(x) \text{ admet une composante irréductible de dimension }> 0 \rbrace$$

Comme $\overline{G_f} \subset X \times \mathbb{P}^1$, la condition précédente se réécrit : 

$$ \mathrm{Ind}(f) = \lbrace x \in X \text{ | } \{x \} \times \mathbb{P}^1 \subset \overline{G_f}\rbrace.$$}

{\lem\label{indetermination} Soient $(X,v_X)$ une $D$-variété irréductible et lisse au-dessus de $(k,0)$ et $f \in k(X)$ une intégrale première rationnelle. Le lieu d'indétermination de $f$, $\mathrm{Ind}(f) \subset (X,v_X)$ est une sous-variété fermée invariante.}

\begin{proof}
On note $Z = \mathrm{Ind}(f) \subset X$ le lieu d'indétermination de $f$ qui est donc une sous-variété fermée  de $X$. On commence par décrire le faisceau d'idéaux associé à $Z$ avant dé vérifier qu'il est invariant.

Considérons $U \subset X$ un ouvert affine. La sous-variété fermée irréductible $G_U := \overline{G_f} \cap (U \times \mathbb{A}^1) \subset U \times \mathbb{A}^1 $ est une hypersurface. Comme $U \times \mathbb{A}^1$ est lisse, il existe $F \in k[U][T]$ tel que $G_U \subset U \times \mathbb{A}^1$  est l'hypersurface d'équation $ F(X,T) = 0 .$
Fixons $t_0 \in k$. D'après la remarque \ref{Maintheorem}, on a alors 
$$ Z = \lbrace x \in X \text{ | } F(x,.) = 0 \rbrace = \lbrace x \in X \text{ | } \forall n \in \mathbb{N} \text{ , } \frac \partial {\partial t^n} F(x,t_0) = 0\rbrace.$$

On en déduit que $I_Z(U) \subset \mathcal O_X(U)$ est l'idéal engendré par $\frac \partial {\partial t^n} F(X,t_0)$ où $n$ parcourt $\mathbb{N}$.

D'après le lemme \ref{graphe}, la sous-variété fermée $G_U \subset (U,v_X) \times (\mathbb{A}^1,0)$ est une sous-variété invariante, donc il existe $h(X,Y) \in k[U][T]$ tel que
$$ \delta_{ X\times \mathbb{A}^1}(F(X,T)) = h(X,T) F(X,T).$$

Les dérivations $\delta_{X \times \mathbb{A}^1}$ et $\frac \partial {\partial t}$ commutent (car $\delta_{X \times \mathbb{A}^1}(T) = 0$) et donc pour tout $n \in \mathbb{N}$,
$$ [\delta_{X \times \mathbb{A}^1} \frac \partial {\partial t^n} F](X,t_0) = [\frac \partial {\partial t^n} \delta_{X \times \mathbb{A}^1}(F)](X,t_0) = [\frac \partial {\partial t^n} (h.F)) ](X,t_0)$$  
Comme $I_Z(U)$ est l'idéal engendré par $\frac \partial {\partial t^n} F(X,t_0)$ où $n$ parcourt $\mathbb{N}$, la règle de Leibniz pour les dérivations supérieures d'un produit permet alors de conclure que $I_Z(U) \subset \mathcal O_X(U)$ est un idéal invariant et donc que $\mathrm{Ind}(f)$ est une sous-variété fermée invariante.
\end{proof}

{\Prop\label{consequenceintegralepremiere} Soit $(X,v_X)$ une $D$-variété irréductible et lisse au dessus de $(k,0)$ admettant une intégrale première rationnelle non constante. En notant $\mathrm{Inv}$ l'ensemble des sous-variétés fermées invariantes strictes de $(X,v_X)$, on a :
$$ X(k) = \bigcup_{Z \in \mathrm{Inv}} Z(k). $$}

\begin{proof}
Considérons $f : X \dashrightarrow \mathbb{P}^1$ une intégrale première rationnelle non constante et notons $U$ son plus grand ouvert de définition.

Considérons $x \in X(k)$. On a deux cas : 
\begin{itemize}
\item Si $x \in U(k)$ alors en notant $a = f(x)$, $f^{-1}(a)$ est une sous-variété fermée invariante de $(U,v_X)$ contenant $x$ et de codimension 1. Sa clôture de Zariski dans $X$ est une sous-variété fermée invariante (d'après le lemme \ref{generique}) et stricte de $(X,v_X)$.
\item Si $x \notin U(k)$ alors $x \in \mathrm{Ind}(f)$ qui est une sous-variété fermée invariante stricte d'après le lemme \ref{indetermination}. \qedhere
\end{itemize}
\end{proof}

\subsection{Critère d'orthogonalité aux constantes pour les $D$-varietes reelles} 

{\lem Soient $(M,v_M)$ une $D$-variété analytique réelle et $K\subset M$ un compact $\phi$-invariant. On note $(U,\phi)$ le flot réel associé. La restriction du flot $\phi$ au sous-ensemble $K$ est un flot complet.}

\begin{proof}
Pour tout $x \in M$, notons $U_x = ]t^-(x) ; t^+(x)[$ où $t^-(x) < 0 < t^+(x)$, le plus grand ouvert de définition de la courbe intégrale $\phi_x$ de $v_M$ en $x \in M$ vérifiant $\phi_x(0) = x$.

Les fonctions $t^+ : M \longrightarrow \R \cup \lbrace -\infty;\infty \rbrace$ et $t^- : M \longrightarrow \R \cup \lbrace -\infty;\infty \rbrace$ sont des fonctions continues. Comme $K$ est compact, le réel $\epsilon = \mathrm{inf}_{x \in K} t^+(x)$ est strictement positif.

On raisonne par l'absurde en considérant $x_0 \in K$ tel que, par exemple, $t_0 = t^+(x) < \infty$. Posons $x_1 = \phi_{t_0- \frac \epsilon 2}(x) \in K$ car $K$ est $\phi$-invariant.

La courbe définie par : 
$$ \overline{\phi_x} : \begin{cases}
]t^-(x) ; t^+(x) + \frac \epsilon 2 [ \longrightarrow M \\
t \mapsto \begin{cases} \phi_{x_0}(t) \text{ si } t \leq t_0 \\
                         \phi_{x_1}(t - t_0 - \frac \epsilon 2) \text{ si } t \geq  t_0 - \frac \epsilon 2
\end{cases}                         
\end{cases}$$
est bien définie par la propriété de semi-groupe du flot associé au champ de vecteurs $v_M$. On vérifie aisément que $\overline{\phi_x}$ est une courbe intégrale de $v_M$ vérifiant $\overline{\phi}_x(0) = x$, ce qui contredit la maximalité de $U_x$.

De même, on obtient que $t^-(x) = - \infty$ pour tout $x \in K$ et donc que la restriction du flot $\phi$ au compact $K$ est complet. 
\end{proof}

{\Thm\label{montheoreme} Soient $X$ une variété absolument irréductible sur $\R$ et $v$ un champ de vecteurs rationnel sur $X$. On note $(M, \phi)$ le flot régulier réel de $(X,v_X)$.
Supposons qu'il existe un compact $K  \subset M$ Zariski-dense dans $X$ et invariant par le flot $\phi$. La restriction du flot $\phi$ à $K$ est alors un flot métrisable complet.

Si $(K,(\phi_{t |K})_{t \in \R})$ est faiblement topologiquement mélangeant alors le type générique de $(X,v)$ est orthogonal aux constantes.}

\begin{proof}
La conclusion reste inchangée lorsque l'on remplace $X$ par $U = X \setminus ( \mathrm{Sing}(X) \cup \mathrm{Sing}(v))$. On peut donc supposer que $(X,v)$ est une $D$-variété absolument irréductible et lisse au dessus de $(\mathbb{R},0)$.

{\lem Soit $(X,v)$ une $D$-variété irréductible lisse au dessus de $(\mathbb{R},0)$. Si le flot réel associé $(M,\phi)$ admet une orbite Zariski-dense dans $X$ alors $(X,v)$ est sans intégrale première rationnelle non constante.}

\begin{proof}
Considérons $x \in M = X(\R)$ dont l'orbite est Zariski-dense dans $X$ et $Z \subset X$ une sous-variété fermée invariante contenant $x$.
D'après la proposition \ref{invariancefinal1}, l'orbite $\mathcal O_x$ de $x \in M$ est contenue dans $Z(\R)$ et comme cette orbite est Zariski-dense dans $X$, on en déduit que $Z = X$.
On en déduit que $x \notin \bigcup_{Z \in \mathrm{Inv}} Z(\R)$. La proposition \ref{consequenceintegralepremiere} permet de conclure que $(X,v)$ est sans  intégrale première rationnelle.
\end{proof}

La fin de la démonstration est alors formelle : Puisque $(K,(\phi_{t |K})_{t \in \R})$ est faiblement topologiquement mélangeant, pour tout $n \in \mathbb{N}$, le flot $(K,(\phi_{t |K})_{t \in \R})^n$ est topologiquement transitif et admet donc une orbite $\mathcal O_n \subset K^n$ dense pour la topologie analytique (Proposition \ref{orbitedense}).

De plus, $K \subset X$ est Zariski-dense et donc $K^n \subset X^n$ aussi pour tout $n \in \mathbb{N}$. La topologie analytique étant plus fine que la topologie de Zariski, on en déduit que pour tout $n \in \mathbb{N}$, $(X,v)^n$ admet $\mathcal O_n$ pour orbite Zariski-dense.

D'après le lemme précédent, la $D$-variété $(X,v)^n$ est sans intégrale première rationnelle non constante pour tout $n \in \mathbb{N}$. Son type générique est donc orthogonal aux constantes d'après le théorème \ref{orthogonaliteconstantes}.
\end{proof}

\subsection{Critère d'orthogonalité aux constantes pour les familles de $D$-variétés} 

{\Thm\label{montheoremeenfamille} Soient $k$ un sous-corps des nombres réels, $S$ une variété algébrique lisse et irréductible au dessus de $k$ et  $f : (\mathcal X,v) \longrightarrow (S,0)$ une famille lisse de $D$-variétés absolument irréductibles à paramètres dans $S$. 

Supposons qu'il existe un point $p \in S(\R)$ et un compact $K \subset \mathcal X_p(\R)^{an}$ Zariski-dense dans $X$ et invariant par le flot $\phi$ du champ de vecteurs $v_{|X_p}$.

Si $(K,(\phi_{t |K})_{t \in \R})$ est faiblement topologiquement mélangeant alors il existe un ensemble dénombrable $\lbrace Z_i : i \in \mathbb{N} \rbrace$ de sous-variétés fermées algébriques strictes (sur le corps des nombres réels) $Z_i$ de $S$ tel que:

$$\forall s \in S(\mathbb C) \setminus \bigcup_{i \in \mathbb{N}} Z_i(\mathbb C)\text{, } (\mathcal X,v)_s \text{ est orthogonal aux constantes}.$$}

\begin{proof}
Notons d'abord qu'on peut toujours supposer que $k$ est finiment engendré et donc un corps dénombrable.

Sous les hypothèses du théorème \ref{montheoremeenfamille}, le théorème \ref{montheoreme} implique que le type générique de la fibre $(\mathcal X,v)_p$ est orthogonal aux constantes.

Notons $\eta$ le point générique de $S$. Le théorème \ref{specialisationtheoreme2}  (ou plus précisément sa contraposée) implique alors que le type générique de $(\mathcal X,v)_\eta$ est orthogonal aux constantes.

En utilisant le lemme \ref{typedefinissable}, pour toutes les realisations

Considérons  maintenant l'ensemble $\mathcal S$ des sous-variétés fermées strictes de $S$ (au dessus de $k$). Puisque $k$ est dénombrable, l'ensemble $\mathcal S$ est dénombrable.

Soit $s \in S(\mathbb C) \setminus \bigcup_{Z  \in \mathcal S} Z(\mathbb C)$. Par construction, $s$ réalise le type générique de S (au sens de $\textbf{ACF}_0$) et donc le type générique de $(S,0)$ (au sens de  $\textbf{DCF}_0$) puisque le corps des constantes est un pur corps algébriquement clos stablement plongé. Le corollaire \ref{typedefinissable} implique donc que $(X,v)_s$ est orthogonal aux constantes.
\end{proof}

\appendix

\section{$D$-variétés et ensembles définissables associés}

Dans cette partie, nous étudions les interactions entre les ensembles définissables dans un corps différentiellement clos et la notion de $D$-variété, introduite par A. Buium dans \cite{Bui}. L'utilisation des $D$-variétés pour l'étude des corps différentiellement clos est fréquente en théorie des modèles (\cite{Itai}, \cite{Pil}). La relation entre la catégorie des $D$-schémas et les ensembles définissables (et les types) dans un corps différentiellement clos est analogue à la relation entre les schémas et les ensembles définissables dans un corps algébriquement clos.

Dans la première section, nous rappelons la définition et les propriétés structurelles des $D$-schémas ainsi que de la notion associée de sous-schéma invariant. Dans la deuxième section, nous nous concentrons sur la relation entre $D$-variété et ensemble définissable dans la théorie $\textbf{DCF}_0$. Enfin, dans la dernière section, nous prouvons, à l'aide du principe de reflexivité de Shelah dans une théorie stable, un premier critère d'orthogonalité aux constantes (Théorème \ref{orthogonaliteconstantes}) qui sera développé et enrichi pour les $D$-variétés définie sur le corps des nombres réels dans la troisième partie de ce texte.

\subsection{$D$-schémas} 

\subsubsection{Catégorie des $D$-schémas} 

{\defn On appelle $\textit{D-schéma}$ tout couple $(X,\delta_X)$ où $X$ est un schéma et $\delta_X : \mathcal O_X \longrightarrow \mathcal O_X$ est une dérivation sur le faisceau structural $\mathcal O_X$ de $X$, c'est-à-dire un morphisme de faisceaux en groupes abéliens satisfaisant à la règle de Leibniz : 
$$ \delta_X(s.t) = \delta_X(s).t + s.\delta_X(t) \text{ pour tout ouvert } U \subset X \text{ et tous } s,t \in \mathcal O_X(U).$$

Si $(X,\delta_X)$  et $(Y,\delta_Y)$ sont deux $D$-schémas, on appelle \textit{morphisme de $D$-schémas de $(X,\delta_X)$ vers $(Y,\delta_Y)$}, tout morphisme de schémas $f = (|f|, f^\sharp) : (X,\mathcal O_X) \longrightarrow (Y,\mathcal O_Y)$  faisant  commuter le diagramme suivant : 

\begin{equation}\label{diagramme}
\begin{gathered}
 \xymatrix{
    \mathcal O_X \ar[r]^{\delta_X} & \mathcal O_X \\
    f^{-1} \mathcal O_Y \ar[u]^{f^{\sharp}} \ar[r]_{f^{-1} \delta_Y} & f^{-1} \mathcal O_Y \ar[u]_{f^{\sharp}}
  }
\end{gathered} \end{equation}}

On a ainsi défini une catégorie appelée \textit{catégorie des $D$-schémas} et notée $\textbf{D-Sch}$.

{\exam Soit $(A,\delta_A)$ un anneau (commutatif) différentiel. La dérivation $\delta_A$ se prolonge uniquement aux localisations de $A$ et induit une dérivation sur le faisceau structural $\mathcal O_{\mathrm{Spec}(A)}$. On obtient ainsi un $D$-schéma noté $(\mathrm{Spec}(A),\delta_A)$.}

{\rem\label{deuxiemedefinitionDschema} Plus généralement, on peut définir la catégorie $\textbf{C}$ des $D$-espaces localement annelés dont : 
\begin{itemize}
\item les objets sont les espaces localement annelés $(X,\mathcal O_X)$ muni d'une dérivation $\delta_X : \mathcal O_X \longrightarrow \mathcal O_X$ du faisceau structural $\mathcal O_X$.
\item les flèches sont les morphismes d'espaces localement annelés faisant commuter le diagramme (\ref{diagramme}).
\end{itemize}
La catégorie des $D$-schémas s'identifie alors à la sous-catégorie pleine de la catégorie $\textbf{C}$ dont les objets sont les $D$-espaces localement annelés $(X,\mathcal O_X,\delta_X)$ localement représentable sous la forme $(\mathrm{Spec}(A),\delta_A)$ où $(A,\delta)$ est un anneau différentiel.}

La catégorie $\textbf{D-Sch}$ admet un foncteur d'oubli naturel vers la catégorie des schémas obtenu en oubliant la dérivation sur le faisceau structural $ F_{oub} : \textbf{D-Sch} \longrightarrow \textbf{Sch}$.

{\nota Pour toute propriété $(P)$ des schémas (resp. des morphismes de schémas), nous dirons qu'un $D$-schéma (resp. un morphisme de $D$-schémas) possède la propriété $(P)$ si le schéma sous-jacent (resp. le morphisme de schémas sous-jacent) possède la propriété $(P)$.

Par exemple, on dira qu'un $D$-schéma $(X,\delta_X)$ est de type fini (resp. séparé, réduit, irréductible) si le schéma sous-jacent $X$ a la même propriété.}

{\Prop\label{produitdef} La catégorie $\textbf{D-Sch}$ admet des produits fibrés et un élément terminal. De plus, la formation des produits fibrés commute au foncteur d'oubli vers la catégorie des schémas.}

\begin{proof}
En effet, si $(C,\delta_C)$ est un anneau différentiel et $(A,\delta_A)$ et $(B,\delta_B)$ sont deux $(C,\delta_C)$-algèbres différentielles, on définit la $(C,\delta)$-algèbre différentielle
$$ (A,\delta_A) \otimes_{(C,\delta_C)} (B,\delta_B) = (A \otimes_C B , \delta_A \otimes_C \mathrm{Id}_B + \mathrm{Id}_A \otimes_C \delta_B).$$

On vérifie alors qu'on a ainsi défini le produit de $(\mathrm{Spec(A)},\delta_A)$ avec $(\mathrm{Spec(B)},\delta_B)$ au dessus de $(\mathrm{Spec(C)},\delta_C)$. 
Pour les $D$-schémas généraux, on procède par recollement à partir d'un recouvrement affine. Cette construction montre que la formation des produits commute au foncteur d'oubli vers la catégorie des schémas.
\end{proof}

{\rem Soit $(K,\delta)$ un corps différentiel. Le couple $(\mathrm{Spec}(K),\delta)$ est un $D$-schéma.  On note $\textbf{D-Sch}/(K,\delta)$, la catégorie des objets au dessus de $(\mathrm{Spec}(K),\delta)$.

La catégorie $\textbf{D-Sch}/(K,\delta)$ admet un élément terminal ainsi que des produits fibrés d'après la proposition \ref{produitdef}.
Cette proposition montre aussi que, si $(K,\delta) \subset (L,\delta_L)$ est une extension de corps différentiels alors, on a un foncteur de changement de base noté 
$$-\times_{(K,\delta)} (L,\delta_L) :
\begin{cases}
\textbf{D-Sch}/(K,\delta) &  \longrightarrow \textbf{D-Sch}/(L,\delta_L) \\
(X,\delta_X) & \mapsto (X,\delta_X)_{(L,\delta_L)} \end{cases}.$$}

Si $(X,\delta_X)$ est un $D$-schéma alors tout ouvert $U \subset X$ est naturellement muni d'une structure de $D$-schéma notée $(U,\delta_U)$ qui fait de l'immersion ouverte un morphisme de $D$-schémas.

{\lem\label{Zariskiouvert} Soient $(X,\delta_X)$ et $(Y,\delta_Y)$ deux $D$-schémas, $f : X \longrightarrow Y$ un morphisme de schémas et $U \subset X$ un ouvert schématiquement dense. Si $f_{|U} : (U,\delta_U) \longrightarrow (Y,\delta_Y)$ est un morphisme de $D$-schémas alors $f : (X,\delta_X) \longrightarrow (Y,\delta_Y)$ est un morphisme de $D$-schémas.}

La réciproque du lemme \ref{Zariskiouvert} est bien-sûr toujours vérifiée même lorsqu'on ne suppose plus que $U \subset X$ est schématiquement dense

\begin{proof}
On veut montrer que le diagramme $(\ref{diagramme})$ est commutatif. Comme l'inclusion $i: (U,\delta_U) \longrightarrow (X,\delta_X)$ est un morphisme de $D$-schémas, et par adjonction entre les foncteurs $i^{-1}$ et $i_\ast$, on a le diagramme commutatif:

\begin{equation}\label{diagrammeouvert}
\begin{gathered}
 \xymatrix{
    i_\ast \mathcal O_U \ar[r]^{i_\ast \delta_U} & i_\ast \mathcal O_U \\
     \mathcal O_X \ar[u] \ar[r]_{ \delta_X} &  \mathcal O_X \ar[u]
  } \end{gathered} \end{equation}
De plus, puisque $U \subset X$ est schématiquement dense, les deux flèches verticales sont injectives. Ainsi, pour vérifier la commutativité du diagramme (\ref{diagramme}), il suffit de vérifier que le diagramme suivant est commutatif :

$$
 \xymatrix{
    i_\ast \mathcal O_U \ar[r]^{i_\ast \delta_U} & i_\ast \mathcal O_U \\
     f^{-1} \mathcal O_Y \ar[u] \ar[r]_{ f^{-1} \delta_Y} &  f^{-1} \mathcal O_Y \ar[u]
  }$$

La commutativité du dernier diagramme est conséquence du fait que $f_{|U} : (U,\delta_U) \longrightarrow (Y,\delta_Y)$ est un morphisme de $D$-schémas et de l'adjonction entre les foncteurs $i^{-1}$ et $i_\ast$.
\end{proof}

{\defn Soient $(X,\delta_X)$ et $(Y,\delta_Y)$ deux $D$-schémas intègres et $f : X \dashrightarrow Y$ un morphisme rationnel. On dit que $f : (X,\delta_X) \dashrightarrow (Y,\delta_Y)$ est \textit{un morphisme rationnel de $D$-schémas} si sa restriction à un ouvert non vide de $X$ est un morphisme de $D$-schémas.} 

D'après le lemme \ref{Zariskiouvert}, si $f : (X,\delta_X) \dashrightarrow (Y,\delta_Y)$ est un morphisme rationnel de $D$-schémas alors $f$ définit un morphisme de $D$-schémas sur son ouvert de définition.

\subsubsection{Sous-schémas invariants}

{\defn Soit $(X,\delta_X)$ un $D$-schéma. On dit qu'un sous-schéma fermé $Y \subset X$ est un \textit{sous-schéma fermé invariant} de $(X,\delta_X)$, si le faisceau d'idéaux $\mathcal I_Y \subset \mathcal O_X$ définissant $Y$ est stable par la dérivation $\delta_X$, c'est-à-dire si pour tout ouvert $U \subset X$,
$$\mathcal I_Y(U) \subset (\mathcal O_X(U),\delta_X) \text{ est un idéal différentiel}.$$}

{\rem\label{miracle}  Soit $(X,\delta_X)$ un $D$-schéma et $i : Y \longrightarrow X$ un sous-schéma fermé invariant.
La dérivation $\delta_X$ passe au quotient en une dérivation $\delta_Y : \mathcal O_Y \longrightarrow \mathcal O_Y$ sur le faisceau structural de $Y$ qui fait de $(Y,\delta_Y)$ un $D$-schéma tel que l'immersion fermée $i : (Y,\delta_Y) \longrightarrow (X,\delta_X)$ est un morphisme de $D$-schémas.

Réciproquement, toute immersion fermée $i : (Y,\delta_Y) \longrightarrow (X,\delta_X)$ de $D$-schémas définit un sous-schéma fermé invariant de $(X,\delta_X)$.}

{\lem\label{changementdebase} Soient $f: (X,\delta_X) \longrightarrow (S,\delta_S)$ un morphisme de $D$-schémas et $(T,\delta_T) \longrightarrow (S,\delta_S)$ un changement de base.

Si $Y \subset X$ un sous-schéma fermé invariant de $(X,\delta_X)$ alors $Y \times_S T$ est un sous-schéma fermé invariant de $(X,\delta_X)\times_{(S,\delta_S)} (T,\delta_T)$.} 

\begin{proof}
L'immersion fermée $i : (Y,\delta_Y) \longrightarrow (X,\delta_X)$ induit après changement de base, un morphisme de $D$-schémas : 
$$ i_B : (Y,\delta_Y)\times_{(S,\delta_S)} (T,\delta_T) \longrightarrow (X,\delta_X)\times_{(S,\delta_S)} (T,\delta_T)$$
La notion d'immersion fermée étant invariante par changement de base, on conclut à l'aide de la remarque \ref{miracle}. 
\end{proof}

En particulier, si $(X,\delta_X)$ est un $D$-schéma, $U \subset X$ un ouvert et $Y \subset X$ un sous-schéma fermé invariant alors $Y\cap U$ est un sous-schéma fermé invariant de $(U,\delta_U)$. Réciproquement, on a le résultat suivant.

{\lem\label{generique} Soient $(X,\delta_X)$ un $D$-schéma et $Y \subset X$ un sous-schéma fermé. Considérons $U \subset X$ un ouvert tel que $Y \cap U$ est schématiquement dense dans $Y$.

Si $Y \cap U$ est un sous-schéma fermé invariant de $(U,\delta_U)$ alors $Y$ est un sous-schéma fermé invariant de $(X,\delta_X)$.}

\begin{proof}
On note $\mathcal I \subset \mathcal O_X$, le faisceau d'idéaux définissant $Y$. L'immersion ouverte $U \subset X$ induit un diagramme commutatif : 
$$
 \xymatrix{
   \mathcal O_X / \mathcal I \ar[r]^j & i_\ast (\mathcal O_X/I)_{|U} \\
   \mathcal O_X \ar[u]^\pi \ar[r] &  i_\ast \mathcal O_{X|U} \ar[u]_{\pi_U}
  }$$
où la première flèche horizontale $j$ est injective car $Y \cap U$ est schématiquement dense dans $Y$.

Montrons que le faisceau d'idéaux $\mathcal I \subset \mathcal O_X$ est stable par la dérivation $\delta_X$.

Soit $V \subset X$ un ouvert et $f \in I(V)$. On a alors $ \pi_U(\delta_X(f)_{|U}) =  \pi_U[\delta_U(f_{|U})] =  0$ car $Y \cap U$ est un sous-schéma fermé invariant de $U$. On en déduit que $(j \circ \pi) (\delta_X(f)) = 0$ donc que $\pi (\delta_X(f)) = 0$ car $j$ est injective, c'est-à-dire $\delta_X(f) \in I(V)$.
 \end{proof}

{\Prop\label{schemainvariantreduit} Soit $(X,\delta_X)$ un $D$-schéma noethérien au dessus de $(\mathbb{Q},0)$ et $Y \subset X$ un sous-schéma fermé invariant. Alors : 
\begin{itemize}
\item[(i)] Le sous-schéma fermé invariant réduit $Y_{red} \subset Y \subset (X,\delta_X)$ est un sous-schéma fermé invariant.
\item[(ii)] Les composantes irréductibles $Y_1, \cdots, Y_n \subset (X,\delta_X)$ de $Y_{red}$ sont des sous-schémas fermés invariants. 
\end{itemize}}

\begin{proof}
Si $(A,\delta)$ est un anneau différentiel de caractéristique $0$ alors le radical de tout idéal différentiel est un idéal différentiel \citep[Chapitre 2, Lemme 1.15]{MTF}. On en déduit que la propriété (i) est vérifiée.

Pour (ii), considérons $Y_i$ une composante irréductible de $Y_{red}$. Considérons un ouvert $U \subset X$ tel que
$$ Y_{red}  \cap U = Y_i \cap U \neq \emptyset$$ (Il suffit de considérer l'ouvert $ U = X \setminus \bigcup_{j \neq i} Y_j$).

Le schéma $Y_i \cap U = Y_{red} \cap U$ est un sous-schéma fermé invariant de $(U,\delta_U)$. Comme de plus $Y_i \cap U \subset Y_i$ est schématiquement dense dans $Y_i$ (car $Y_i$ est irréductible et $U \cap Y_i$ est non vide), le sous-schéma fermé $Y_i \subset (X,\delta_X)$ est invariant d'après le lemme \ref{generique}. 
\end{proof}

\subsection{Ensemble définissable associé à une $D$-variété} 
On fixe $(\U,\delta_U)$ un modèle saturé de la théorie des corps différentiellement clos. Tous les corps différentiels considérés seront des sous-corps différentiels de $(\U,\delta_U)$.

{\nota Soit $(K,\delta)$ un corps différentiel. On appelle \textit{$D$-variété au dessus de $(K,\delta)$}, tout $D$-schéma au dessus de $(K,\delta)$ séparé, de type fini et réduit.

En particulier, les $D$-variétés au dessus de $(K,\delta)$ ne sont pas supposées irréductibles et on parlera de $D$-variétés irréductibles (resp. absolument irréductibles) au dessus  de $(K,\delta)$  pour désigner les $D$-variétés au dessus de $(K,\delta)$ dont le $K$-schéma sous-jacent est irréductible (resp. absolument irréductible).

{\rem La notion de $D$-variété est préservée par changement de base par un corps différentiel puisque les corps différentiels sont de caractéristique $0$, c'est-à-dire si $(K,\delta) \subset (L,\delta_L)$ est une extension de corps différentiels et $(X,\delta_X)$ est une $D$-variété au dessus de $(K,\delta)$, alors $(X,\delta_X)_{(L,\delta_L)}$ est une $D$-variété au dessus de $(L,\delta_L)$.}

\subsubsection{Interprétation d'une $D$-variété dans $\textbf{DCF}_0$} Soit $(K,\delta)$ un corps différentiel.

{\nota Soient $(K,\delta) \subset (L,\delta_L)$ une extension de corps différentiels et $(X,\delta_X)$ un $D$-schéma au dessus de $(K,\delta)$. On appelle \textit{ensemble des  $(L,\delta_L)$-points différentiels de ($X,\delta_X)$}, l'ensemble 
$$(X,\delta_X)^{(L,\delta_L)} = \mathrm{Hom}_{D-Sch/(K,\delta)}[(\Spec L, \delta_L) ;(X,\delta_X)].$$}

On montre dans la suite de cette partie que l'ensemble $(X,\delta_X)^{(U,\delta_U)}$ peut être muni d'une structure d'ensemble définissable dans $\textbf{DCF}_0$.

{\exam Soit $(S)$ un système d'équations différentielles algébriques à paramètres dans $(K,\delta)$, c'est-à-dire un système d'équations différentielles algébriques de la forme :
$$
(S) : \begin{cases} P_1(x_1,\cdots,x_r,\delta(x_1), & \cdots , \delta^k(x_1),\cdots \delta^k(x_r)) = 0 \\ & \vdots \\ P_r(x_1,\cdots,x_r,\delta(x_1),& \cdots , \delta^k(x_1),  \cdots \delta^k(x_r)) = 0 \end{cases}$$
où les $P_i \in K[X_1^{(0)},\cdots, X_r^{(0)}, \cdots , X_1^{(k)}, \cdots, X_r^{(k)}]$ sont des polynômes.

Le foncteur $\mathcal S$ "solutions de $(S)$" de la catégorie des $(K,\delta)$-algèbres différentielles vers la catégorie des ensembles est représentable par un $D$-schéma affine $(X,\delta_X)$ au dessus de $(K,\delta)$.

En effet, considérons $K\lbrace X_1 , \cdots, X_r \rbrace$ la $(K,\delta)$-algèbre différentielle libre engendrée par les indéterminées $X_1,\cdots , X_r$. Pour tout $i \leq r$, on a  
$$P_i(X_1,\cdots,X_r,\delta(X_1),  \cdots , \delta^k(X_1),\cdots \delta^k(X_r)) \in K\lbrace X_1, \cdots , X_r \rbrace.$$

Notons $I \subset K\lbrace X_1, \cdots X_r \rbrace$ l'idéal différentiel engendré par ces éléments. Le foncteur $\mathcal S$ est alors représentable par le $D$-schéma affine associé à la $(K,\delta)$-algèbre différentielle  $K\lbrace X_1,\cdots, X_n \rbrace / I$. En général, ce $D$-schéma n'est ni réduit ni de type fini.}

{\rem On en déduit que si $(X,\delta_X)$ représente le foncteur solutions d'un système d'équations différentielles algébriques $(S)$ à paramètres dans $(K,\delta)$ alors
$(X,\delta_X)^{(\U,\delta_\U)} = \mathcal S(\U,\delta_U)$ s'identifie à un ensemble $K$-définissable de $(\U,\delta_\U)$.}

{\cons Soit $(X,\delta_X)$ une $D$-variété au dessus de $(K,\delta)$. Considérons $\mathcal V = (V_i,f_i)_{i = 1 ,\cdots , n}$ où $(V_i)$ est un recouvrement ouvert affine de $X$ et où les $f_i : V_i \longrightarrow \mathbb{A}^n$ sont des immersions fermées.

On peut construire une formule $\phi_\mathcal V(x)$ sans quantificateurs à paramètres dans $K$ telle que pour toute extension de corps différentiels $(K,\delta) \subset (L,\delta_L)$, on a l'identification suivante :  
$$(X,\delta_X)^{(L,\delta_L)}  = \phi_\mathcal V(L,\delta_L).$$ }

\begin{proof}
Pour toute extension de corps différentiels $(K,\delta) \subset (L,\delta_L)$, puisque $(V_i)_{i = 1}^n$ est un recouvrement ouvert de $X$, on a une application surjective : 
$$ \pi : \bigsqcup_{i =1}^n (V_i,\delta_{V_i})^{(L,\delta_L)} \longrightarrow (X,\delta_X)^{(L,\delta_L)}.$$
qui identifie $(X,\delta_X)^{(L,\delta_L)}$ au quotient de $\bigsqcup_{i =1}^n (V_i,\delta_{V_i})^{(L,\delta_L)}$ par une relation d'équivalence $E$.

Pour tout $i \leq n$, les immersions fermées $f_i : V_i \longrightarrow \mathbb{A}^n$ permettent d'identifier $(V_i,\delta_i)^{(L,\delta_L)}$ à l'ensemble des solutions d'un système d'équations différentielles algébriques $(S_i)$ à $n$ indéterminées et donc à un ensemble $K$-définissable sans quantificateurs dans $(L,\delta_L)$.

On en déduit que $\bigsqcup_{i =1}^n (V_i,\delta_{V_i})^{(L,\delta_L)}$ s'identifie à un ensemble $K$-définissable de $(L,\delta_L)$ sans quantificateurs.

{\asser La relation d'équivalence $E$ est $K$-définissable sans quantificateurs dans le langage des anneaux.}

En effet, sa définition est donnée par les fonctions de transitions du recouvrement ouvert $(U_i)_{i \leq n}$ et donc définissable sans quantificateurs dans le langage des anneaux (et donc dans le langage des anneaux différentiels).
  
Par élimination des imaginaires dans la théorie des corps algébriquement clos, il existe une formule sans quantificateurs $\phi_\mathcal V(\overline{x})$ telle que pour toute extension de corps différentiels $(K,\delta) \subset (L,\delta_L)$: 
$$ (X,\delta_X)^{(L,\delta_L)} = \bigsqcup_{i =1}^n (V_i,\delta_{V_i})^{(L,\delta_L)}/E = \phi_{\mathcal V}(L,\delta_L).$$
\end{proof}

Le lemme suivant montre que la structure définissable induite ne dépend pas du recouvrement ouvert choisi.

{\lem Soit $\mathcal E$ l'ensemble des données $\mathcal V = (V_i,f_i)_{i = 1 ,\cdots , n}$ où $(V_i)$ est un recouvrement ouvert affine de $V$ et où les $f_i : V_i \longrightarrow \mathbb{A}^n$ sont des immersions fermées.

Si $\mathcal V, \mathcal V' \in \mathcal E$ alors pour toute extension de corps différentiels $(K,\delta) \subset (L,\delta_L)$, les ensembles $\phi_\mathcal V(L,\delta_L)$ et $\phi_{\mathcal V'}(L,\delta_L)$ sont en bijection définissable sans quantificateurs à paramètres dans $(K,\delta)$.}

\begin{proof}
Deux recouvrement ouverts admettent toujours un raffinement commun. On peut donc supposer que $\mathcal V'$ est un recouvrement plus fin que $\mathcal V$.
Il suffit alors de vérifier que si $X$ est un $D$-schéma affine et $(V_i)_{i \leq n}$ un recouvrement affine de $X$ alors l'application 
$$ \pi : \bigsqcup_{i =1}^n (V_i,\delta_{V_i})^{(L,\delta_L)} \longrightarrow (X,\delta_X)^{(L,\delta_L)}$$ est $K$-définissable sans quantificateurs, ce qui est immédiat.
\end{proof}
Suivant l'usage en théorie des modèles, si $(K,\delta) \subset (\U,\delta_\U)$ est une extension différentiellement close, on identifiera $(X,\delta_X)^{(\U,\delta_\U)}$ à un ensemble $K$-définissable dans $(\U,\delta_\U)$ sans préciser le recouvrement affine choisi.
On a la description intrinsèque suivante de la structure induite par $(\U,\delta_U)$ sur $(X,\delta_X)^{(\U,\delta_\U)}$.
 
{\cor\label{description des formules} Soient $(K,\delta)$ un corps différentiel et $(X,\delta_X)$ une $D$-variété au-dessus de $(K,\delta)$.
Les sous-ensembles $K$-définissables  de $(X,\delta_X)^{(\U,\delta_\U)}$ sont les combinaisons booléennes d'ensembles de la forme $(Y,\delta_Y)^{(\U,\delta_\U)}$ où $Y \subset (X,\delta_X)$ est une sous-variété invariante.}

\begin{proof}
Supposons que le $D$-schéma $(X,\delta_X)$ est affine et considérons une immersion fermée $X \subset \mathbb{A}^n$. Par élimination des quantificateurs, les sous-ensembles $K$-définissables de $\U^n$ sont les combinaisons booléennes de sous-ensembles de la forme :
$$ V(I) = \lbrace x \in \U^n \text{ | } f(x) = 0 \text{ , } \forall f \in I \rbrace \subset \U^n.$$

où $I \subset K\lbrace X_1,\cdots,X_n \rbrace$ est un idéal différentiel de la $(K,\delta)$-algèbre libre d'indéterminées $X_1,\cdots X_n$.
On en déduit que les sous-ensemble définissables de $(X,\delta_X)^{(\U,\delta_\U)}$ sont les combinaisons booléennes des ensembles $V(I)$ où $I \subset K\lbrace X_1,\cdots,X_n \rbrace$ est un idéal différentiel contenant l'idéal $I_X$ définissant $X$.
Ces derniers sont en correspondance avec les sous-schémas fermés invariants de $(X,\delta_X)$. La proposition \ref{schemainvariantreduit} permet alors de conclure.

On en déduit le cas général à l'aide du lemme \ref{generique} et de la construction précédente. 

\end{proof}

{\nota Soient $(X,\delta_X)$ une $D$-variété au dessus d'un corps différentiel $(K,\delta)$ et $\Sigma = (X,\delta_X)^{(\U,\delta_\U)}$. On note $S_\Sigma(K)$ l'ensemble des types vivants sur $\Sigma$ (c'est-à-dire vérifiant $p_\Sigma(x) \models x \in \Sigma$) à paramètres dans $K$.}

{\cor\label{description des types} Soient $(K,\delta)$ un corps différentiel et $(X,\delta_X)$ une $D$-variété au dessus de $(K,\delta)$. On pose $\Sigma = (X,\delta_X)^{(\U,\delta_\U)}$. Pour toute sous-variété fermée invariante $Y \subset (X,\delta_X)$, il existe un unique type complet $p_{(Y,\delta_Y)}(x)$ vivant sur $\Sigma$  vérifiant : 
$$ p_{(Y,\delta_Y)}(x) \models \begin{cases} x \in (Y,\delta_Y)^{(\U,\delta_\U)} \\ x \notin (Y',\delta_Y')^{(\U,\delta_\U)} \text{ pour } Y' \nsubseteq Y \text{ une sous-variété invariante stricte} \end{cases}$$

De plus, en notant  $\mathrm{Inv}_\delta (X,\delta_X)$ l'ensemble des sous-variétés fermées irréductibles et invariantes de  $(X,\delta_X)$, l'application $\mathrm{Inv}_\delta (X,\delta_X)  \longrightarrow S_\Sigma(K)$ ainsi définie est une bijection.}

Le cas affine est donné par le corollaire  \ref{descriptiondestypesaffine} et la cas général s'en déduit aisément à l'aide du lemme \ref{generique}.

\subsubsection{Type générique d'une $D$-variété irréductible} Soit $(X,\delta_X)$ une $D$-variété irréductible au dessus de $(K,\delta)$. Le corollaire précédent montre qu'il existe un type $p \in S(K)$ à paramètres dans $K$ et vivant sur $(X,\delta_X)^{(\U,\delta_\U)}$ dont les réalisations coïncident avec les réalisations du point générique (au sens de la géométrie algébrique) de $(X,\delta_X)$.

{\defn\label{type generique def} Soient $(K,\delta)$ un corps différentiel  et $(X,\delta_X)$ une $D$-variété irréductible au dessus de $(K,\delta)$. On appelle \textit{type générique de $(X,\delta_X)$}, le type $p_{(X,\delta_X)} \in S(K)$ correspondant à la variété $X$ elle-même dans le corollaire \ref{description des types}.}

{\lem\label{typegenerique} Soient $(K,\delta)$ un corps différentiel  et $(X,\delta_X)$ une $D$-variété irréductible au dessus de $(K,\delta)$. Posons $\Sigma = (X,\delta_X)^{(\U,\delta_\U)}$.
\begin{itemize}
\item[(i)] Le type $p_{(X,\delta_X)}$ est l'unique type d'ordre maximal vivant sur $\Sigma$. De plus, il vérifie
$$\mathrm{ord}(p_{(X,\delta_X)}) = \mathrm{dim}(X).$$

\item[(ii)] Le type $p_{(X,\delta_X)}$ est stationnaire si et seulement si la $D$-variété $(X,\delta_X)$ est absolument irréductible. Dans ce cas, son unique  extension non-déviante à une extension $(K,\delta) \subset (L,\delta_L)$ est le type générique de $(X,\delta_X)_{(L,\delta_L)}$.
\end{itemize}}

\begin{proof}
Soit $a \models p_{(X,\delta_X)}$ une réalisation du type $p_{(X,\delta_X)}$. Par définition du type $p_{(X,\delta_X)}$, ses réalisations sont les réalisations du point générique de $(X,\delta_X)$ et donc le sous-corps différentiel  $K\langle a \rangle$ engendré par $a$ dans $(\U,\delta_U)$ est isomorphe au corps différentiel $(K(X),\delta_X)$.
On en déduit que :
$$ \mathrm{ord}(p) = \mathrm{tr}(K\langle a \rangle/K) = \mathrm{tr}(K(X)/K) = \mathrm{dim}(X).$$

La propriété (i) est donc vérifiée. Montrons la propriété (ii). Soit $(K,\delta) \subset (L,\delta_L)$ une extension de corps différentiels.

Le corollaire \ref{description des types} montre que les extensions de $p_{(X,\delta_X)}$ à $(L,\delta_L)$ correspondent aux sous-variétés irréductibles invariantes de $(X,\delta_X)_{(L,\delta_L)}$ dont la projection vers $X$ est dominante.

D'après le corollaire \ref{ordredeviation} et la propriété $(i)$ au-dessus, les extensions non déviantes de $p_{(X,\delta_X)}$ à $(L,\delta_L)$ correspondent aux sous-variétés irréductibles invariantes de $(X,\delta_X)_{(L,\delta_L)}$ de dimension $\mathrm{dim}(X)$ et dont la projection vers $X$ est dominante.

On en déduit que les extension non-déviantes de $p_{(X,\delta_X)}$ à $(L,\delta_L)$ correspondent aux composantes irréductibles de $(X,\delta_X)_{(L,\delta_L)}$.
\end{proof}

{\rem La propriété $(i)$ du lemme précédent montre qu'il est facile de calculer l'ordre du type générique $p_{(X,\delta_X)}$ si l'on connait la $D$-variété $(X,\delta_X)$. En particulier l'ordre de $p_{(X,\delta_X)}$ ne dépend que de la variété $X$ et non de la structure de $D$-variété sur $X$ considérée.

Cependant, on ne sait en général pas déterminer $\mathrm{RM}(p_{(X,\delta_X)})$ et $\mathrm{RU}(p_{(X,\delta_X)})$. De plus, ces deux quantités dépendent de la structure de $D$-variété sur $X$.}

{\lem\label{morphismerationnel} Soient $(K,\delta)$ un corps différentiel, $(X,\delta_X)$ et $(Y,\delta_Y)$ des  $D$-variétés irréductibles au dessus de $(K,\delta)$. Considérons $a \models p_{(X,\delta_X)}$, $b \models p_{(Y,\delta_Y)}$ des réalisations des types génériques de $(X,\delta_X)$ et $(Y,\delta_Y)$ respectivement. Les propriétés suivantes sont équivalentes : 
\begin{itemize}
\item[(i)] $b \in \mathrm{dcl}(K,a)$
\item[(ii)] Il existe un morphisme rationnel de $D$-variétés
$$ \phi : (X,\delta_X) \dashrightarrow (Y,\delta_Y) \text{ tel que } \phi^\U(a) = b.$$
où $\phi^\U$ désigne le morphisme induit  par $\phi$ après passage au points différentiels dans $(\U,\delta_\U)$.
\end{itemize}}

\begin{proof}
$(ii) \Longrightarrow (i)$ est immédiat puisque l'application $\phi^\U$ est $K$-définissable et $\phi^\U(a) = b$.
Montrons que $(i) \longrightarrow (ii)$. Comme $b \in \mathrm{dcl}(K,a)$, on en déduit que : 
$$ K\langle b \rangle = \mathrm {dcl}(K,b) \subset K \langle a \rangle = \mathrm{ dcl}(K,a) \subset (\U,\delta_U).$$
Les isomorphismes de corps différentiels au dessus de $(K,\delta)$, $(K(X),\delta_X) \simeq K\langle a \rangle  $ et $(K(Y),\delta_Y) \simeq  K\langle b \rangle$ définissent alors une injection de corps différentiels $i : (K(Y),\delta_Y) \longrightarrow (K(X),\delta_X)$ et donc un morphisme rationnel de $D$-variétés  
$ \phi : (X,\delta_X) \dashrightarrow (Y,\delta_Y) $
tel que $\phi^\U(a) = b$.
\end{proof}
\bibliographystyle{alpha}
\bibliography{bibliographie}

\end{document}